\documentclass[10pt]{article}

\RequirePackage{amsmath}
\RequirePackage{amsfonts}
\RequirePackage{amssymb}
\RequirePackage{hyperref}
\RequirePackage{multirow}
\RequirePackage{algorithm}
\RequirePackage{algpseudocode}
\RequirePackage{graphicx}
\RequirePackage{xcolor}
\RequirePackage{amsthm}
\RequirePackage{amsmath}
\RequirePackage{amsfonts}
\RequirePackage{amssymb}
\RequirePackage{hyperref}
\RequirePackage{multirow}
\RequirePackage{algorithm}
\RequirePackage{algpseudocode}
\RequirePackage{graphicx}
\RequirePackage{xcolor}
\usepackage{enumitem}
\usepackage{subcaption}
\allowdisplaybreaks

\usepackage{dpr_fec,ifthen,dsfont}
\newcommand{\ones}{\mathds{1}}

\newtheorem{remark}{Remark}[section]

\newcommand{\edit}[1]{#1}
\newcommand{\editcolor}{}

\usepackage{isetechreport}
\def\reportyear{25T}
\def\reportno{015}
\def\revisionno{1}
\def\originaldate{August 30, 2025}
\def\revisiondate{\today}
\coraltrue
\cvcrfalse

\begin{document}

\title{Active-Set Identification in Noisy and Stochastic Optimization}

\author{Frank E.~Curtis\thanks{E-mail: frank.e.curtis@lehigh.edu}}
\author{Daniel P.~Robinson\thanks{E-mail: daniel.p.robinson@lehigh.edu}}
\author{Lara Zebiane\thanks{E-mail: lara.zebiane@lehigh.edu}}
\affil{Department of Industrial and Systems Engineering, Lehigh University}
\titlepage

\maketitle

\begin{abstract}
  Identifying active constraints from a point near a local minimizer is important both theoretically and practically in constrained continuous optimization as it can help identify optimal Lagrange multipliers and essentially reduce an inequality-constrained problem to an equality-constrained one. Traditional active-set identification guarantees have been proved under various types of constraint qualifications. Importantly, they also assume exact function and derivative values. This work extends such results to settings when both objective and constraint function and derivative values have deterministic or stochastic noise. \edit{Two strategies are analyzed. For each strategy, it is shown that, under constraint qualifications and other standard conditions, the active set of a local minimizer can be identified correctly when a point is close enough to the minimizer and the noise is sufficiently small. It is also shown that regularity of a minimizer is critical. Under strict complementarity, a fixed noise level suffices, whereas in the presence of a weakly active constraint the noise must vanish at least linearly in the distance to the minimizer in general, a rate that is also shown to be sufficient for correct identification.  Guarantees are also stated for the use of active-set identification strategies within a stochastic algorithm, with explicit sample-size conditions. The theoretical findings are demonstrated} with two simple illustrative examples and a more realistic constrained neural-network training task.
\end{abstract}

\newcommand{\eapprox}{\etilde}
\newcommand{\capprox}{\ctilde}
\newcommand{\fapprox}{\ftilde}
\newcommand{\egradapprox}{\widetilde{\nabla e}}
\newcommand{\cgradapprox}{\widetilde{\nabla c}}
\newcommand{\fgradapprox}{\widetilde{\nabla f}}
\newcommand{\dapprox}{\dtilde}
\newcommand{\lambdaapprox}{\tilde{\lambda}}
\newcommand{\muapprox}{\tilde{\mu}}
\newcommand{\laracomment}[1]{\textcolor{purple}{\Lara{}: #1}}
\newcommand{\frankcomment}[1]{\textcolor{blue}{Frank: #1}}
\newcommand{\epsilonerror}{\epsilon_{\text{err}}}
\newcommand{\PP}{\Pcal\!\Pcal}
\newcommand{\PN}{\Pcal\!\Ncal}
\newcommand{\NP}{\Ncal\!\Pcal}
\newcommand{\NN}{\Ncal\!\Ncal}

\section{Introduction}\label{sec.introduction}

A critical feature of an algorithm for solving constrained continuous optimization problems is its ability to identify the inequality constraints that are active at a local minimizer.  Suppose that, with respect to a function $c : \R{n} \to \R{q}$, a problem defined with respect to a decision variable $x \in \R{n}$ involves inequality constraints of the form $c(x) \leq 0$.  The aim of active-set identification is for an algorithm to be able to determine, with respect to a local minimizer \edit{and/or KKT point} $x_* \in \R{n}$, the set $\{i \in \{1,\dots,q\} : c_i(x_*) = 0\}$ of inequality constraints that hold at equality---i.e., are \emph{active}---at $x_*$, even if the algorithm only has access to a solution estimate that is within a local neighborhood of $x_*$.  Active-set identification is critical in multiple ways.  One way is in the context of algorithms that, in each iteration, compute a search direction by (a) estimating the active set, and then (b) computing the search direction with respect to the estimated active set (while ignoring all of the constraints that are estimated to be \emph{inactive}) \cite{ByrdGoulNoceWalt2003,ByrdGoulNoceWalt2005}.  Another manner in which active-set identification is critical is in the computation of Lagrange multipliers, which in turn can be useful algorithmically and/or for sensitivity analysis of the \edit{solution of interest}.

Various strong active-set identification guarantees have been proved for settings in which a problem's objective and constraint functions are continuously differentiable and an algorithm has access to exact function and derivative values; see, e.g., \cite{Burk1990,BurkMore1988,FaccFiscKanz1998,HareLewi2004,Lewi2002,OberWrig2006,Wrig1993}.  Typically, these guarantees require that a constraint qualification---such as the linear independence constraint qualification (LICQ) or Mangasarian-Fromovitz constraint qualification (MFCQ) \cite{MangFrom1967}---\edit{and other standard conditions hold at the local minimizer or Karush-Kuhn-Tucker (KKT) point $x_*$, such as strict complementarity or a second-order sufficiency condition}.  Each such guarantee states that, if a point~$x$ is sufficiently close to $x_*$, then the active set at $x_*$ can be identified correctly by solving a subproblem involving problem function and derivative values at~$x$.

The contributions of this paper are to show that such guarantees can be extended to settings in which the problem function and derivative values contain errors.  The design and analysis of optimization algorithms when function and derivative values are subject to errors are active areas of research, which motivates such extensions.  Errors can arise in multiple settings.  One setting of modern interest is that of \emph{deterministic noise}, which we use to refer to situations in which the evaluation of a problem function or derivative always results in the same approximate value whenever the evaluation occurs multiple times at the same point.  This may represent situations in which the function is defined by a simulation that is subject to computational noise.  We refer the reader to the recent articles \cite{berahas2019derivative,curtis2025interior,dezfulian2024convergence,oztoprak2023constrained,sun2023trust,xie2020analysis} on the design and analysis of algorithms for solving unconstrained or constrained optimization problems subject to deterministic noise.  Another setting of modern interest is that of \emph{stochastic noise}.  This refers to situations in which a problem function value is defined for any point in its domain by, e.g., a large finite sum and/or an expectation over a random variable argument.  The noise arises from the fact that, in practice, it is intractable to evaluate the value or derivative of such a function exactly, and instead one obtains an approximate value defined by a random sample of inputs.  There are a large number of articles on the design and analysis of algorithms for solving unconstrained optimization problems when evaluations are subject to stochastic noise; see, e.g., the references in \cite{BottCurtNoce2018,CurtSche17}.  For recent articles on solving constrained problems of this type, see \cite{BeraCurtONeiRobi2023,BeraCurtRobiZhou2021,CurtJianWang2024,CurtONeiRobi2024,CurtRobiZhou2024b,CurtRobiZhou2024,Fang2024Fully,Na2022adaptive,na2023inequality}.

Despite these efforts on algorithm design and analysis, little has been done to extend active-set identification guarantees to settings with errors.  (An exception is~\cite{na2023inequality} wherein an active-set algorithm is devised for solving problems when objective function and derivative values are subject to stochastic noise while constraint function and derivative values are exact.  However, in that setting, the active-set estimation is tied to a specific algorithm and the results do not allow noise in the constraint function or derivative values.)  In this paper, we consider strategies for active-set identification that may be employed irrespective of any particular algorithm, and when function and derivative values may be subject to noise.  Following primarily~\cite{OberWrig2006}, we show that when a point~$x$ is sufficiently close to a local minimizer~$x_*$ and the noise is sufficiently small, the active set at~$x_*$ can be identified correctly through solving a few types of subproblems at~$x$.  In this respect, it is particularly comforting to confirm that state-of-the-art active-set identification techniques are robust enough to extend readily to settings involving function and derivative noise.  We also state how our active-set identification guarantees can be employed for accurate active-set identification in the context of a stochastic algorithm for constrained optimization.  Finally, we demonstrate our theoretical guarantees through a pair of illustrative numerical experiments, one involving a couple of small-dimensional examples where the results can be seen visually and another involving a constrained neural-network training problem.

\edit{Our contributions are best understood within the context of perturbation theory of nonlinear optimization; see, e.g., \cite{BonnShap2000,KlatKumm2002,Fiac1976,DontRock2009,IzmaSolo2014,Robi1974,Robi1980,Koji1980}. That theory describes how the solutions and optimal Lagrange multipliers of a constrained optimization problem may change with perturbations of the problem data. Broadly speaking, perturbation theory has shown that, under regularity conditions, the active set of a perturbed solution coincides with that of the unperturbed one. Active-set identification in the presence of noise asks a related, but distinct question. Since a point $x$ at which one may attempt to perform active-set identification might not be the solution of any perturbation of the original optimization problem of interest and the errors in the function and derivative values at $x$ do not correspond to a perturbation of the problem in the classical sense, the errors that may be observed at two different points can in general be unrelated, so there might be no perturbed problem whose solution one could follow as $x$ approaches a minimizer~$x_*$. Rather, what one can do is consider a subproblem at a given point $x$ with noisy data as a perturbation of the same type of subproblem at~$x_*$ with exact data, whose solution is known, and ask when the active-set estimate produced by the former equals the active-set at $x_*$. This is the approach taken in some of the active-set identification literature when employing exact data, and it is the approach that we take as well.}

\edit{There are two features of our theoretical results for noisy settings that deserve emphasis, since they would not be evident from a direct appeal to classical perturbation theory or prior results on active-set identification with exact values. First, our results quantify the level of noise that is admissible for correct identification. In particular, for one of our proposed strategies that relies on strict complementarity, we state an admissible noise level explicitly in terms of the conditioning of the KKT system at $x_*$ (see Lemma~\ref{lem.explicit}). On the other hand, for another proposed approach, we show a rate at which it suffices for the noise to shrink as $x \to x_*$, a result that is unique to our noisy setting (see Corollary~\ref{cor.rate}).  Second, we show that the dependence on regularity is critical.  In particular, we show that a fixed noise level suffices under strict complementarity, whereas in the presence of a weakly active constraint the noise must, in general, vanish at least linearly in the distance $\|x - x_*\|_2$, which is the rate that we show suffices (see Theorem~\ref{th.impossible} and Corollary~\ref{cor.rate}).  Overall, our analysis shows that the two strategies behave differently under noise, where the primary reason for the difference is not the strategies themselves, but rather possible degeneracy at $x_*$.}

\subsection{Notation}

We use $\R{}$ to denote the set of real numbers and use $\R{}_{\geq r}$ to denote the set of real numbers greater than or equal to $r \in \R{}$.  Real-valued $n$-vectors and $m$-by-$n$ matrices are denoted by $\R{n}$ and $\R{m \times n}$, respectively.  Given an integer $q \geq 1$, we define the corresponding set of integers $[q] := \{1,\dots,q\}$.  Given $\Acal \subseteq [q]$, we use $\Acal^c \subseteq [q]$ to denote the complement of $\Acal$ with respect to the set $[q]$, namely, $\Acal^c := [q] \setminus \Acal$. We use $\ones$ to denote the vector with all components equal to~1; in each instance, the length of the vector is determined by the context in which it appears.  Similarly, we use $I$ to denote the identity matrix whose size is determined by its context.  Given any vector (function) $v$, its $i$th component is written as~$v_i$.  In addition, given such $v$, we define $[v]_+ := \max\{v,0\}$ and $[v]_- := \max\{-v,0\}$, where in each case the maximum is taken component-wise.  Given a pair of vectors $(a,b) \in \R{n} \times \R{n}$, we use $a \perp b$ to indicate that $a_i b_i = 0$ for all $i \in [n]$.  Given sets $\Acal \subseteq \R{n}$ and $\Bcal \subseteq \R{n}$, we denote the $\ell_2$-norm-induced distance between them as
\bequation \label{eq.def_dist}
  \dist(\Acal,\Bcal) := \inf_{(a,b) \in \Acal \times \Bcal}\ \|a - b\|_2.
\eequation

\subsection{Organization}

In~\S\ref{sec.preliminary}, we state the form of the optimization problems that we consider throughout the paper, and provide background on optimality conditions and related concepts that are referenced.  In~\S\ref{sec.exact}, we discuss active-set identification guarantees in the context of having exact problem function and derivative values.  This section includes one new contribution in the form of a modification of a result from \cite{OberWrig2006} from the setting of employing a linear optimization problem (LP) with a trust-region constraint to one employing a quadratic optimization problem (QP).  \edit{Our main contributions are contained in~\S\ref{sec.noisy}.  In \S\ref{sec.lagrange} and~\S\ref{sec.primal_step}, we prove that a Lagrange multiplier-based strategy and a primal-step-based strategy, respectively, identify the active set correctly with explicit tolerances on the allowable noise.  In \S\ref{sec.sharp}, we show that the different noise requirements of the two strategies reflect the regularity, or lack thereof, of the optimization problem at the local minimizer of interest.  In particular, we show that in the presence of a weakly active constraint at a minimizer, no rule can identify the active set under a fixed noise level, and that the linear rate at which our analysis requires the noise to vanish cannot be improved in general.  In \S\ref{sec.stochastic_algorithm}, we demonstrate the use of our results for a stochastic algorithm, with explicit sample-size conditions.}  We provide the results of illustrative numerical experiments in~\S\ref{sec.numerical} and provide concluding remarks in~\S\ref{sec.conclusion}.

\section{Preliminaries}\label{sec.preliminary}

Given an objective function $f : \R{n} \to \R{}$, equality-constraint function $e : \R{n} \to \R{p}$, and inequality-constraint function $c : \R{n} \to \R{q}$, where all functions (i.e., $f$, $e$, and $c$) are continuously differentiable, consider the constrained continuous optimization problem
\bequation\label{eq.opt}
  \min_{x \in \R{n}}\ f(x)\ \st\ e(x) = 0\ \text{and}\ c(x) \leq 0.
\eequation
Let the Lagrangian $\Lcal : \R{n} \times \R{p} \times \R{q} \to \R{}$ of \eqref{eq.opt} be defined by $\Lcal(x,y,z) = f(x) + e(x)^Ty + c(x)^Tz$, where $y \in \R{p}$ and $z \in \R{q}$ are Lagrange multiplier variables.  Under a constraint qualification, first-order necessary conditions for \eqref{eq.opt} are that at a local minimizer $x \in \R{n}$ there exist $y \in \R{p}$ and $z \in \R{q}$ such that
\bequation\label{eq.kkt}
  \nabla f(x) + \nabla e(x)y + \nabla c(x)z = 0,\ \ e(x) = 0,\ \ \text{and}\ \ 0 \leq z \perp c(x) \leq 0,
\eequation
where the first two equations correspond to the conditions $\nabla_x \Lcal(x,y,z) = 0$ and $\nabla_y \Lcal(x,y,z) = 0$.  We refer to any tuple of vectors $(x,y,z) \in \R{n} \times \R{p} \times \R{q}$ satisfying \eqref{eq.kkt} as a KKT point.  We refer to $x \in \R{n}$ as a primal KKT point if and only if there exists $(y,z) \in \R{p} \times \R{q}$ such that $(x,y,z)$ is a KKT point.

At a local minimizer $x \in \R{n}$ of \eqref{eq.opt}, a pair $(y,z)$ such that $(x,y,z)$ satisfies \eqref{eq.kkt} is not necessarily unique.  Therefore, for any $x \in \R{n}$, let us define $\Dcal(x) := \{(y,z) \in \R{p} \times \R{q} : (x,y,z)\ \text{satisfies}\ \eqref{eq.kkt}\}$. That is, for any $x \in \R{n}$, the set $\Dcal(x)$ (which may be empty) contains all Lagrange multiplier pairs that yield a KKT point along with~$x$.  \edit{At a local minimizer~$x$ of \eqref{eq.opt} the set $\Dcal(x)$ is nonempty and bounded if and only if the MFCQ holds at~$x$ \cite{Gauv1977}, and if the LICQ holds at a local minimizer~$x$ then $\Dcal(x)$ is a singleton.}  Strict complementarity holds for a KKT point $(x,y,z)$ if and only if $z_i > 0$ for all $i \in [q]$ such that $c_i(x) = 0$.

Let us indicate the set of active constraints at any primal point $x \in \R{n}$ as the set of indices of the inequality constraints that are satisfied at equality at $x$, namely, $\Acal(x) := \{i \in [q] : c_i(x) = 0\}$. Let $(x,y,z)$ be a KKT point.  Of the constraints whose indices are contained in~$\Acal(x)$, the set of \emph{weakly} active constraints at $x$ are those such that the corresponding Lagrange multiplier is equal to zero for all elements of $\Dcal(x)$, i.e., $\Acal_0(x) := \{i \in \Acal(x) : z_i = 0\ \text{for all}\ (y,z) \in \Dcal(x)\}$. Later, we reference a particular second-order optimality condition.  In particular, let us define the critical cone at a local minimizer $x \in \R{n}$ as
\bequationNN
  \baligned
    \Ccal(x) := \{ d \in \R{n} :\ \nabla e_i(x)^T d = 0\ \text{for all}\ i \in [p],\ \ \nabla c_i(x)^T d &\leq 0\ \text{for all}\ i \in \Acal_0(x),\ \text{and} \\
    \nabla c_i(x)^T d &= 0\ \text{for all}\ i \in \Acal(x) \setminus \Acal_0(x)\}.
  \ealigned
\eequationNN
Then, the second-order condition that we reference requires that, at $x$, one has
\bequation\label{eq:second_order_cond}
  d^T \nabla_{xx}^2 \Lcal(x, y, z) d > 0\ \ \text{for all}\ \ d \in \Ccal(x) \setminus \{0\}\ \ \text{and all}\ \ (y, z) \in \Dcal(x).
\eequation
\edit{Condition \eqref{eq:second_order_cond} is the second-order condition employed, e.g., in \cite{OberWrig2006,Wrig2002}.  It is related to the standard second-order condition as follows.  For a fixed $(y,z) \in \Dcal(x)$, the standard critical cone is
\bequationNN
  \baligned
    \Ccal_{(y,z)}(x) := \{ d \in \R{n} :\ \nabla e(x)^Td = 0,\ \ \nabla c_i(x)^Td &\leq 0\ \text{for all}\ i \in \Acal(x)\ \text{with}\ z_i = 0,\ \text{and} \\
    \nabla c_i(x)^Td &= 0\ \text{for all}\ i \in \Acal(x)\ \text{with}\ z_i > 0 \},
  \ealigned
\eequationNN
and in fact $\Ccal_{(y,z)}(x) = \Ccal(x)$ for every $(y,z) \in \Dcal(x)$.  Indeed, $\Ccal(x) \subseteq \Ccal_{(y,z)}(x)$ since $\Acal_0(x) \subseteq \{i \in \Acal(x) : z_i = 0\}$.  Conversely, any $d \in \Ccal_{(y,z)}(x)$ satisfies $\nabla f(x)^Td = -\sum_{i \in \Acal(x)} z_i \nabla c_i(x)^Td = 0$ by \eqref{eq.kkt}, so for any $(\ybar,\zbar) \in \Dcal(x)$ one has $\sum_{i \in \Acal(x)} \zbar_i \nabla c_i(x)^Td = 0$ with every term nonpositive, which yields $\nabla c_i(x)^Td = 0$ for all $i \in \Acal(x) \setminus \Acal_0(x)$.  Thus,~\eqref{eq:second_order_cond} is the standard second-order sufficient condition \cite{NoceWrig2006} imposed for every multiplier in $\Dcal(x)$, rather than for some multiplier, and so it is stronger in general.  For example, for the problem to minimize $x_1 + \thalf\|x\|_2^2$ subject to $-x_1 \leq 0$ and $-x_1 - x_2^2 \leq 0$, the origin is a strict global minimizer at which the MFCQ holds, $\Dcal(0) = \{z \in \R{2}_{\geq0} : z_1 + z_2 = 1\}$, $\Acal_0(0) = \emptyset$, and $\Ccal(0) = \{d \in \R{2} : d_1 = 0\}$.  Since the Hessian of the Lagrangian at $(0,z)$ is $\diag(1,1-2z_2)$, the standard condition holds with $z = (1,0)$, whereas \eqref{eq:second_order_cond} fails for every $z \in \Dcal(0)$ with $z_2 \geq \thalf$.  Condition \eqref{eq:second_order_cond} is implied by the strong second-order condition in~\cite{Robi1980} holding for every $(y,z) \in \Dcal(x)$ since that requires positive curvature on a subspace containing $\Ccal(x)$.  Under the LICQ, the optimal multiplier set $\Dcal(x)$ is a singleton and $\Acal_0(x) = \{i \in \Acal(x) : z_i = 0\}$, so condition~\eqref{eq:second_order_cond} reduces to the standard second-order sufficient condition in this case.}

We close this section with an observation and a lemma that will be referenced later on.  First, recall that by Farkas' theorem, one has that for any $g \in \R{n}$, $A \in \R{n \times p}$, and $B \in \R{n \times \ell}$ where $\ell \in [q]$, one and only one of the following two systems of equations and inequalities has a solution:
\begin{align*}
  g + Ay + B\zhat = 0\ \text{and}\ \zhat \geq 0 &\ \text{for some}\ (y,\zhat) \in \R{p} \times \R{\ell};\ \text{or} \\
  A^Td = 0,\ B^Td \leq 0,\ \text{and}\ g^Td < 0 &\ \text{for some}\ d \in \R{n}.
\end{align*}
Thus, at any $x \in \R{n}$ that is feasible for \eqref{eq.opt} in the sense that $e(x) = 0$ and $c(x) \leq 0$, either $\Dcal(x)$ is nonempty or for any trust-region radius $\Delta \in (0,\infty)$ the LP given by
\bequation\label{eq.improving_direction_LP}
  \min_{d \in \R{n}}\ \nabla f(x)^Td\ \ \st\ \ e(x) + \nabla e(x)^Td = 0,\ \ c(x) + \nabla c(x)^Td \leq 0,\ \ \text{and}\ \ \|d\|_\infty \leq \Delta
\eequation
has a negative optimal value.  Similarly, by the same reasoning, at any such feasible $x$, either $\Dcal(x)$ is nonempty or for any $\theta \in (0,\infty)$ the QP given by
\bequation\label{eq.improving_direction}
  \min_{d \in \R{n}}\ \nabla f(x)^Td + \thalf \theta \|d\|_2^2\ \ \st\ \ e(x) + \nabla e(x)^Td = 0\ \text{and}\ c(x) + \nabla c(x)^Td \leq 0
\eequation
has a negative optimal value.  In addition, recall the following lemma pertaining to solutions of perturbed linear systems; see, e.g., \cite[Theorem~7.29]{BurdFair2006} or \cite[Theorem~2.1.2]{OrteRhei2014}.

\blemma\label{lem.SysPertMat}
  Suppose $K \in \R{n \times n}$ is nonsingular and $b \in \R{n} \setminus \{0\}$.  Let $v \in \R{n}$ be the unique solution of $Kv = b$ $($meaning $v \neq 0$ since $b \neq 0$$)$ and let $\vtilde$ satisfy $\Ktilde \vtilde = \btilde$, where $\Ktilde = K + \delta K$ and $\btilde = b + \delta b$.  Then, for any vector-induced norm $\|\cdot\|$ and with $\cond(K) := \|K\| \|K^{-1}\|$, one has that
  \bequationNN
    \|\delta K\| \|K^{-1}\| < 1\ \ \implies\ \ \tfrac{\|v - \vtilde\|}{\|v\|} \leq \tfrac{\cond(K)}{1 - \|\delta K\| \|K^{-1}\|} \( \tfrac{\|\delta K\|}{\|K\|} + \tfrac{\|\delta b\|}{\|b\|} \).
  \eequationNN
\elemma

{\editcolor
\edit{We also make use of the following absolute form of a perturbation bound.}

\blemma\label{lem.SysPertAbs}
  \edit{Suppose $K \in \R{n \times n}$ is nonsingular, and let $\Ktilde = K + \delta K$, where for some vector-induced norm $\|\cdot\|$ one has $\|\delta K\| \|K^{-1}\| \leq \thalf$.  Then, $\Ktilde$ is nonsingular with $\|\Ktilde^{-1}\| \leq 2 \|K^{-1}\|$, and if $Kv = b$ and $\Ktilde \vtilde = \btilde$ for some $(v,b,\vtilde,\btilde) \in \R{n} \times \R{n} \times \R{n} \times \R{n}$, then one has that}
  \bequationNN
    \|\vtilde - v\| \leq 2 \|K^{-1}\| ( \|\btilde - b\| + \|\delta K\| \|v\| ).
  \eequationNN
\elemma
\bproof{Proof.}
  \edit{Since $\|K^{-1} \delta K\| \leq \|\delta K\| \|K^{-1}\| \leq \thalf$, by the properties of Neumann series $I + K^{-1} \delta K$ is nonsingular with an inverse having norm at most~2.  Thus, $\Ktilde^{-1} = (K (I + K^{-1} \delta K))^{-1} = (I + K^{-1} \delta K)^{-1} K^{-1} \implies \|\Ktilde^{-1}\| \leq 2\|K^{-1}\|$. Since subtracting $\Ktilde v = b + \delta K v$ from $\Ktilde \vtilde = \btilde$ yields $\Ktilde(\vtilde - v) = (\btilde - b) - \delta K v$, multiplying this equation by $\Ktilde^{-1}$ and employing the bound above yields the desired conclusion.}
  \qed
\eproof
}

\section{Exact Function and Derivative Values}\label{sec.exact}

The purpose of this section is to state three active-set identification results for the setting of having exact function and derivative values.  These are the results upon which we build for our new theoretical guarantees for noisy settings.

Let $x_* \in \R{n}$ be a local minimizer of problem~\eqref{eq.opt}.  Each active-set identification result that we present is based on solving a subproblem that is designed to determine an estimate of the optimal active set with respect to $x_*$, namely, $\Acal(x_*)$, at some distinct point $x \in \R{n}$.  In particular, the results that we present state that, under certain conditions, the estimate of the optimal active set is correct as long as~$x$ is sufficiently close to~$x_*$.  The first subproblem that we introduce in this section, an~LP defined with respect to problem function and derivative values at~$x$, involves computing an estimate of an element of $\Dcal(x_*)$ that can, in turn, be employed to estimate~$\Acal(x_*)$.  The second and third subproblems that we introduce, an LP and a QP also defined with respect to function and derivative values at $x$, each involve computing a step in the primal space that in turn can be employed to estimate~$\Acal(x_*)$.  Both of these LP techniques are based directly on ones established in \cite{OberWrig2006}.  The QP technique is similar to the latter LP technique from \cite{OberWrig2006}, except that we consider a QP subproblem.  Our motivation for investigating this QP technique is that this allows active-set identification to be done through subproblems that are commonly solved for computing search directions within the context of recently designed algorithms for noisy and stochastic settings.

Let $x \in \R{n}$ be given.  A strategy for estimating the active set at a nearby local minimizer $x_*$ involves first obtaining Lagrange multiplier values $(y,z) \in \R{p} \times \R{q}$, and then evaluating the KKT error given by the function $\psi : \R{n} \times \R{p} \times \R{q}$ defined as
\bequation\label{eq.psi}
  \psi(x,y,z) = \|\nabla f(x) + \nabla e(x)y + \nabla c(x)z\|_1 + \|e(x)\|_1 + \|\min\{z,-c(x)\}\|_1,
\eequation
where the minimum is taken component-wise; recall \eqref{eq.kkt}.  The motivation for obtaining Lagrange multiplier values and evaluating the KKT error can be seen through the following from~\cite{OberWrig2006}; see also \cite{FaccFiscKanz1998,HageSeet1999,Wrig2002}.

\btheorem\label{th.OberWrig2.1} \cite[Theorem~2.1]{OberWrig2006}
  Suppose that, with respect to~\eqref{eq.opt}, the MFCQ and the second-order condition~\eqref{eq:second_order_cond} hold at a local minimizer $x_* \in \R{n}$.  Then, there exist $\epsilon \in (0,1]$ and $C \in \R{}_{>0}$ such that, for all $(x,y,z) \in \R{n} \times \R{p} \times \R{q}$ with $z \geq 0$ and $\dist(\{(x,y,z)\},\{x_*\} \times \Dcal(x_*)) \leq \epsilon$, one has that $C^{-1} \psi(x,y,z) \leq \dist(\{(x,y,z)\}, \{x_*\} \times \Dcal(x_*)) \leq C\psi(x,y,z)$.
\etheorem

Motivated by Theorem~\ref{th.OberWrig2.1}, a technique is proposed in \cite{OberWrig2006} for estimating the active set from a primal point $x \in \R{n}$.  At $x$, the idea is to compute Lagrange multipliers by minimizing $\psi(x,y,z)$ with respect to $(y,z)$ subject to $z \geq 0$, i.e., to evaluate $\omega : \R{n} \to \R{}$ at $x$ that is defined by
\bequation\label{eq.def_omega}
  \omega(x) := \min_{(y,z) \in \R{p} \times \R{q}_{\geq0}}\ \psi(x,y,z).
\eequation
This can be done by solving a linear program with equilibrium constraints (LPEC), which in turn can be solved as a mixed-integer linear program (MIP).  However, since it can be computationally expensive to solve such a problem, a technique is proposed in \cite{OberWrig2006} to solve an LP approximation to the LPEC, the solution of which can in turn be used to estimate the active set.  Let us refer to this as the LP-LPEC approach.  This is the approach that we extend to the noisy setting in the next section.

To state the LP-LPEC approach, let us first introduce $\kappa : \R{n} \times \R{p} \times \R{q} \to \R{}$, $\rho : \R{n} \times \R{p} \times \R{q}_{\geq0} \to \R{}$, and $\bar\rho : \R{n} \times \R{p} \times \R{q}_{\geq0} \to \R{}$ defined respectively by
\bsubequations\label{def.kappa_rhos}
  \begin{align}
    \kappa(x,y,z) &:= \|\nabla f(x) + \nabla e(x)y + \nabla c(x) z\|_1 + \|e(x)\|_1, \\
    \rho(x,y,z) &:= \kappa(x,y,z) + \sum_{\{i : c_i(x)<0\}} -c_i(x) z_i + \sum_{\{i : c_i(x) \geq 0\}} c_i(x), \label{eq.def_rho} \\ \text{and}\ \ 
    \bar{\rho}(x,y,z) &:= \kappa(x,y,z) + \sum_{\{i : c_i(x) < 0\}} (-c_i(x) z_i)^{1/2} + \sum_{\{i : c_i(x) \geq 0\}} c_i(x). \label{eq.def_bar_rho}
  \end{align}
\esubequations
At $x \in \R{n}$, an approximate solution of the aforementioned LPEC can then be obtained by solving the following subproblem for the given constant $M \in \R{}_{>0}$:
\bequation\label{eq.lp}
  \baligned
    &\ \min_{(y,z) \in \R{p} \times \R{q}}\ \rho(x,y,z)\ \ \st\ \ 0 \leq z \leq M \ones \\
    \text{where}\ \ &\ M := \max \left\{ \|c(x_*)\|_\infty, \max_{(y_*,z_*) \in \Dcal(x_*)} \left\| \bbmatrix y_* \\ z_* \ebmatrix \right\|_\infty \right\} + 1.
  \ealigned
\eequation
\edit{The constant $M$ is finite if and only if the MFCQ holds at $x_*$; see~\cite{Gauv1977}. A variant of \eqref{eq.lp} where the upper bound is dropped is considered in \cite{OberWrig2006} and in our experiments it proved to be effective, although this variant is not covered by the theory in this paper. The following remark concerns a practical implementation of~\eqref{eq.lp}.}

{\editcolor
\begin{remark}\label{rem:M}
  \textit{
    Since the value of $M$ in~\eqref{eq.lp} is usually unknown in practice, one could set an initial estimate for $M$ and double it every time a bound constraint involving this estimate is active at the computed solution of \eqref{eq.lp}, or of its noisy counterpart \eqref{prob.lp-lpec} in \S\ref{sec.lagrange} (similar to penalty-parameter updates for exact penalty methods \cite{Burk1991,ByrdNoceWalt2008,HanMang1979}).  If this strategy terminates, then none of these bound constraints is active at the final computed solution, so by convexity of the subproblem this solution also solves the subproblem with these bound constraints removed.  Consequently, if the final estimate is at most the value of $M$ defined in~\eqref{eq.lp}, then the final computed solution also solves the subproblem with that value of $M$, and our identification results apply to it.  If instead the final estimate, call it $M'$, exceeds $M$, then for any $\Mbar \geq M'$ the upcoming results of \S\ref{sec.lagrange} still hold with $M$ replaced by $M'$ in \eqref{prob.lp-lpec}, with thresholds that depend on $\Mbar$ but not on $M'$, since the proof of Theorem~\ref{the.LP_main_theorem} uses the bounds in \eqref{prob.lp-lpec} only to ensure that certain points with $\ell_\infty$-norm at most $M - \thalf$ are feasible and that Lemma~\ref{lem.basic_lemma_2} applies at the computed solution.
  }
\end{remark}
}

A solution of \edit{subproblem~\eqref{eq.lp}} can be obtained by solving the LP
\bequation
  \baligned
    \min_{(y,z,r,s) \in \R{p} \times \R{q} \times \R{n} \times \R{n}} &\ \ones^Tr + \ones^Ts + \sum_{\{i : c_i(x)<0\}} -c_i(x) z_i \\
    \st &\ \nabla f(x) + \nabla e(x)y + \nabla c(x) z = r - s,\ r \geq 0,\ s \geq 0,\ \text{and}\ M \ones \geq z \geq 0.
  \ealigned
\eequation
Given $x$, let $(y_x,z_x)$ be an optimal solution of \eqref{eq.lp}.  An active-set estimate can then be obtained through evaluation of the function $\bar\rho$; in particular, let the estimate be
\bequation\label{eq.lovelylara}
  \Acal_{{\rm LP}}(y_x,z_x;x) := \{i \in [q] : c_i(x) \geq -(\beta \bar\rho(x,y_x,z_x))^\sigma\},
\eequation
where $\beta \in \R{}_{>0}$ and $\sigma \in (0,1)$ are user-defined constants.  \edit{Theorem~\ref{th.lplpec} below, from~\cite{OberWrig2006}, shows that this LP-LPEC approach can be employed under certain assumptions to provide an accurate active-set estimate from $x$ near a local minimizer $x_*$.  This theorem (and Theorem~\ref{the.LP_main_theorem} in the noisy setting) holds for any $\beta \in \R{}_{>0}$ and $\sigma \in (0,1)$, and the radius of the neighborhood on which identification is guaranteed depends on both of these parameters.  The parameter $\sigma$, in particular, yields a tradeoff.  Since $\bar\rho(x,y_x,z_x) \to 0$ as $x \to x_*$, a value of $\sigma$ near~1 makes the threshold $(\beta \bar\rho(x,y_x,z_x))^\sigma$ small, and hence rejects inactive constraints on a larger neighborhood. On the other hand, a value of $\sigma$ near~0 keeps the threshold large relative to this residual, which makes it more robust to errors in $\bar\rho(x,y_x,z_x)$. In \S\ref{sec.noisy}, this robustness is what is needed to identify active constraints subject to noise.  The parameter $\beta$ should be scaled to the problem data; see~\S\ref{sec.numerical}.}

\btheorem\label{th.lplpec} \cite[Theorem~3.7]{OberWrig2006}
  Suppose that, with respect to~\eqref{eq.opt}, the MFCQ and the second-order condition~\eqref{eq:second_order_cond} hold at a local minimizer $x_* \in \R{n}$.  Then, with $\epsilon \in (0,1]$ defined as in Theorem~\ref{th.OberWrig2.1} and any $\beta \in \R{}_{>0}$ and $\sigma \in (0,1)$, there exists $\hat\epsilon \in (0,\epsilon]$ such that for all $x \in \R{n}$ with $\|x - x_*\|_2 \leq \hat\epsilon$, the following hold.
  \benumerate
    \item[(a)] There exists a solution $(y_x,z_x)$ of \eqref{eq.lp} such that $\dist(\{(y_x,z_x)\}, \Dcal(x_*)) \leq \tfrac12 \epsilon$.
    \item[(b)] Any $(y_x,z_x)$ solving \eqref{eq.lp} has $M^{-1} \rho(x,y_x,z_x) \leq \omega(x) \leq \bar\rho(x,y_x,z_x)$.
    \item[(c)] Any $(y_x,z_x)$ solving \eqref{eq.lp} yields $\Acal_{{\rm LP}}(y_x,z_x;x) = \Acal(x_*)$.
  \eenumerate
\etheorem
This completes our discussion of the LP-LPEC approach that we extend in the next section.

Let us now introduce the other LP-based approach and a related QP-based approach that we extend to the noisy setting as well.  \edit{If a local minimizer $x_*$ is a primal KKT point and $x = x_*$, then both~\eqref{eq.improving_direction_LP} and~\eqref{eq.improving_direction} are feasible with an optimal value of zero.  In particular, $d_* = 0$ is an optimal solution of each, and it is the unique optimal solution of \eqref{eq.improving_direction} by strong convexity.  (For the LP \eqref{eq.improving_direction_LP}, any feasible $d$ in the critical cone with $\nabla f(x_*)^Td = 0$ is also optimal, so uniqueness cannot be expected.)  At a point $x \neq x_*$, a negative optimal value is to be expected and is precisely what a primal search direction requires. The difficulty, however, is that the linearized constraints may be inconsistent, in which case \eqref{eq.improving_direction_LP} and \eqref{eq.improving_direction} are infeasible.  The standard remedy, which goes back to the $S\ell_1$QP method of Fletcher \cite{Flet1985} and the robust SQP methods of Burke and Han \cite{Burk1989,BurkHan1989}, and which is used in certain solvers in so-called elastic mode \cite{GillMurrSaun2002,ByrdNoceWalt2006}, is to penalize violation of the linearized constraints with an exact $\ell_1$-norm penalty.  To that end, and in particular} to compute a primal search direction that can also be used to determine an optimal active-set estimate at~$x$, one can choose $(\Delta,\nu) \in \R{}_{>0} \times \R{}_{>0}$ and solve the nonsmooth subproblem
\bequation\label{eq.FEEL_LP}
  \min_{d \in \R{n}}\ \nabla f(x)^Td + \nu \|e(x) + \nabla e(x)^Td\|_1 + \nu \| [c(x) + \nabla c(x)^Td]_+ \|_1\ \ \st\ \ \|d\|_\infty \leq \Delta
\eequation
or choose $(\theta,\nu) \in \R{}_{>0} \times \R{}_{>0}$ and solve the nonsmooth subproblem
\bequation\label{eq.FEEL}
  \min_{d \in \R{n}}\ \nabla f(x)^Td + \thalf \theta \|d\|_2^2 + \nu \|e(x) + \nabla e(x)^Td\|_1 + \nu \| [c(x) + \nabla c(x)^Td]_+ \|_1.
\eequation
Alternatively, one can solve the constrained LP subproblem
\bequation\label{eq.sub_LP}
  \baligned
    \min_{(d,r,s,t) \in \R{n} \times \R{p} \times \R{p} \times \R{q}} &\ \nabla f(x)^Td + \nu(\ones^Tr + \ones^Ts + \ones^Tt) \\
    \st &\ e(x) + \nabla e(x)^Td = r - s,\ \ c(x) + \nabla c(x)^Td \leq t, \\
        &\ r \geq 0,\ \ s \geq 0,\ \ t \geq 0,\ \ \text{and}\ \ \|d\|_\infty \leq \Delta
  \ealigned
\eequation
that is equivalent to \eqref{eq.FEEL_LP}, or one can solve the constrained QP subproblem
\bequation\label{eq.sub}
  \baligned
    \min_{(d,r,s,t) \in \R{n} \times \R{p} \times \R{p} \times \R{q}} &\ \nabla f(x)^Td + \nu(\ones^Tr + \ones^Ts + \ones^Tt) + \thalf \theta \|d\|_2^2 \\
    \st &\ e(x) + \nabla e(x)^Td = r - s,\ \ c(x) + \nabla c(x)^Td \leq t, \\
        &\ r \geq 0,\ \ s \geq 0,\ \ \text{and}\ \ t \geq 0
  \ealigned
\eequation
that is equivalent to \eqref{eq.FEEL}. These subproblems are feasible for any such $(x,\Delta,\theta,\nu)$; indeed, $d = 0$, $r = [e(x)]_+$, $s = [e(x)]_-$, and $t = [c(x)]_+$ form a feasible point for either of the subproblems.  \edit{The equivalence of nonsmooth \eqref{eq.FEEL_LP}--\eqref{eq.FEEL} with smooth constrained \eqref{eq.sub_LP}--\eqref{eq.sub}, and the form of the duals \eqref{eq.sub_dual_LP}--\eqref{eq.sub_dual} below, are classical; see \cite{Flet1985,Burk1989}, and see \cite{ByrdNoceWalt2008} for strategies for updating~$\nu$.  Our interest in \eqref{eq.FEEL} specifically is that it is the subproblem solved in several recently proposed methods for noisy and stochastic constrained optimization \cite{BeraCurtONeiRobi2023,BeraCurtRobiZhou2021,CurtJianWang2024,CurtONeiRobi2024,CurtRobiZhou2024b,CurtRobiZhou2024}.}  
The dual of~\eqref{eq.sub_LP} can be written as the following (equivalent to a constrained LP)
\bequation\label{eq.sub_dual_LP}
  \baligned
    \max_{(\alpha,\beta) \in \R{p} \times \R{q}} &\ e(x)^T\alpha + c(x)^T\beta - \Delta \|\nabla f(x) + \nabla e(x)\alpha + \nabla c(x)\beta\|_1 \\
    \st &\ - \nu \ones \leq \alpha \leq \nu \ones\ \text{and}\ 0 \leq \beta \leq \nu \ones,
  \ealigned
\eequation
whereas the dual of~\eqref{eq.sub} can be written as the QP
\bequation\label{eq.sub_dual}
  \baligned
    \max_{(\alpha,\beta) \in \R{p} \times \R{q}} &\ e(x)^T\alpha + c(x)^T\beta - \tfrac{1}{2} \theta^{-1} \|\nabla f(x) + \nabla e(x)\alpha + \nabla c(x)\beta\|_2^2 \\
    \st &\ - \nu \ones \leq \alpha \leq \nu \ones\ \text{and}\ 0 \leq \beta \leq \nu \ones.
  \ealigned
\eequation
These dual subproblems are also feasible for any $(x,\Delta,\theta,\nu)$; indeed, $(\alpha,\beta) = (0,0)$ is feasible for both dual subproblems.  Consequently, by weak duality, the optimal values of \eqref{eq.sub_LP}, \eqref{eq.sub}, \eqref{eq.sub_dual_LP}, and \eqref{eq.sub_dual} are all finite, and in each case these values are attained by some primal solution $(d,r,s,t)$ and some dual solution $(\alpha,\beta)$, respectively; see, e.g., \cite[Proposition~1.4.12]{Bert2009}.

Given any $(x,\Delta,\nu)$ and a corresponding solution $(d_{x,\Delta,\nu},r_{x,\Delta,\nu},s_{x,\Delta,\nu},t_{x,\Delta,\nu})$ of subproblem~\eqref{eq.sub_LP}, an estimate of the optimal active set is given by
\bequation\label{eq.active_LPtr}
  \Acal_{{\rm LP}}(d_{x,\Delta,\nu}; x) := \{i \in [q] : c_i(x) + \nabla c_i(x)^Td_{x,\Delta,\nu} \geq 0\}.
\eequation
Similarly, given any $(x,\theta,\nu)$ and a corresponding solution $(d_{x,\theta,\nu},r_{x,\theta,\nu},s_{x,\theta,\nu},t_{x,\theta,\nu})$ of subproblem~\eqref{eq.sub}, an estimate of the optimal active set is given by
\bequation\label{eq.active_QP}
  \Acal_{{\rm QP}}(d_{x,\theta,\nu}; x) := \{i \in [q] : c_i(x) + \nabla c_i(x)^Td_{x,\theta,\nu} \geq 0\}.
\eequation
Theorem~\ref{th.OberWrig} below shows that this QP-based approach is successful when $x$ is sufficiently close to a local minimizer of \eqref{eq.opt} at which the LICQ and strict complementarity hold.  The same kind of result can be proved for the LP-based approach; indeed, this is what has been proved in Lemma~3.1, Theorem 3.2, and Corollary 3.3 in~\cite{OberWrig2006}.  (For the sake of brevity, we do not restate those results here, since they are identical to our theorem below, except for the use of an LP subproblem with a trust-region constraint instead of our QP.)  Our proof of the following theorem follows closely the analysis in~\cite{OberWrig2006}, although details are different due to our use of a QP subproblem\edit{; see Appendix~\ref{app.OberWrig}}.

\btheorem\label{th.OberWrig}\label{th.QP_exact}
  Suppose that, with respect to~\eqref{eq.opt}, the MFCQ holds at a local minimizer $x_* \in \R{n}$.
  \benumerate
    \item[(a)] There exists $(\epsilon,\underline\theta) \in \R{}_{>0} \times \R{}_{>0}$ such that, for any $(x,\theta,\nu) \in \R{n} \times \R{}_{>\underline\theta} \times \R{}_{>0}$ with $\|x - x_*\|_2 \leq \epsilon$, the solution $d_{x,\theta,\nu}$ of \eqref{eq.FEEL} yields $\Acal_{{\rm QP}}(d_{x,\theta,\nu}; x) \subseteq \Acal(x_*)$.
    \item[(b)] \edit{For} any $\zeta \in (0,1)$, there exists $(\epsilon,\underline\theta,\underline\nu) \in \edit{\R{}_{>0} \times \R{}_{>0} \times \R{}_{>0}}$ such that for any $(x,\theta,\nu) \in \R{n} \times \R{}_{>\underline\theta} \times \R{}_{>\underline\nu}$ with $\|x - x_*\|_2 \leq \epsilon$ and $\theta \leq \|x - x_*\|_2^{-\zeta}$, the solution $d_{x,\theta,\nu}$ of \eqref{eq.FEEL} yields
    \bequation\label{eq.ran_pos_stat}
      -\nabla f(x_*)\in {\rm range}[\nabla e(x_*)] + {\rm pos}[\{\nabla c_i(x_*) : i \in \Acal_{{\rm QP}}(d_{x,\theta,\nu}; x)\}],
    \eequation
    where $\text{range}[\cdot]$ denotes the range space of the matrix input and $\text{pos}[\cdot]$ denotes the set of nonnegative combinations of the input set of vectors.
    \item[(c)] If the LICQ holds at $x_*$, then for any $(\zeta,\epsilon,\underline\theta,\underline\nu,x,\theta,\nu)$ satisfying the conditions of part (b), the solution $d_{x,\theta,\nu}$ of \eqref{eq.FEEL} yields $\Acal(x_*) \setminus \Acal_0(x_*) \subseteq \Acal_{{\rm QP}}(d_{x,\theta,\nu}; x) \subseteq \Acal(x_*)$. If strict complementarity also holds at $x_*$, then $\Acal_{{\rm QP}}(d_{x,\theta,\nu}; x) = \Acal(x_*)$.
  \eenumerate
\etheorem

The proof of Theorem~\ref{th.QP_exact} reveals the influence of the parameters of the QP approach. Intuitively, the lower threshold $\underline\theta$ ensures that $d_{x,\theta,\nu}$ is sufficiently small in norm such that the step in \eqref{eq.active_QP} does not ``go too far'' to include indices of constraints that are inactive at $x_*$.  The lower threshold~$\underline\nu$, on the other hand, is a typical type of bound for exact penalization; it ensures that, when possible, the solution of \eqref{eq.sub} will have the appropriate auxiliary variable values equal to zero.  \edit{Part~(b) of Theorem~\ref{th.QP_exact} requires $\underline\theta < \theta \leq \|x-x_*\|_2^{-\zeta}$, which is a range that is nonempty for $x$ sufficiently close to $x_*$. The lower bound ensures that the step is sufficiently small while the upper bound ensures that the regularization does not dominate the linear model.}

\section{Noisy Function and Derivative Values}\label{sec.noisy}

\edit{We now assume} that, at any $x \in \R{n}$, one only has access to approximations of the values of $f$, $e$, and $c$ and their first-order derivatives.  More precisely, we assume that there exist positive error bounds $(\epsilon_f, \epsilon_e, \epsilon_c, \epsilon_{\nabla f}, \epsilon_{\nabla e}, \epsilon_{\nabla c})$ such that at all $x \in \R{n}$ one can only compute $\fapprox(x) \in \R{}$, $\eapprox(x) \in \R{p}$, $\capprox(x) \in \R{q}$, $\fgradapprox(x) \in \R{n}$, $\egradapprox(x) \in \R{n \times p}$, and $\cgradapprox(x) \in \R{n \times q}$ satisfying
\bequation \label{eq.error_bounds}
  \begin{aligned}
    \|f(x) - \fapprox(x)\|_2 \leq \epsilon_f,\ \ \|e(x) - \eapprox(x)\|_2 &\leq \epsilon_e,\ \ \|c(x) - \capprox(x)\|_2 \leq \epsilon_c, \\ 
    \|\nabla f(x) - \fgradapprox(x)\|_2 \leq \epsilon_{\nabla f},\ \ \|\nabla e(x) - \egradapprox(x)\|_2 &\leq \epsilon_{\nabla e},\ \ \|\nabla c(x) - \cgradapprox(x)\|_2 \leq \epsilon_{\nabla c}.
  \end{aligned}
\eequation
\edit{Here and throughout, vectors are measured in the $\ell_2$-norm and matrices in the induced (spectral) norm, so, e.g., $|c_i(x) - \capprox_i(x)| \leq \epsilon_c$ for every $i \in [q]$.  The results of \S\ref{sec.lagrange} and \S\ref{sec.primal_step} are statements about any realization of the approximate values satisfying \eqref{eq.error_bounds}. They hold regardless of whether the noise is deterministic or stochastic, and no property of the error distribution other than \eqref{eq.error_bounds} is used.  In \S\ref{sec.stochastic_algorithm}, our probabilistic results are based on the probability with which \eqref{eq.error_bounds} holds.}  For convenience, let us define the vector
\bequation\label{eq.epsilon_error}
  \epsilonerror :=[\epsilon_f, \epsilon_e, \epsilon_c, \epsilon_{\nabla f}, \epsilon_{\nabla e}, \epsilon_{\nabla c}]^T.
\eequation

\subsection{Identification through Lagrange Multipliers}\label{sec.lagrange}

We first describe a technique for identifying the active set by estimating Lagrange multipliers at a primal point $x \in \R{n}$. This approach extends the LP-LPEC method described in \S\ref{sec.exact}, adapting it to noisy settings.

With the goal of again constructing a two-sided bound on $\omega(x)$, we follow \eqref{def.kappa_rhos} to define the noisy counterparts $\tilde\kappa : \R{n} \times \R{p} \times \R{q} \to \R{}$, $\tilde\rho : \R{n} \times \R{p} \times \R{q}_{\geq0} \to \R{}$, and $\tilde{\bar\rho} : \R{n} \times \R{p} \times \R{q}_{\geq0} \to \R{}$ by
\bsubequations\label{def.kappa_rhos_noisy}
  \begin{align}
    \tilde\kappa(x,y,z) &:= \|\fgradapprox(x) + \egradapprox(x) y + \cgradapprox(x) z\|_1 + \|\eapprox(x)\|_1, \\
    \tilde\rho(x,y,z) &:= \tilde\kappa(x,y,z) + \sum_{\{i : \capprox_i(x)<0\}} -\capprox_i(x) z_i + \sum_{\{i : \capprox_i(x) \geq 0\}} \capprox_i(x), \label{eq.def_tilde_rho} \\ \text{and}\ \ 
    \tilde{\bar{\rho}}(x,y,z) &:= \tilde\kappa(x,y,z) + \sum_{\{i : \ctilde_i(x) < 0\}} (-\ctilde_i(x) z_i)^{1/2} + \sum_{\{i : \ctilde_i(x) \geq 0\}} \ctilde_i(x). \label{eq.def_tilde_bar_rho}
  \end{align}
\esubequations

The following relationship can be established for $\tilde\rho$ and $\tilde{\bar{\rho}}$, drawing on \cite[Lemma 3.5]{OberWrig2006}.

\blemma\label{lem.basic_lemma_1}
  For any $(x,y,z) \in \R{n} \times \R{p} \times \R{q}_{\geq0}$, $\tilde{\bar{\rho}}(x,y,z) \leq \tilde{\rho}(x,y,z) + q^{1/2} \tilde{\rho}(x,y,z)^{1/2}$.
\elemma
\bproof{Proof.}
  \edit{The proof follows the same logic as \cite[Lemma~3.5]{OberWrig2006}.}
  \qed
\eproof

A relationship between $\tilde{\rho}$, $\tilde{\bar{\rho}}$, and $\psi$ (recall \eqref{eq.psi}) is given in our next lemma.  The existence of the error bounds in the lemma is crucial for our subsequent analysis for the noisy setting.

{\editcolor
\blemma\label{lem.basic_lemma_2}
  Let $\Mbar \in \R{}_{\geq1}$, define the error bounds
  \bequation\label{eq.eps_rho_defs}
    \begin{aligned}
      \epsilon_\kappa &:= p^{1/2}\epsilon_e + n^{1/2}(\epsilon_{\nabla f} + p^{1/2}\Mbar\epsilon_{\nabla e} + q^{1/2}\Mbar\epsilon_{\nabla c}),\ \ \epsilon_\rho := \epsilon_\kappa + q\Mbar\epsilon_c,\ \ \text{and}\ \ \epsilon_\psi := \epsilon_\kappa + q\epsilon_c,
    \end{aligned}
  \eequation
  and let $\tilde\psi(x,y,z) := \tilde\kappa(x,y,z) + \|\min\{z,-\capprox(x)\}\|_1$ denote the KKT error \eqref{eq.psi} evaluated with the approximations in place of the exact values.  Then, $\epsilon_\rho \to 0$ and $\epsilon_\psi \to 0$ as $\epsilonerror \to 0$ $($recall \eqref{eq.epsilon_error}$)$. Furthermore, for all $(x,y,z) \in \R{n} \times \R{p} \times \R{q}$ with $z \geq 0$ and $\|(c(x),y,z)\|_\infty \leq \Mbar$ one has with these definitions that
  \begin{subequations}
  \begin{align}
    |\rho(x,y,z) - \tilde\rho(x,y,z)| \leq \epsilon_\rho\ \ \text{and}\ \ |\psi(x,y,z) - \tilde\psi(x,y,z)| &\leq \epsilon_\psi, \label{eq.rho_pert} \\
    \tilde\psi(x,y,z) &\leq \tilde{\bar\rho}(x,y,z), \label{eq.psitilde_rhobar} \\
    \text{which further imply}\ \ \Mbar^{-1}\(\tilde\rho(x,y,z) - \epsilon_\rho\) \leq \psi(x,y,z) &\leq \tilde{\bar\rho}(x,y,z) + \epsilon_\psi. \label{eq.ine_rhos}
  \end{align}
  \end{subequations}
  Moreover, if such a tuple $(x,y,z)$ has $\rho(x,y,z) > 2\epsilon_\rho$, then the left-hand side of \eqref{eq.ine_rhos} is strictly positive.
\elemma
\bproof{Proof.}
  The fact that $\epsilon_\rho \to 0$ and $\epsilon_\psi \to 0$ as $\epsilonerror \to 0$ follows directly from \eqref{eq.epsilon_error} and \eqref{eq.eps_rho_defs}. Now consider $(x,y,z)$ as stated in the lemma and, for notational convenience, abbreviate $c_i := c_i(x)$ and $\capprox_i := \capprox_i(x)$ so that $|c_i - \capprox_i| \leq \epsilon_c$ for all $i \in [q]$ by \eqref{eq.error_bounds}.  For all $i \in [q]$, define the function $\phi_i : \R{} \to \R{}$ by $\phi_i(t) = -tz_i$ for all $t<0$ and $\phi_i(t) = t$ for all $t \geq 0$ so that, by \eqref{def.kappa_rhos} and \eqref{def.kappa_rhos_noisy}, one has that
  \bequationNN
    \rho(x,y,z) = \kappa(x,y,z) + \sum_{i \in [q]} \phi_i(c_i)\ \ \text{and}\ \ \tilde\rho(x,y,z) = \tilde\kappa(x,y,z) + \sum_{i \in [q]} \phi_i(\capprox_i).
  \eequationNN
  The leading-term difference can be bounded since \eqref{eq.error_bounds}, $\|y\|_2 \leq p^{1/2}\Mbar$, and $\|z\|_2 \leq q^{1/2}\Mbar$ give
  \begin{align*}
    &\ |\kappa(x,y,z) - \tilde\kappa(x,y,z)| \\
    \leq&\ \|e(x) - \eapprox(x)\|_1 + \|(\nabla f(x) - \fgradapprox(x)) + (\nabla e(x) - \egradapprox(x))y + (\nabla c(x) - \cgradapprox(x))z\|_1 \\
    \leq&\ p^{1/2} \epsilon_e + n^{1/2}\(\epsilon_{\nabla f} + \epsilon_{\nabla e}\|y\|_2 + \epsilon_{\nabla c}\|z\|_2\)\ \leq\ \epsilon_\kappa.
  \end{align*}
  Next, observe that $\phi_i$ is piecewise linear and continuous with slope either $-z_i$ or~1 at any point, hence it is Lipschitz with constant $\max\{1,z_i\} \leq \Mbar$.  Therefore $|\phi_i(c_i) - \phi_i(\capprox_i)| \leq \Mbar\epsilon_c$, and summing over $i \in [q]$ and combining with the bound above gives the first inequality in \eqref{eq.rho_pert}.  It remains to establish the second inequality in \eqref{eq.rho_pert} along with \eqref{eq.psitilde_rhobar}, \eqref{eq.ine_rhos}, and the final claim.  With $\tilde\psi$ as defined, for each $i \in [q]$ one has $|\min\{z_i,-\capprox_i\}| = \capprox_i$ if $\capprox_i \geq 0$ and $|\min\{z_i,-\capprox_i\}| \leq (-\capprox_iz_i)^{1/2}$ if $\capprox_i < 0$ since the minimum of two nonnegative numbers is at most their geometric mean. Hence $\tilde\psi(x,y,z) \leq \tilde{\bar\rho}(x,y,z)$, which is \eqref{eq.psitilde_rhobar}.  Moreover, since $\psi(x,y,z) = \kappa(x,y,z) + \sum_{i \in [q]} |\min\{z_i,-c_i\}|$ and $t \mapsto |\min\{z_i,-t\}|$ is Lipschitz with constant~1 for each $i \in [q]$, one has $|\psi(x,y,z) - \tilde\psi(x,y,z)| \leq |\kappa(x,y,z) - \tilde\kappa(x,y,z)| + \sum_{i \in [q]} |c_i - \capprox_i| \leq \epsilon_\kappa + q\epsilon_c = \epsilon_\psi$, which is the second inequality in \eqref{eq.rho_pert}.  Thus, with \eqref{eq.psitilde_rhobar}, $\psi(x,y,z) \leq \tilde\psi(x,y,z) + \epsilon_\psi \leq \tilde{\bar\rho}(x,y,z) + \epsilon_\psi$, which is the upper bound in \eqref{eq.ine_rhos}.  Now, for the lower bound in \eqref{eq.ine_rhos}, one has $\psi(x,y,z) \geq \Mbar^{-1}\rho(x,y,z)$ (see \cite[Lemma~3.6]{OberWrig2006}) since $\Mbar \geq 1$, the terms of $\rho(x,y,z)$ and $\psi(x,y,z)$ with $c_i \geq 0$ coincide, and $c_i < 0$ implies $-c_iz_i \leq \max\{z_i,-c_i\}\min\{z_i,-c_i\} \leq \Mbar\min\{z_i,-c_i\}$ by $\|(c(x),y,z)\|_\infty \leq \Mbar$; combining with the first inequality in \eqref{eq.rho_pert} yields $\psi(x,y,z) \geq \Mbar^{-1}(\tilde\rho(x,y,z) - \epsilon_\rho)$.  Finally, if $\rho(x,y,z) > 2\epsilon_\rho$, then \eqref{eq.rho_pert} gives $\tilde\rho(x,y,z) - \epsilon_\rho \geq \rho(x,y,z) - 2\epsilon_\rho > 0$, so the left-hand side of \eqref{eq.ine_rhos} is strictly positive.
  \qed
\eproof
}

The LP-LPEC-based active-set identification strategy can now be established as follows.  First, one computes a solution $(y_x,z_x)$ of the optimization problem
\bequation\label{prob.lp-lpec}
  \min_{(y,z) \in \R{p} \times \R{q}}\ \tilde{\rho}(x,y,z)\ \st\ \edit{-M\ones \leq y \leq M\ones,}\ \ 0\leq z \leq M\ones,
\eequation
where $M$ is defined in \eqref{eq.lp}. \edit{Concerning the challenge of not (typically) knowing the value of $M$ in practice, see Remark~\ref{rem:M}.} The solution of \eqref{prob.lp-lpec} can be obtained by solving the LP
\begin{align*}
  \min_{(y,z,r,s) \in \R{p} \times \R{q} \times \R{n} \times \R{n}} &\ \ones^Tr + \ones^Ts + \sum_{\{i : \capprox_i(x)<0\}} -\capprox_i(x) z_i \\
  \st &\ \fgradapprox(x) + \egradapprox(x)y + \cgradapprox(x) z = r - s, \\
  &\ r \geq 0,\ \ s \geq 0,\ \ \edit{-M\ones \leq y \leq M\ones,}\ \ 0 \leq z \leq M\ones.
\end{align*}
The active-set estimate is then defined, similarly to \eqref{eq.lovelylara}, by
\bequation\label{eq.active_LP}
  \widetilde\Acal_{{\rm LP}}(y_x, z_x; x) := \{ i \in [q] : \capprox_i(x)\geq-(\beta\tilde{\bar{\rho}}(x,y_x,z_x))^\sigma\}
\eequation
where $\beta \in \R{}_{>0}$ and $\sigma \in (0,1)$ are user-defined constants.

In the following theorem, we demonstrate that this methodology yields a two-sided bound on the value of $\omega(x)$ and can accurately identify the optimal active set at a point near a local minimizer.

{\editcolor
\btheorem\label{the.LP_main_theorem}
  Suppose that, with respect to~\eqref{eq.opt}, the MFCQ and~\eqref{eq:second_order_cond} hold at a local minimizer $x_* \in \R{n}$.  Let $\epsilon$ and $C$ be defined as in Theorem~\ref{th.OberWrig2.1}, let $\epsilon' \in (0,\epsilon/2]$ be the radius as in \cite[Theorem~3.4]{OberWrig2006}, let $\Mbar \geq M$ with $M$ defined in~\eqref{eq.lp}, and let~$\epsilon_\rho$ and $\epsilon_\psi$ be defined by \eqref{eq.eps_rho_defs}.  Consider any $\beta \in \R{}_{>0}$ and $\sigma \in (0,1)$.  Then, there exist $\hat\epsilon \in (0,\epsilon']$ and $\bar\epsilon \in \R{}_{>0}$ such that for any $x \in \R{n}$ with $\|x - x_*\|_2 \leq \hat\epsilon$ and any approximations satisfying \eqref{eq.error_bounds} with $\|\epsilonerror\|_\infty \leq \bar\epsilon$, the following hold, where $(y_x,z_x)$ denotes any solution of \eqref{prob.lp-lpec}.
  \benumerate[label=$($\alph*$)$]
    \item If $\rho(x,y_x,z_x) > 2\epsilon_\rho$, then $0 < \Mbar^{-1}\(\tilde\rho(x,y_x,z_x) - \epsilon_\rho\) \leq \omega(x) \leq \tilde{\bar\rho}(x,y_x,z_x) + \epsilon_\psi$.
    \item If $\epsilon_\psi \leq (2C)^{-1}\|x - x_*\|_2$, then $\widetilde\Acal_{{\rm LP}}(y_x, z_x; x) = \Acal(x_*)$.
  \eenumerate
\etheorem
\bproof{Proof.}
  Throughout the proof, let us restrict attention to $x \in \R{n}$ satisfying $\|x - x_*\|_2 \leq \epsilon'$ and approximations satisfying \eqref{eq.error_bounds}. The (potentially smaller) radius  $\hat\epsilon \in (0,\epsilon']$ and the (potentially tighter) bound on the approximations $\bar\epsilon$ that are required by the theorem are developed in the proof below.

  Let us begin by observing that there exists a minimizer of $\psi(x,\cdot,\cdot)$ that is feasible for \eqref{prob.lp-lpec}.  Indeed, by \cite[Theorem~3.4]{OberWrig2006}, there exists $(\yhat,\zhat)$ attaining the minimum in \eqref{eq.def_omega}, where in particular $\zhat \geq 0$ and $\omega(x) = \psi(x,\yhat,\zhat)$, with $\dist(\{(\yhat,\zhat)\},\Dcal(x_*)) \leq \epsilon/2$.  Since $\Dcal(x_*)$ is nonempty and closed, this distance is attained, and since every element of $\Dcal(x_*)$ has $\ell_\infty$-norm at most $M - 1$ by \eqref{eq.lp}, it follows with $\epsilon \leq 1$ that $\|(\yhat,\zhat)\|_\infty \leq M - 1 + \epsilon/2 \leq M - \thalf$.  Thus, $(\yhat,\zhat)$ is feasible for \eqref{prob.lp-lpec}, as claimed.

  Next, let us derive a lower bound for $\bar\rho(x,y_x,z_x)$ that does not depend on the noise, namely,
  \bequation\label{eq.rhobar_lower}
    \bar\rho(x,y_x,z_x) \geq \psi(x,y_x,z_x) \geq \omega(x) = \psi(x,\yhat,\zhat) \geq C^{-1}\|x - x_*\|_2.
  \eequation
  For the first inequality, observe that $\kappa(x,y,z)$ in~\eqref{def.kappa_rhos} is the sum of the first two terms of $\psi(x,y,z)$ in~\eqref{eq.psi}, so all that remains is to consider the final terms in the definitions of $\bar\rho(x,y,z)$ and $\psi(x,y,z)$. If $c_i(x) \geq 0$, then $|\min\{z_i,-c_i(x)\}| = c_i(x)$ and the terms match. If $c_i(x) < 0$, then $z_i$ and $-c_i(x)$ are both nonnegative, and since the minimum of two nonnegative numbers is at most their geometric mean, it follows that $|\min\{z_i,-c_i(x)\}| \leq (-c_i(x)z_i)^{1/2}$.  Overall, the first inequality in \eqref{eq.rhobar_lower} follows. The second inequality in \eqref{eq.rhobar_lower} holds since $(y_x,z_x)$ is feasible for~\eqref{eq.def_omega}, and the equality holds by the definition of $(\yhat,\zhat)$.  For the last inequality in~\eqref{eq.rhobar_lower}, consider Theorem~\ref{th.OberWrig2.1} at $(x,\yhat,\zhat)$, the conditions of which hold since $\zhat \geq 0$ and
  \bequationNN
    \dist(\{(x,\yhat,\zhat)\},\{x_*\} \times \Dcal(x_*)) \leq \|x - x_*\|_2 + \dist(\{(\yhat,\zhat)\},\Dcal(x_*)) \leq \epsilon' + \epsilon/2 \leq \epsilon.
  \eequationNN
  Theorem~\ref{th.OberWrig2.1} then gives $\|x - x_*\|_2 \leq \dist(\{(x,\yhat,\zhat)\},\{x_*\} \times \Dcal(x_*)) \leq C\psi(x,\yhat,\zhat)$, where the first inequality holds since every element of $\{x_*\} \times \Dcal(x_*)$ has $x_*$ as its first block, which gives the final inequality in~\eqref{eq.rhobar_lower}.  (As an aside, note that Theorem~\ref{th.OberWrig2.1} cannot be applied directly at $(x,y_x,z_x)$, since nothing ensures that a solution of \eqref{prob.lp-lpec}, which is computed from approximate values, lies within distance $\epsilon$ of $\{x_*\} \times \Dcal(x_*)$. This is why the bound is obtained here through $(\yhat,\zhat)$ attaining the minimum in \eqref{eq.def_omega}.)

  Our aim now is to introduce conditions that define the thresholds $\hat\epsilon$ and $\bar\epsilon$ such that the desired conclusions of the theorem hold.  Toward this end, let $L_c \in \R{}_{>0}$ be a Lipschitz constant for $c$ on $\{x \in \R{n} : \|x - x_*\|_2 \leq \epsilon'\}$, let $m := \min_{i \notin \Acal(x_*)} |c_i(x_*)| > 0$ (with $m := \infty$ if $\Acal(x_*) = [q]$), and define
  \bequationNN
    \tau := \min\left\{ 2, \tfrac{1}{4\beta}\(\tfrac{m}{4}\)^{1/\sigma} \right\},\ \ t_\tau := \min\left\{\tfrac{\tau}{2},\ \tfrac{\tau^2}{4q}\right\},\ \ \text{and}\ \ \bar\delta := \( \tfrac{(\beta/(2C))^\sigma}{2\max\{L_c,(2Cq)^{-1}\}} \)^{1/(1-\sigma)}.
  \eequationNN
  Fix $(y_*,z_*) \in \Dcal(x_*)$ and choose $\hat\epsilon \in (0,\min\{\epsilon',\bar\delta\}]$ and $\bar\epsilon \in \R{}_{>0}$ sufficiently small that, for all $x$ with $\|x - x_*\|_2 \leq \hat\epsilon$ and all approximations with $\|\epsilonerror\|_\infty \leq \bar\epsilon$, one has (i)~$\|c(x)\|_\infty \leq \Mbar$; (ii)~$c_i(x) \leq \thalf c_i(x_*) < 0$ for all $i \notin \Acal(x_*)$; (iii)~$\rho(x,y_*,z_*) \leq \thalf t_\tau$ and $\epsilon_\rho \leq \tfrac14 t_\tau$; and (iv)~$\epsilon_c \leq \tfrac14 m$.  Conditions (i), (ii), and the first part of (iii) can be ensured by reducing $\hat\epsilon$ toward zero, as needed. In particular, with respect to (i) and~(ii), this can be done by continuity of $c$, and with respect to the first part of (iii) it follows since $\rho(\cdot,y_*,z_*)$ is continuous (the functions $\phi_i$ for all $i \in [q]$ in the proof of Lemma~\ref{lem.basic_lemma_2} are continuous) and $\rho(x_*,y_*,z_*) = 0$ by $(y_*,z_*) \in \Dcal(x_*)$ and \eqref{eq.kkt}.  The second part of (iii) and (iv), on the other hand, can be ensured by reducing~$\bar\epsilon$ as needed since $\epsilon_\rho$ is a linear function of $\epsilonerror$ according to \eqref{eq.eps_rho_defs} and $\epsilon_c \leq \|\epsilonerror\|_\infty$.

  Let us now derive upper bounds that hold for all such $x$ and all such approximations (satisfying (i)--(iv)).  The point $(y_*,z_*)$ is feasible for \eqref{prob.lp-lpec} since $\|(y_*,z_*)\|_\infty \leq M - 1$, and both $(x,y_x,z_x)$ and $(x,y_*,z_*)$ satisfy the conditions of Lemma~\ref{lem.basic_lemma_2} by (i) and the box constraints in \eqref{prob.lp-lpec}.  Hence, by the former inequality in~\eqref{eq.rho_pert}, the optimality of $(y_x,z_x)$ with respect to \eqref{prob.lp-lpec}, and (iii), it follows that
  \bequation\label{eq.stumpy}
    \rho(x,y_x,z_x) - \epsilon_\rho \leq \tilde\rho(x,y_x,z_x) \leq \tilde\rho(x,y_*,z_*) \leq \rho(x,y_*,z_*) + \epsilon_\rho \leq \tfrac34 t_\tau.
  \eequation
  In particular, $\rho(x,y_x,z_x) \leq t_\tau \leq \thalf\tau \leq 1$.  (This fact is used later in Corollary~\ref{cor.rate}.)  Moreover, by Lemma~\ref{lem.basic_lemma_1} and since $t \mapsto t + (qt)^{1/2}$ is increasing, \eqref{eq.stumpy} gives $\tilde{\bar\rho}(x,y_x,z_x) \leq \tilde{\rho}(x,y_x,z_x) + q^{1/2} \tilde{\rho}(x,y_x,z_x)^{1/2} \leq t_\tau + (qt_\tau)^{1/2} \leq \thalf\tau + \thalf\tau = \tau$, so that, by the definition of $\tau$, one obtains
  \bequation\label{eq.thm41_eps3}
    \(\beta\tilde{\bar\rho}(x,y_x,z_x)\)^\sigma \leq (\beta\tau)^\sigma \leq 4^{-\sigma}\tfrac14 m < \tfrac14 m.
  \eequation

  Let us now prove part (a) of the theorem.  By (i), the box constraints in \eqref{prob.lp-lpec}, and $\|(\yhat,\zhat)\|_\infty \leq M \leq \Mbar$ (proved earlier), Lemma~\ref{lem.basic_lemma_2} applies at both $(x,y_x,z_x)$ and $(x,\yhat,\zhat)$.  The upper bound in part~(a) follows from \eqref{eq.def_omega} and \eqref{eq.ine_rhos} at $(x,y_x,z_x)$, namely, $\omega(x) \leq \psi(x,y_x,z_x) \leq \tilde{\bar\rho}(x,y_x,z_x) + \epsilon_\psi$.  For the lower bound, \eqref{eq.ine_rhos} at $(x,\yhat,\zhat)$ and the optimality of $(y_x,z_x)$ for \eqref{prob.lp-lpec}, for which $(\yhat,\zhat)$ is feasible, give
  \bequationNN
    \omega(x) = \psi(x,\yhat,\zhat) \geq \Mbar^{-1} (\tilde\rho(x,\yhat,\zhat) - \epsilon_\rho) \geq \Mbar^{-1}(\tilde\rho(x,y_x,z_x) - \epsilon_\rho),
  \eequationNN
  and the last quantity is positive by the conditions of part~(a) and the final part of Lemma~\ref{lem.basic_lemma_2} at $(x,y_x,z_x)$.

  Finally, let us prove part (b). Let $i \notin \Acal(x_*)$.  By \eqref{eq.error_bounds}, (ii), (iv), and \eqref{eq.thm41_eps3}, it follows that
  \bequationNN
    \capprox_i(x) \leq c_i(x) + \epsilon_c \leq \thalf c_i(x_*) + \epsilon_c \leq -\thalf m + \tfrac14 m = -\tfrac14 m < -(\beta\tilde{\bar\rho}(x,y_x,z_x))^\sigma,
  \eequationNN
  so $i \notin \widetilde\Acal_{{\rm LP}}(y_x,z_x;x)$ by \eqref{eq.active_LP}.  Now let $i \in \Acal(x_*)$.  By \eqref{eq.rho_pert}--\eqref{eq.psitilde_rhobar}, \eqref{eq.rhobar_lower}, and $\epsilon_\psi \leq (2C)^{-1}\|x - x_*\|_2$,
  \bequationNN
    \tilde{\bar\rho}(x,y_x,z_x) \geq \tilde\psi(x,y_x,z_x) \geq \psi(x,y_x,z_x) - \epsilon_\psi \geq C^{-1}\|x - x_*\|_2 - \epsilon_\psi \geq (2C)^{-1}\|x - x_*\|_2.
  \eequationNN
  Moreover, $\epsilon_c \leq q^{-1}\epsilon_\psi \leq (2Cq)^{-1}\|x - x_*\|_2$ and $|c_i(x)| = |c_i(x) - c_i(x_*)| \leq L_c\|x - x_*\|_2$.  Since $\|x - x_*\|_2 \leq \hat\epsilon \leq \bar\delta$, one has $\max\{L_c,(2Cq)^{-1}\}\|x - x_*\|_2^{1-\sigma} \leq \thalf(\beta/(2C))^\sigma$, and therefore
  \bequationNN
    |\capprox_i(x)| \leq |c_i(x)| + \epsilon_c \leq \(L_c + (2Cq)^{-1}\)\|x - x_*\|_2 \leq \(\tfrac{\beta\|x - x_*\|_2}{2C}\)^\sigma \leq (\beta\tilde{\bar\rho}(x,y_x,z_x))^\sigma,
  \eequationNN
  so $\capprox_i(x) \geq -(\beta\tilde{\bar\rho}(x,y_x,z_x))^\sigma$ and $i \in \widetilde\Acal_{{\rm LP}}(y_x,z_x;x)$ by \eqref{eq.active_LP}.  Overall, $\widetilde\Acal_{{\rm LP}}(y_x,z_x;x) = \Acal(x_*)$.
  \qed
\eproof
}

\edit{
The condition in part~(a) of Theorem~\ref{the.LP_main_theorem} involves $(y_x,z_x)$, meaning that it depends on the realized noise. On the other hand, the condition in part~(b) depends on $x$ only through $\|x - x_*\|_2$.  The following corollary gives sufficient conditions of the latter form for both parts of the theorem.
}

{\editcolor
\bcorollary\label{cor.rate}
  Suppose that the conditions of Theorem~\ref{the.LP_main_theorem} hold, which define $C$, $\hat\epsilon$, and $\bar\epsilon$.  Then, for every $x \in \R{n}$ with $\|x - x_*\|_2 \leq \hat\epsilon$ and any approximations satisfying \eqref{eq.error_bounds} with $\|\epsilonerror\|_\infty \leq \bar\epsilon$, it follows that
  \bequation\label{eq.rho_lower}
    \bar\rho(x,y_x,z_x) \geq \tfrac{\|x-x_*\|_2}{C}\ \ \text{and}\ \ \rho(x,y_x,z_x) \geq \( \tfrac{\|x-x_*\|_2}{C(1+q^{1/2})} \)^2.
  \eequation
  Also, there exists $(\Lambda',\Lambda) \in \R{}_{>0} \times \R{}_{>0}$, dependent only on $(C,\Mbar,n,p,q,\hat\epsilon,\bar\epsilon)$, such that $\Lambda'\hat\epsilon \leq \bar\epsilon$, $\Lambda\hat\epsilon^2 \leq \bar\epsilon$, and for each $x$ with $0 < \|x - x_*\|_2 \leq \hat\epsilon$ and any approximations satisfying \eqref{eq.error_bounds}, Theorem~\ref{the.LP_main_theorem}(b) holds when
  \bequation\label{eq.rate_lin}
    \|\epsilonerror\|_\infty \leq \Lambda'\|x - x_*\|_2,
  \eequation
  and Theorem~\ref{the.LP_main_theorem}(a)--(b) hold when
  \bequation\label{eq.rate}
    \|\epsilonerror\|_\infty \leq \Lambda\|x - x_*\|_2^2.
  \eequation
\ecorollary
\bproof{Proof.}
  The first bound in \eqref{eq.rho_lower} follows from \eqref{eq.rhobar_lower}, which was established in the proof of Theorem~\ref{the.LP_main_theorem}.  For the second bound in \eqref{eq.rho_lower}, recall from the same proof that $\rho(x,y_x,z_x) \leq 1$ for every~$x$ with $\|x - x_*\|_2 \leq \hat\epsilon$ and any approximations with $\|\epsilonerror\|_\infty \leq \bar\epsilon$.  (See the comment after \eqref{eq.stumpy}.)  Hence, by \cite[Lemma~3.5]{OberWrig2006},
  \bequationNN
    \bar\rho(x,y_x,z_x) \leq \rho(x,y_x,z_x) + q^{1/2}\rho(x,y_x,z_x)^{1/2} \leq (1+q^{1/2})\rho(x,y_x,z_x)^{1/2},
  \eequationNN
  so $\rho(x,y_x,z_x) \geq \bar\rho(x,y_x,z_x)^2/(1+q^{1/2})^2$. Combining this with the first bound in \eqref{eq.rho_lower} yields the second.

  It remains to establish the existence of $(\Lambda',\Lambda)$.  Consider arbitrary $x$ with $0 < \|x - x_*\|_2 \leq \hat\epsilon$ and define
  \bequationNN
    c_\psi := c_\kappa + q\ \ \text{and}\ \ c_\rho := c_\kappa + q\Mbar,\ \ \text{where}\ \ c_\kappa := p^{1/2} + n^{1/2}(1 + p^{1/2}\Mbar + q^{1/2}\Mbar),
  \eequationNN
  so that, by \eqref{eq.eps_rho_defs}, $\epsilon_\psi \leq c_\psi \|\epsilonerror\|_\infty$ and $\epsilon_\rho \leq c_\rho \|\epsilonerror\|_\infty$.  Let
  \bequationNN
    \Lambda' := \min\left\{ \tfrac{1}{2Cc_\psi},\ \tfrac{\bar\epsilon}{\hat\epsilon} \right\}
    \ \ \text{and}\ \
    \Lambda := \min\left\{ \tfrac{1}{4c_\rho C^2(1+q^{1/2})^2},\ \tfrac{\Lambda'}{\hat\epsilon} \right\},
  \eequationNN
  so that $\Lambda'\hat\epsilon \leq \bar\epsilon$ and $\Lambda\hat\epsilon^2 \leq \Lambda'\hat\epsilon \leq \bar\epsilon$.  Suppose first that \eqref{eq.rate_lin} holds.  Then $\|\epsilonerror\|_\infty \leq \Lambda'\hat\epsilon \leq \bar\epsilon$, so the conditions of Theorem~\ref{the.LP_main_theorem} hold and $\epsilon_\psi \leq c_\psi \|\epsilonerror\|_\infty \leq c_\psi\Lambda'\|x - x_*\|_2 \leq (2C)^{-1}\|x - x_*\|_2$. Hence Theorem~\ref{the.LP_main_theorem}(b) holds.  Suppose now that $\|\epsilonerror\|_\infty \leq \Lambda\|x - x_*\|_2^2$.  Then $\|\epsilonerror\|_\infty \leq (\Lambda'/\hat\epsilon)\|x - x_*\|_2^2 \leq \Lambda'\|x - x_*\|_2$ since $\|x - x_*\|_2 \leq \hat\epsilon$, so Theorem~\ref{the.LP_main_theorem}(b) holds as shown.  Moreover, the bounds in \eqref{eq.rho_lower} hold since $\|\epsilonerror\|_\infty \leq \bar\epsilon$, so
  \bequationNN
    \epsilon_\rho \leq c_\rho \|\epsilonerror\|_\infty \leq c_\rho\Lambda\|x - x_*\|_2^2 \leq \tfrac{\|x - x_*\|_2^2}{4C^2(1+q^{1/2})^2} \leq \tfrac14\rho(x,y_x,z_x) < \thalf\rho(x,y_x,z_x),
  \eequationNN
  which is the condition in Theorem~\ref{the.LP_main_theorem}(a).  Hence, Theorem~\ref{the.LP_main_theorem}(a)--(b) hold, as claimed.
  \qed
\eproof

Corollary~\ref{cor.rate} shows that correct identification requires the noise to vanish only linearly in $\|x - x_*\|_2$, whereas the two-sided bound on $\omega(x)$ in Theorem~\ref{the.LP_main_theorem}(a) requires a quadratic rate.  The latter is due to the square-root terms in $\bar\rho$, which make $\rho$ of order $\bar\rho^2$ near $x_*$.  In upcoming Theorem~\ref{th.impossible}(e), it is shown that the linear rate in \eqref{eq.rate_lin} cannot be improved in general. That said, we do not know whether the quadratic rate for ensuring that part~(a) holds is tight.  Since \eqref{eq.rate_lin} depends on $x$ only through $\|x - x_*\|_2$, it is the form of the noise requirement that we use for the setting of a stochastic algorithm in \S\ref{sec.stochastic_algorithm}.
}

\subsection{Identification through a Primal Step}\label{sec.primal_step}

We now present a technique for active-set identification in noisy settings that is based on computing a primal step through solving either an LP or a QP.  The QP approach is an extension of the one that offers Theorem~\ref{th.QP_exact}, and as previously stated a similar result holds when one solves an LP with a trust-region constraint.  Our primary tool for extending these approaches to noisy settings is perturbation theory for the solution of linear systems, in particular, \edit{Lemmas~\ref{lem.SysPertMat}--\ref{lem.SysPertAbs}}.

For the sake of brevity in this section, we state only the QP-based approach in detail.  All of the problems statements and results hold similarly for an LP-based approach with a trust-region constraint.

Analogously to \eqref{eq.sub}, consider the QP subproblem
\bequation\label{eq.sub_approx_nu}
  \baligned
    \min_{(d,r,s,t) \in \R{n} \times \R{p} \times \R{p} \times \R{q}} &\ \fgradapprox(x)^Td + \nu(\ones^Tr + \ones^Ts + \ones^Tt) + \thalf \theta \|d\|_2^2 \\
    \st &\ \eapprox(x) + \egradapprox(x)^Td = r - s,\ \ \capprox(x) + \cgradapprox(x)^Td \leq t, \\
    &\ r \geq 0,\ \ s \geq 0,\ \ \text{and}\ \ t \geq 0.
  \ealigned
\eequation
Like~\eqref{eq.sub}, this QP is always feasible and has a finite optimal value.  \edit{Later in our analysis we use a fact from exact penalty theory \cite[Theorem~4.4]{HanMang1979}, namely, that if \eqref{eq.sub_approx} below is feasible and its (unique) solution admits multipliers $(\alpha,\beta)$, then for every $\nu \geq \|(\alpha,\beta)\|_\infty$ the solution of \eqref{eq.sub_approx_nu} has $(r_{x,\theta,\nu},s_{x,\theta,\nu},t_{x,\theta,\nu}) = (0,0,0)$ and its $d$-component coincides with the (unique) solution $d_{x,\theta,\nu}$ of}
\bequation\label{eq.sub_approx}
  \min_{d \in \R{n}}\ \fgradapprox(x)^Td + \thalf \theta \|d\|_2^2\ \ \st\ \ \eapprox(x) + \egradapprox(x)^Td = 0\ \ \text{and}\ \ \capprox(x) + \cgradapprox(x)^Td \leq 0;
\eequation
\edit{that is, under these conditions the auxiliary variables in \eqref{eq.sub_approx_nu} all vanish at the solution.}

Our first lemma is based on the observation that if the LICQ holds at a \edit{primal KKT point} $x_*$ of~\eqref{eq.opt}, a point~$x$ is sufficiently close to~$x_*$, and the noise in certain constraint derivative values at~$x$ is sufficiently small, then the corresponding constraint derivative matrix has full rank.

\blemma\label{lem.invertibility}
  Suppose that, with respect to~\eqref{eq.opt}, the LICQ holds at a \edit{primal KKT point} $x_* \in \R{n}$.  Let the active set at $x_*$ be $\Acal_* := \Acal(x_*)$ and let its cardinality be $\ell := |\Acal_*|$.  Then, there exist $\epsilon_x \in \R{}_{>0}$, $\epsilon_{\nabla e} \in \R{}_{>0}$, and $\epsilon_{\nabla c} \in \R{}_{>0}$ such that any $(x, \egradapprox(x), \cgradapprox_{\Acal_*}(x)) \in \R{n} \times \R{n \times p} \times \R{n \times \ell}$ with
  \bequation\label{eq.derivative_errors1}
    \|x - x_*\|_2 \leq \epsilon_x,\ \|\egradapprox(x) - \nabla e(x)\|_2 \leq \epsilon_{\nabla e},\ \ \text{and}\ \ \|\cgradapprox_{\Acal_*}(x) - \nabla c_{\Acal_*}(x)\|_2 \leq \epsilon_{\nabla c}
  \eequation
  yields that the matrix $\bbmatrix \egradapprox(x) & \cgradapprox_{\Acal_*}(x) \ebmatrix$ has full column rank.
\elemma
\bproof{Proof.}
  Under the conditions of the lemma, specifically that the LICQ for~\eqref{eq.opt} holds at $x_*$, the matrix $\bbmatrix \nabla e(x_*) & \nabla c_{\Acal_*}(x_*) \ebmatrix$ has full column rank.  Then, since $\{e_i\}_{i\in[p]}$ and $\{c_i\}_{i\in\Acal_*}$ are continuously differentiable, it follows that $\bbmatrix \nabla e(x) & \nabla c_{\Acal_*}(x) \ebmatrix$ is continuous as a function of $x$.  Since the set of full-column-rank matrices is open, there exists $\epsilon_x \in \R{}_{>0}$ such that for all $x \in \R{n}$ with $\|x - x_*\|_2 \leq \epsilon_x$, the matrix $\bbmatrix \nabla e(x) & \nabla c_{\Acal_*}(x) \ebmatrix$ has full column rank.  Consequently, for any such $\epsilon_x$, there exist sufficiently small $\epsilon_{\nabla e} \in \R{}_{>0}$ and $\epsilon_{\nabla c} \in \R{}_{>0}$ such that $\bbmatrix \egradapprox(x) & \cgradapprox_{\Acal_*}(x) \ebmatrix$ has full column rank as well, which completes the proof.
  \qed
\eproof

\edit{Lemma~\ref{lem.invertibility} is consistent with a well-known phenomenon, namely, that constraint qualifications are stable under small perturbations of the problem data.  Such results go back at least to \cite{Robi1976}, and they appear, in far greater generality, as special cases of the stability theory for convex composite optimization (see, e.g., \cite[Lemma~3.2]{BurkFerr1995}) and of the characterizations of strong regularity for variational inequalities over polyhedral convex sets in \cite{DontRock1996}.  We include the short proof of Lemma~\ref{lem.invertibility} for the sake of self-containedness.  We also note that the same argument applies to the MFCQ.  That is, the MFCQ at $x_*$ is equivalent to positive linear independence of $\{\nabla e_i(x_*)\}_{i \in [p]} \cup \{\nabla c_i(x_*)\}_{i \in \Acal_*}$, which is an open condition, so under the MFCQ the conclusion of Lemma~\ref{lem.invertibility} holds with ``full column rank'' replaced by ``positive linear independence.''}

Our next lemma shows that if $x$ is close to a \edit{primal KKT point} $x_*$ and the noise is sufficiently small, then the QP~\eqref{eq.sub_approx} is feasible and its primal-dual solution is given by the solution of a linear system defined with respect to the active set at $x_*$.  A similar result would hold for an LP subproblem with a trust-region constraint, since at a \edit{primal KKT point} the solution of the subproblem would have no activity with respect to the trust region, and in a neighborhood of the \edit{primal KKT point} the subproblem solutions would remain inactive with respect to the trust region.  For the LP case, \eqref{eq.PerSys} would have a zero matrix in the (1,1)-block.

\blemma\label{lem.keynoisylemma}
  Suppose that, with respect to~\eqref{eq.opt}, the LICQ and strict complementarity hold at \edit{a primal KKT point} $x_* \in \R{n}$.  Let $\Acal_* := \Acal(x_*)$ and $\ell := |\Acal_*|$.  Then, for any $\theta \in \R{}_{>0}$, there exist $\epsilon_x \in \R{}_{>0}$, $\epsilon_e \in \R{}_{>0}$, $\epsilon_c \in \R{}_{>0}$, $\epsilon_{\nabla f} \in \R{}_{>0}$, $\epsilon_{\nabla e} \in \R{}_{>0}$, and $\epsilon_{\nabla c} \in \R{}_{>0}$ such that for all $(x,\eapprox(x),\capprox(x),\fgradapprox(x),\egradapprox(x),\cgradapprox(x))$ satisfying $\|x - x_*\|_2 \leq \epsilon_x$ and~\eqref{eq.error_bounds} the linear system
  \bequation\label{eq.PerSys}
    \bbmatrix \theta I & \egradapprox(x) & \cgradapprox_{\Acal_*}(x) \\ \egradapprox(x)^T & 0 & 0 \\ \cgradapprox_{\Acal_*}(x)^T & 0 & 0 \ebmatrix \bbmatrix d \\ \alpha \\ \beta_{\Acal_*} \ebmatrix = -\bbmatrix \fgradapprox(x) \\ \eapprox(x) \\ \capprox_{\Acal_*}(x) \ebmatrix
  \eequation
  has a unique solution $(d, \alpha, \beta_{\Acal_*}) \in \R{n} \times\R{p} \times \R{\ell}$ which, along with $\beta_{\Acal_*^c} = 0$, satisfies the KKT conditions of subproblem~\eqref{eq.sub_approx}; in particular, one has
  \bsubequations\label{eq.sub_approx_KKT}
    \begin{align}
      \fgradapprox(x) + \theta d + \egradapprox(x) \alpha + \edit{\cgradapprox(x)} \beta &= 0, \label{eq.sub_approx_KKT_1} \\
      \eapprox(x) + \egradapprox(x)^T d = 0,\ \ \capprox_{\Acal_*}(x) + \cgradapprox_{\Acal_*}(x)^T d &= 0,\label{eq.sub_approx_KKT_2} \\
      \capprox_{\Acal_*^c}(x) + \cgradapprox_{\Acal_*^c}(x)^T d &< 0, \label{eq.sub_approx_KKT_3} \\
      \beta_{\Acal_*} > 0,\ \ \beta_{\Acal_*^c} &= 0, \label{eq.sub_approx_KKT_4} \\
      \text{and}\ \ \capprox(x) + \cgradapprox(x)^T d &\perp \beta. \label{eq.sub_approx_KKT_5}
    \end{align}
  \esubequations
\elemma
\bproof{Proof.}
  Under the stated conditions, the conditions of Lemma~\ref{lem.invertibility} hold, which implies that there exist $\epsilon_x \in \R{}_{>0}$, $\epsilon_{\nabla e} \in \R{}_{>0}$, and $\epsilon_{\nabla c} \in \R{}_{>0}$ such that the matrix in~\eqref{eq.PerSys}, call it $\Ktilde$, is invertible and the linear system in~\eqref{eq.PerSys} has a unique solution.  Let such $(\epsilon_x, \epsilon_{\nabla e}, \epsilon_{\nabla c})$ be given.  Now note that, by construction of~\eqref{eq.PerSys} and with $\beta_{\Acal_*^c} = 0$, the conditions  \eqref{eq.sub_approx_KKT_1}, \eqref{eq.sub_approx_KKT_2}, and \eqref{eq.sub_approx_KKT_5} are all satisfied.  All that remains is to show that, when $x$ is sufficiently close to $x_*$ and all of the noise in the function and derivative values is sufficiently small, the conditions \eqref{eq.sub_approx_KKT_3} and \eqref{eq.sub_approx_KKT_4} hold as well.  Toward this end, let us define
  \bequationNN 
    K := \bbmatrix \theta I & \nabla e(x_*) & \nabla c_{\Acal_*}(x_*) \\ \nabla e(x_*)^T & 0 & 0 \\ \nabla c_{\Acal_*}(x_*)^T & 0 & 0 \ebmatrix\ \text{and}\ \bbmatrix d_* \\ \alpha_* \\ [\beta_{\Acal_*}]_{*} \ebmatrix = - K^{-1} \bbmatrix \nabla f(x_*) \\ e(x_*) \\ c_{\Acal_*}(x_*) \ebmatrix.
  \eequationNN
  Under the conditions of the lemma, it follows by Farkas' theorem that $d_* = 0$ and, along with $[\beta_{\Acal_*^c}]_* = 0$, the pair $(\alpha_*,\beta_*) =: (y_*,z_*)$ is the unique Lagrange multiplier such that $(x_*,y_*,z_*)$ is a KKT point for problem~\eqref{eq.opt}. Indeed, by the strict complementarity assumption in the lemma, it follows that $[\beta_{\Acal_*}]_{*} > 0$. Now note that the matrix~$K$ is invertible (e.g., by setting $x = x_*$ in Lemma~\ref{lem.invertibility}), and that one has $\Ktilde = K + \delta K$, where---since $f$, $e$, and $c$ are continuously differentiable---one finds that $\delta K \to 0$ as $x \to x_*$ and $\epsilonerror \to 0$, where $\epsilonerror$ is defined in~\eqref{eq.epsilon_error}.  Consequently, $\|\delta K\|_2 \|K^{-1}\|_2 \to 0$ as $x \to x_*$ and $\epsilonerror \to 0$, which means that, for $x$ sufficiently close to $x_*$ and $\epsilonerror$ sufficiently small, one has from Lemma \ref{lem.SysPertMat} that \eqref{eq.PerSys} yields
  \bequationNN 
    \tfrac{\left\| \bbmatrix d_* - d \\ \alpha_* - \alpha \\ [\beta_{\Acal_*}]_* - \beta_{\Acal_*} \ebmatrix \right\|_2}{\left\| \bbmatrix d_* \\ \alpha_* \\ [\beta_{\Acal_*}]_* \ebmatrix \right\|_2} \leq \tfrac{\cond(K)}{1 - \|\delta K\|_2 \|K^{-1}\|_2} \left( \tfrac{\|\delta K\|_2}{\|K\|_2} + \tfrac{\left\|\bbmatrix \fgradapprox(x) \\ \eapprox(x) \\ \capprox_{\Acal_*}(x) \ebmatrix - \bbmatrix \nabla f(x_*) \\ e(x_*) \\ c_{\Acal_*}(x_*) \ebmatrix \right\|_2}{\left\| \bbmatrix \nabla f(x_*) \\ e(x_*) \\ c_{\Acal_*}(x_*) \ebmatrix \right\|_2} \right).
  \eequationNN
  Overall, since $d_* = 0$ and $[\beta_{\Acal_*}]_{*} > 0$, it follows that for $\epsilon_x \in \R{}_{>0}$, $\epsilon_e \in \R{}_{>0}$, $\epsilon_c \in \R{}_{>0}$, $\epsilon_{\nabla f} \in \R{}_{>0}$, $\epsilon_{\nabla e} \in \R{}_{>0}$, and $\epsilon_{\nabla c} \in \R{}_{>0}$ all sufficiently small, any $(x,\eapprox(x),\capprox(x),\fgradapprox(x),\egradapprox(x),\cgradapprox(x))$ satisfying $\|x - x_*\|_2 \leq \epsilon_x$ and the error bounds~\eqref{eq.error_bounds} yield that all conditions in \eqref{eq.sub_approx_KKT} hold, as desired.
  \qed
\eproof

{\editcolor
The tolerances asserted by Lemma~\ref{lem.keynoisylemma} can be made explicit in terms of the conditioning of the matrix~$K$ defined in its proof. This is shown next.  In what follows, using the notation of Lemma~\ref{lem.keynoisylemma} for such a primal KKT point~$x_*$, let us define $A_* := \displaystyle \bbmatrix \nabla e(x_*) & \nabla c_{\Acal_*}(x_*) \ebmatrix \in \R{n \times (p+\ell)}$, $a_* := \|A_*\|_2$, $\sigma_* := \sigma_{\min}(A_*) > 0$ (where positivity follows under the LICQ), $s_* := \min_{i \notin \Acal_*}(-c_i(x_*)) > 0$, and $\zeta_* := \min_{i \in \Acal_*}[z_*]_i > 0$ (where positivity follows by strict complementarity).  Furthermore, suppose that $\nabla f$, $\nabla e$, $\nabla c$, $e$, and $c$ are Lipschitz continuous over $\{x \in \R{n} : \|x - x_*\|_2 \leq 1\}$ with $L \geq 1$ being a common Lipschitz constant for these functions over this set, and let $G := \max_{i \in [q]} \max_{x \in \R{n}} \{\|\nabla c_i(x)\|_2 : \|x - x_*\|_2 \leq 1\}$.

\blemma\label{lem.explicit}
  Suppose that the conditions of Lemma~\ref{lem.keynoisylemma} hold and consider arbitrary $\theta \in \R{}_{>0}$.  Then, the matrix $K$ defined in the proof of Lemma~\ref{lem.keynoisylemma} satisfies
  \bequation\label{eq.condK}
    \baligned
      \|K\|_2 &= \tfrac{\theta + \sqrt{\theta^2 + 4a_*^2}}{2} \leq \theta + a_* \\
      \text{and}\ \ \|K^{-1}\|_2 &\leq \kappa_\theta := \max\left\{ \tfrac{1}{\theta},\ \tfrac{\theta + \sqrt{\theta^2+4\sigma_*^2}}{2\sigma_*^2} \right\} \leq \tfrac{1}{\theta} + \tfrac{1}{\sigma_*} + \tfrac{\theta}{\sigma_*^2},
    \ealigned
  \eequation
  from which it follows that $\cond(K) \leq (\theta + a_*)\kappa_\theta$.  Moreover, defining $\delta_K := \sqrt2 (L\epsilon_x + \epsilon_{\nabla e} + \epsilon_{\nabla c})$ and $\delta_b := \epsilon_{\nabla f} + \epsilon_e + \epsilon_c + 3L\epsilon_x$, it follows that if $\epsilon_x \leq 1$,
  \bequation\label{eq.explicit_cond}
    \kappa_\theta \delta_K \leq \thalf,\ \ \text{and}\ \ \Theta := 2\kappa_\theta\( \delta_b + \delta_K \|(y_*,[z_*]_{\Acal_*})\|_2 \)
    < \min\left\{ 1, \zeta_*,\ \tfrac{s_* - L\epsilon_x - \epsilon_{\nabla c} - \epsilon_c}{G} \right\},
  \eequation
  then the conclusions of Lemma~\ref{lem.keynoisylemma} hold; i.e., \eqref{eq.PerSys} has a unique solution $(d,\alpha,\beta_{\Acal_*})$ and, with $\beta_{\Acal_*^c} = 0$, the conditions \eqref{eq.sub_approx_KKT} are all satisfied.  In particular, \eqref{eq.explicit_cond} holds whenever $\epsilon_x$ and $\|\epsilonerror\|_\infty$ are bounded above by a constant determined by $\theta$, $\sigma_*$, $a_*$, $L$, $G$, $\zeta_*$, $s_*$, and $\|(y_*,z_*)\|_2$.
\elemma
\bproof{Proof.}
  Let $r := p + \ell$, let $A_* = U\Sigma V^T$ be a (thin) singular value decomposition of $A_*$ with $U \in \R{n \times r}$, $\Sigma = \diag(\sigma_1,\dots,\sigma_r)$, and orthogonal $V \in \R{r \times r}$, and let $U_\perp \in \R{n \times (n-r)}$ complete $U$ to give an orthonormal basis of $\R{n}$.  With the orthogonal matrix $Q := \diag(\bbmatrix U & U_\perp \ebmatrix, V)$ one finds
  \bequationNN
    Q^TKQ = \bbmatrix \theta I_r & 0 & \Sigma \\ 0 & \theta I_{n-r} & 0 \\ \Sigma & 0 & 0 \ebmatrix\ \ \text{with identity matrices $I_r \in \R{r \times r}$ and $I_{n-r} \in \R{(n-r) \times (n-r)}$},
  \eequationNN
  which, after a symmetric permutation, is the sum of the $2 \times 2$ blocks $\bbmatrix \theta & \sigma_i \\ \sigma_i & 0 \ebmatrix$ for $i \in [r]$ and $\theta I_{n-r}$.  Thus, the eigenvalues of $K$ are $\thalf(\theta \pm (\theta^2 + 4\sigma_i^2)^{1/2})$ for $i \in [r]$, together with $\theta$ with multiplicity $n-r$.  Since $K$ is symmetric, $\|K\|_2$ is the largest of these in absolute value, which is $\thalf(\theta + (\theta^2+4a_*^2)^{1/2})$ because $a_* = \max_i \sigma_i$ and $\thalf(\theta + (\theta^2+4\sigma^2)^{1/2}) \geq \theta$ for every $\sigma \geq 0$. The inequality $\thalf(\theta + (\theta^2+4a_*^2)^{1/2}) \leq \theta + a_*$ in \eqref{eq.condK} follows from $(\theta^2+4a_*^2)^{1/2} \leq \theta + 2a_*$.  On the other hand, $\|K^{-1}\|_2$ is the reciprocal of the smallest eigenvalue of $K$ in absolute value.  That smallest value is $\min\{\theta, \thalf((\theta^2+4\sigma_*^2)^{1/2} - \theta)\}$, where the second entry is attained at $\sigma_* = \min_i \sigma_i$ since $\sigma \mapsto (\theta^2+4\sigma^2)^{1/2} - \theta$ is increasing.  Overall, $2/((\theta^2+4\sigma_*^2)^{1/2} - \theta) = (\theta + (\theta^2+4\sigma_*^2)^{1/2})/(2\sigma_*^2)$, which gives $\kappa_\theta$ in \eqref{eq.condK} with equality when $r<n$ and with inequality when $r=n$, since in that case the eigenvalue $\theta$ is absent. The end of \eqref{eq.condK} follows from $(\theta^2+4\sigma_*^2)^{1/2} \leq \theta + 2\sigma_*$.

  Now let $w := (d,\alpha,\beta_{\Acal_*})$ solve \eqref{eq.PerSys} and let $w_* := (0,y_*,[z_*]_{\Acal_*})$, which solves $Kw_* = -(\nabla f(x_*),e(x_*),c_{\Acal_*}(x_*)) = -(\nabla f(x_*),0,0)$ as established in the proof of Lemma~\ref{lem.keynoisylemma}.  The matrix $\Ktilde$ of \eqref{eq.PerSys} satisfies $\Ktilde = K + \delta K$, where $\delta K$ has zero diagonal blocks and off-diagonal blocks $E$ and $E^T$ with $E := \bbmatrix \egradapprox(x) & \cgradapprox_{\Acal_*}(x) \ebmatrix - A_*$, so that $\|\delta K\|_2 = \|E\|_2$.  Writing $E = \bbmatrix P & Q \ebmatrix$ with $P := \egradapprox(x) - \nabla e(x_*)$ and $Q := \cgradapprox_{\Acal_*}(x) - \nabla c_{\Acal_*}(x_*)$, one has $\|E\|_2^2 = \|PP^T + QQ^T\|_2 \leq \|P\|_2^2 + \|Q\|_2^2$, while $\|P\|_2 \leq L\epsilon_x + \epsilon_{\nabla e}$ and $\|Q\|_2 \leq L\epsilon_x + \epsilon_{\nabla c}$ by \eqref{eq.error_bounds} and the definition of $L$. Hence, $\|\delta K\|_2 \leq \sqrt2(L\epsilon_x + \epsilon_{\nabla e} + \epsilon_{\nabla c}) = \delta_K$.  Similarly the right-hand sides differ by at most $\delta_b$.  By the first condition in \eqref{eq.explicit_cond} and Lemma~\ref{lem.SysPertAbs}, the matrix $\Ktilde$ is nonsingular---so \eqref{eq.PerSys} has a unique solution---and $\|w - w_*\|_2 \leq \Theta$.  Conditions \eqref{eq.sub_approx_KKT_1}, \eqref{eq.sub_approx_KKT_2}, and \eqref{eq.sub_approx_KKT_5} hold by construction, as in the proof of Lemma~\ref{lem.keynoisylemma}.  For \eqref{eq.sub_approx_KKT_4}, one has $[\beta_{\Acal_*}]_i \geq [z_*]_i - \Theta \geq \zeta_* - \Theta > 0$ by \eqref{eq.explicit_cond}.  For \eqref{eq.sub_approx_KKT_3}, let $i \notin \Acal_*$. Then, $\|d\|_2 = \|d - 0\|_2 \leq \Theta < 1$ yields
  \bequationNN
    \capprox_i(x) + \cgradapprox_i(x)^Td \leq c_i(x_*) + L\epsilon_x + \epsilon_c + (G + \epsilon_{\nabla c})\Theta
    \leq -s_* + L\epsilon_x + \epsilon_{\nabla c} + \epsilon_c + G\Theta < 0
  \eequationNN
  by the second condition in \eqref{eq.explicit_cond}.  The final claim follows since $\delta_K$ and $\delta_b$ are nondecreasing linear functions of the tolerances $\epsilon_x$ and $\|\epsilonerror\|_\infty$ that vanish at the origin.
  \qed
\eproof

The quantity $\kappa_\theta$ in \eqref{eq.condK} is a Lipschitz constant of the solution map of the KKT system of \eqref{eq.sub_approx} with the active set fixed at $\Acal_*$, and the simpler bound $\kappa_\theta \leq \theta^{-1} + \sigma_*^{-1} + \theta\sigma_*^{-2}$ shows its dependence on the problem data. It degrades as the active constraint gradients approach linear dependence (i.e., as $\sigma_* \searrow 0$).  For fixed~$\sigma_*$, the precise value for $\kappa_\theta$ in \eqref{eq.condK} is minimized at $\theta = \sigma_*/\sqrt2$, which yields $\kappa_\theta = \sqrt2/\sigma_*$.  For such $\theta$, or any $\theta$ proportional to $\sigma_*$, condition \eqref{eq.explicit_cond} states that the admissible noise level is proportional to $\sigma_*\min\{1,\zeta_*,s_*/G\}/(1 + \|(y_*,z_*)\|_2)$, whereas for fixed $\theta$ it scales like $\sigma_*^2$ as $\sigma_* \searrow 0$. Both of the terms $\sigma_*$ and $\min\{\zeta_*, s_*/G\}$ in this noise threshold vanish in degenerate situations that are excluded by the conditions of the lemma, which shows the importance of nondegeneracy to reach its conclusion.
}

We are now prepared to prove our main theorem of this subsection, for which we refer to the active-set estimate defined similarly to \eqref{eq.active_QP} as
\bequation\label{eq.active_QP_noisy}
  \widetilde\Acal_{{\rm QP}}(d; x) := \{i \in [q] : \capprox_i(x) + \cgradapprox_i(x)^Td \geq 0\}.
\eequation
As claimed, a similar result would hold for the LP-based subproblem and the active-set estimate
\bequation\label{eq.active_LP_noisy}
  \widetilde\Acal_{{\rm LP}}(d; x) := \{i \in [q] : \capprox_i(x) + \cgradapprox_i(x)^Td \geq 0\}.
\eequation

\btheorem\label{th.PertGra}
  Suppose that, with respect to~\eqref{eq.opt}, the LICQ and strict complementarity hold at a \edit{primal KKT point} $x_* \in \R{n}$.  Let $\Acal_* := \Acal(x_*)$ and $\ell := |\Acal_*|$.  Then, for any $\theta \in \R{}_{>0}$, there exist $\epsilon_x \in \R{}_{>0}$, $\epsilon_e \in \R{}_{>0}$, $\epsilon_c \in \R{}_{>0}$, $\epsilon_{\nabla f} \in \R{}_{>0}$, $\epsilon_{\nabla e} \in \R{}_{>0}$, $\epsilon_{\nabla c} \in \R{}_{>0}$, and $\underline\nu \in \R{}_{>0}$ such that for all $(x,\eapprox(x),\capprox(x),\fgradapprox(x),\egradapprox(x),\cgradapprox(x),\nu)$ with $\|x - x_*\|_2 \leq \epsilon_x$, \eqref{eq.error_bounds}, and $\nu \geq \underline\nu$, the solution $(d, r, s, t)$ of~\eqref{eq.sub_approx_nu} is unique and the active-set estimate defined by \eqref{eq.active_QP_noisy} yields $\widetilde\Acal_{{\rm QP}}(d; x) = \Acal_*$.
\etheorem
\bproof{Proof.}
  Under the stated conditions, it follows by Lemma~\ref{lem.keynoisylemma} that for any $\theta \in \R{}_{>0}$ one has for $x$ sufficiently close to $x_*$ and sufficiently small errors that~\eqref{eq.sub_approx} is feasible and has a unique minimizer at which the LICQ holds.  Thus, by \cite[Theorem~4.4]{HanMang1979}, it follows that for any $\theta \in \R{}_{>0}$ there exist $(\epsilon_x,\epsilon_e,\epsilon_c,\epsilon_{\nabla f},\epsilon_{\nabla e},\epsilon_{\nabla c})$ as stated such that with any $\nu \geq \underline\nu := M$ (where $M$ is defined as in \eqref{eq.lp}) the solution $(d,r,s,t)$ of \eqref{eq.sub_approx_nu} has $(r,s,t) = 0$, which in turn means by Lemma~\ref{lem.keynoisylemma} that $d$ is the unique solution of \eqref{eq.sub_approx}.  All that remains is to observe that by Lemma~\ref{lem.keynoisylemma} one can conclude that for $x$ sufficiently close to $x_*$, sufficiently small noise, and $\nu$ sufficiently large, one also has that \eqref{eq.active_QP_noisy} yields $\Acal_*$.  In particular, this fact follows by the latter equation in \eqref{eq.sub_approx_KKT_2} along with \eqref{eq.sub_approx_KKT_3}.  \edit{Note that the threshold may be taken to be $\underline\nu := \|(y_*,z_*)\|_\infty + 1 \leq M$, since by Lemma~\ref{lem.explicit} the multipliers of \eqref{eq.sub_approx} satisfy $\|(\alpha,\beta) - (y_*,z_*)\|_2 \leq \Theta < 1$ under \eqref{eq.explicit_cond}.}
  \qed
\eproof

\edit{Theorem~\ref{th.PertGra} requires neither optimality of $x_*$ nor any second-order condition, in contrast with Theorem~\ref{the.LP_main_theorem}.  The reason is that \eqref{eq.sub_approx_nu} is strongly convex for every $\theta \in \R{}_{>0}$, so its solution is characterized by the nonsingular linear system \eqref{eq.PerSys} irrespective of the Hessian of the Lagrangian. The second-order condition \eqref{eq:second_order_cond} enters Theorem~\ref{the.LP_main_theorem} only through the error bound of Theorem~\ref{th.OberWrig2.1}, which has no analogue here.  Thus, Theorem~\ref{th.PertGra} identifies the active set of whichever primal KKT point the iterates approach.}

\subsection{\edit{Effect of Degeneracy}}\label{sec.sharp}

{\editcolor
Theorems~\ref{the.LP_main_theorem} and~\ref{th.PertGra} differ in a critical manner. Theorem~\ref{th.PertGra}, which assumes the LICQ and strict complementarity at $x_*$, guarantees correct identification for all $x$ in a neighborhood of $x_*$ under a \emph{fixed} bound on the noise. Theorem~\ref{the.LP_main_theorem}, on the other hand, which assumes only the MFCQ and the second-order condition \eqref{eq:second_order_cond}, and therefore permits weakly active constraints, requires that the noise vanishes linearly as $x \to x_*$; see Corollary~\ref{cor.rate}. It is natural to ask whether the latter requirement might be an artifact of the LP-LPEC strategy or of our analysis. In fact, it is neither. Our next theorem demonstrates two important things for a simple instance with a weakly active constraint.  First, no identification rule can be correct under a fixed noise level uniformly over all instances, since an instance can be paired with a second one whose data are indistinguishable at that noise level, but whose active set differs.  Second, even for that single instance, any reasonable rule can fail at a point $x$ when the error in the objective gradient may be as large as $2\|x - x_*\|_2$, so the linear rate in Corollary~\ref{cor.rate} cannot be improved in general.

To state this precisely, let us call a map $\Rcal$ an \emph{identification rule} if it assigns to each $x \in \R{n}$ and each collection of approximate values and derivatives at a finite number of points that is admissible in the sense of~\eqref{eq.error_bounds} a (possibly random) subset of $[q]$, using no information about the problem other than those values. We call an identification rule \emph{pointwise} if its output at $x$ depends only on the approximate values and derivatives at $x$, and \emph{consistent} if, whenever its input consists of the exact values and derivatives at a local minimizer $\bar x$ of an instance of \eqref{eq.opt} at which the MFCQ and \eqref{eq:second_order_cond} hold, it returns $\Acal(\bar x)$ (almost surely, if it is randomized).  Consistency is a minimal requirement, since a rule that fails it can err even with exact data at the minimizer itself.  All of our aforementioned rules are pointwise, and both the LP-LPEC rule \eqref{eq.active_LP} and the QP rule \eqref{eq.active_QP_noisy} return the correct active set for the instances used in part~(e) below.

\btheorem\label{th.impossible}
  Consider two instances of \eqref{eq.opt} with $n = 2$, $p = 0$, and $q = 1$ defined by
  \bequationNN
    {\rm (P)}:\ \ \left\{\begin{aligned} f(x) &= (x_1-1)^2 + x_2^2 \\ c(x) &= x_1 - 1 \end{aligned} \right\}
    \ \ \text{and}\ \ 
    \text{\rm ($\Phat$)}:\ \ \left\{\begin{aligned} \fhat(x) &= f(x) \\ \chat(x) &= x_1 - 1 - \epsilon_c \end{aligned} \right\}
  \eequationNN
  for some constraint evaluation error bound $\epsilon_c \in \R{}_{>0}$.  Then, the following hold.
  \benumerate[label=$($\alph*$)$]
    \item Both instances have the same unique global minimizer $x_* = (1,0)$ at which the LICQ, the MFCQ, and the second-order condition \eqref{eq:second_order_cond} hold, and both have a unique Lagrange multiplier equal to zero.
    \item $\Acal(x_*) = \{1\}$ for {\rm (P)}, with the constraint being weakly active, and $\hat\Acal(x_*) = \emptyset$ for {\rm ($\Phat$)}.
    \item At every $x \in \R{n}$, the objective and all first-order derivatives of the two instances coincide, and every constraint value in $[c(x) - \epsilon_c, c(x)]$ is admissible under \eqref{eq.error_bounds} for both {\rm (P)} and {\rm ($\Phat$)}. For example, $\capprox(x) := x_1 - 1 - \thalf\epsilon_c$ is admissible for both problems for all $x \in \R{n}$.
    \item No identification rule is correct for both instances. In particular, for every rule $\Rcal$ there is an instance and an admissible realization of the data for which $\Rcal$ does not return the correct active set at $x_*$.  If, moreover, a rule $\Rcal$ uses a single evaluation of the constraint value, whose noise is uniform on $[-\epsilon_c,\epsilon_c]$, and the noise in all other quantities is zero, then for every rule $\Rcal$, one has
    \bequationNN
      \P_{\rm (P)} \left[\Rcal \neq \Acal(x_*)\right] + \P_{\rm (\hat{P})} \left[\Rcal \neq \hat{\Acal}(x_*)\right] \geq \thalf,
    \eequationNN
    so that $\Rcal$ errs with probability at least $\tfrac14$ on one of the two instances at every $x$.
    \item Let $\Rcal$ be a pointwise and consistent identification rule.  Then, for every $x \in \R{2}$ with $x \neq x_*$ and $c(x) \leq 0$, and for any error bounds satisfying $\epsilon_{\nabla f} \geq 2\|x - x_*\|_2$, there is a realization of the approximations satisfying \eqref{eq.error_bounds} for {\rm (P)} under which $\Rcal$ applied at $x$ does not return $\Acal(x_*)$.  Consequently, any guarantee of correct identification by such a rule at the feasible points near $x_*$ requires that $\epsilon_{\nabla f} < 2\|x - x_*\|_2$, i.e., that the noise vanish at least linearly as $x \to x_*$.
  \eenumerate
\etheorem
\bproof{Proof.}
  Consider part (a). For both instances, $x_* = (1,0)$ is the unconstrained global minimizer of $f$ and is feasible, so it is the unique global minimizer of both instances. Furthermore, $\nabla c(x_*) = \nabla \chat(x_*) = (1,0)$, so the LICQ and MFCQ hold in both instances, and for both the unique Lagrange multiplier is $z_* = 0$ since $\nabla f(x_*) = 0$. Finally, $\nabla^2_{xx}\Lcal = 2I$ is positive definite for both instances, so \eqref{eq:second_order_cond} holds for both.
  
  Now consider part (b). The active sets are immediate from $c(x_*) = 0$ and $\chat(x_*) = -\epsilon_c$. Since $c(x_*) = 0$ and $z_* = 0$, the constraint of {\rm (P)} is weakly active, i.e., $\Acal_0(x_*) = \Acal(x_*) = \{1\}$.
  
  Next, consider part (c). The two instances differ only by the constant $\epsilon_c$ in the constraint value, so all first-order derivatives agree and, using the example, $|\capprox(x) - c(x)| = |\capprox(x) - \chat(x)| = \thalf\epsilon_c \leq \epsilon_c$. 

  Next, consider part (d). The first claim follows from (b) and (c). For example, a rule sees the same input under both instances when the example in (c) occurs, and hence returns the same output, which cannot equal both $\{1\}$ and $\emptyset$.  For the second claim, let $\mu$ and $\hat{\mu}$ denote the laws of the observed constraint value at an arbitrary point~$x$ under instances {\rm (P)} and {\rm ($\Phat$)}, respectively; namely, these are the uniform distributions on $[c(x) - \epsilon_c, c(x) + \epsilon_c]$ and on $[c(x) - 2\epsilon_c, c(x)]$, which have variational distance between them equal to $\thalf$.  For any (possibly randomized) rule $\Rcal$, writing $\Ecal := \{\Rcal = \{1\}\}$ and $\Ecal^c := \{\Rcal = \emptyset\}$, one finds
  \bequationNN
    \P_{{\rm (P)}}[\Rcal \neq \{1\}] + \P_{\text{($\Phat$)}}[\Rcal \neq \emptyset] = \P_{{\rm (P)}}[\Ecal^c] + \P_{\text{($\Phat$)}}[\Ecal] = 1 - (\P_{{\rm (P)}}[\Ecal] - \P_{\text{($\Phat$)}}[\Ecal]) \geq 1 - \text{TV}(\mu,\hat{\mu}) = \thalf,
  \eequationNN
  where TV denotes variational distance (or total variation).

  Finally, consider part (e).  Let $x \in \R{2}$ satisfy $x \neq x_*$ and $c(x) \leq 0$, let $\eta \in (0,\epsilon_c]$, and consider the instance {\rm (P$_x$)} of \eqref{eq.opt} defined by $f_x(w) := \|w - x\|_2^2 + f(x)$ and $c_x(w) := c(w) - \eta$.  Since $c_x(x) = c(x) - \eta < 0$, the point $x$ is feasible for {\rm (P$_x$)} and is the unconstrained global minimizer of $f_x$, so $x$ is the unique global minimizer of {\rm (P$_x$)}, with unique Lagrange multiplier equal to zero and active set $\Acal_x(x) = \emptyset$.  Moreover, the LICQ, the MFCQ, and \eqref{eq:second_order_cond} hold at $x$ for {\rm (P$_x$)}, since $\nabla c_x(x) = (1,0) \neq 0$ and $\nabla^2 f_x = 2I$.  Now consider approximations for {\rm (P)} that are exact at every point other than $x$ and that at $x$ are given by
  \bequationNN
    \baligned
      \fapprox(x) &:= f(x) = f_x(x), & \capprox(x) &:= c(x) - \eta = c_x(x), \\
      \fgradapprox(x) &:= 0 = \nabla f_x(x), & \text{and}\ \ \cgradapprox(x) &:= \nabla c(x) = \nabla c_x(x).
    \ealigned
  \eequationNN
  At $x$, the errors are $|\capprox(x) - c(x)| = \eta \leq \epsilon_c$ and $\|\fgradapprox(x) - \nabla f(x)\|_2 = \|2(x - x_*)\|_2 = 2\|x - x_*\|_2 \leq \epsilon_{\nabla f}$, with all other errors equal to zero, so these approximations satisfy \eqref{eq.error_bounds}.  Since $\Rcal$ is pointwise, its output at $x$ depends only on these values, which coincide with the exact values and derivatives of {\rm (P$_x$)} at its minimizer $x$.  Since $\Rcal$ is consistent, it returns $\Acal_x(x) = \emptyset$, whereas $\Acal(x_*) = \{1\}$ for {\rm (P)}.  For our proposed strategies, this can also be checked directly.  For the LP-LPEC rule, the objective of \eqref{prob.lp-lpec} with these values is $|z| + (\eta - c(x))z$, which is minimized over $z \in [0,M]$ at $z = 0$ with value zero, so the solution is $z_x = 0$, $\tilde{\bar\rho}(x,y_x,z_x) = 0$, and \eqref{eq.active_LP} excludes the constraint since $\capprox(x) < 0$.  For the QP rule, since $\fgradapprox(x) = 0$ and $\capprox(x) < 0$, the unique solution of \eqref{eq.sub_approx_nu} has $d = 0$ and $t = 0$, so \eqref{eq.active_QP_noisy} again excludes the constraint.
  \qed
\eproof

The probabilistic statement in Theorem~\ref{th.impossible}(d) concerns a single noisy evaluation.  If, instead, $\Rcal$ uses~$m$ independent evaluations of the constraint value, then the variational distance between the corresponding product laws is $1 - 2^{-m}$ and the same argument gives the lower bound $2^{-m}$ in place of $\thalf$.

The construction in parts~(a)--(d) is not specific to the example in Theorem~\ref{th.impossible}. More generally, suppose that $x_*$ is a local minimizer at which the MFCQ and \eqref{eq:second_order_cond} hold and that $\Acal_0(x_*) = \{i_0\}$ for some index $i_0 \in [q]$.  Replacing the constraint function $c_{i_0}$ by $c_{i_0} - \delta$ for any $\delta \in (0,\epsilon_c]$ yields a second problem instance for which $x_*$ is again a KKT point with the same set of optimal Lagrange multipliers and with active set $\Acal(x_*) \setminus \{i_0\}$.  Moreover, its critical cone at $x_*$ is $\{d : \nabla e(x_*)^Td = 0,\ \nabla c_i(x_*)^Td = 0\ \text{for all}\ i \in \Acal(x_*) \setminus \{i_0\}\} = \Ccal(x_*) \cup (-\Ccal(x_*))$, with respect to which~\eqref{eq:second_order_cond} still holds.  Hence, $x_*$ is a strict local minimizer of the second instance whose data are indistinguishable from those of the first instance at the noise level $\epsilon_c$, even though its active set at $x_*$ is strictly smaller than that of the first instance.

Parts~(a)--(d) of Theorem~\ref{th.impossible} thus show that correct identification under a fixed noise level, uniformly over instances, requires a condition ruling out weakly active constraints, such as strict complementarity.  Part~(e) shows more. In particular, it shows that even for a single instance with a weakly active constraint, the noise must vanish at least linearly in $\|x - x_*\|_2$ as $x \to x_*$.  For the instance {\rm (P)}, Corollary~\ref{cor.rate} ensures correct identification by the LP-LPEC rule whenever $\|\epsilonerror\|_\infty \leq \Lambda'\|x - x_*\|_2$, where $\Lambda' \leq \thalf$ by the proof of that corollary (since $C \geq 1$ and $c_\psi \geq 1$), whereas part~(e) of Theorem~\ref{th.impossible} shows failure whenever $\epsilon_{\nabla f} \geq 2\|x - x_*\|_2$.  Hence, for {\rm (P)}, the noise vanishing at a linear rate is both sufficient and necessary for correct identification, so the rate in \eqref{eq.rate_lin} cannot be improved.  We note that part~(e) concerns a worst-case realization of the noise in which the error in the objective gradient cancels $\nabla f(x)$ exactly.  Under random noise, such a cancellation is unlikely, and correct identification may then be observed at larger noise levels.  Viewed from the perspective of perturbation theory, one can say that in noisy settings the active set is not a stable feature of the problem data at a degenerate local minimizer.

The two strategies analyzed in this paper demonstrate this fact.  The LP-LPEC strategy of \S\ref{sec.lagrange} has been analyzed for the degenerate case, and it pays for this generality with the requirement of Corollary~\ref{cor.rate} that the noise vanish linearly, which by Theorem~\ref{th.impossible}(e) cannot be avoided, or improved in rate, in general.  The primal-step strategies of \S\ref{sec.primal_step}, on the other hand, assume the LICQ and strict complementarity and, consequently, can tolerate a fixed noise level.  Overall, the difference between Theorems~\ref{the.LP_main_theorem} and~\ref{th.PertGra} should be read as a difference between regularity assumptions rather than between active-set identification schemes.
}

\subsection{Active-Set Identification within a Stochastic Algorithm}\label{sec.stochastic_algorithm}

Let us now consider the use of active-set identification procedures within the context of a stochastic algorithm for solving~\eqref{eq.opt}.  Let $(\Omega,\Fcal,\P)$ be a probability space that captures all outcomes of a run of the algorithm.  That is, each run of the algorithm is associated with an outcome $\omega \in \Omega$, an infinite-dimensional tuple whose $k$th element determines the outcome of any random quantity in iteration $k \in \N{}$.  In this manner, the algorithm defines a stochastic process $\{X_k\}$, where for all $k \in \N{}$ the random variable $X_k$ is the primal iterate.  Given the initial condition that $X_1(\omega) = x_1$ for some $x_1 \in \R{n}$ for all $\omega \in \Omega$, let $\Fcal_k$ denote the sub-$\sigma$-algebra of $\Fcal$ corresponding to the initial condition and all of the stochastic quantities prior to iteration $k$.  In this manner $\{\Fcal_k\}$ is a filtration.

\edit{The conditions of our main results in this section rely on assumptions that say: ``conditioned on the behavior of an algorithm up to iteration $k \in \N{}$, the error bounds in~\eqref{eq.error_bounds} hold at the current iterate with probability at least $p$.''  Thus, before stating those results, let us present the following lemma that shows that such conditions are achievable with sufficiently large sample sizes.  Recall that \eqref{eq.error_bounds} is a condition on a realization, meaning that nothing in \S\ref{sec.lagrange} or \S\ref{sec.primal_step} requires the errors to be bounded almost surely.}

{\editcolor
\blemma\label{lem.sample}
  Consider arbitrary $k \in \N{}$ and suppose that a sample size $b_k \in \N{}$ and the error bounds in \eqref{eq.error_bounds} are $\Fcal_k$-measurable.  Suppose also that each of the six approximations in \eqref{eq.error_bounds} with $x = X_k$ is computed as an average of $b_k$ estimators that, conditioned on $\Fcal_k$, are independent and unbiased, and whose conditional mean-squared errors---measured in the $\ell_2$-norm for vectors and in the Frobenius norm for matrices---are bounded by $\gamma_f^2$, $\gamma_e^2$, $\gamma_c^2$, $\gamma_{\nabla f}^2$, $\gamma_{\nabla e}^2$, and $\gamma_{\nabla c}^2$, respectively.  Then,
  \bequationNN
    \baligned
      \P[\text{\eqref{eq.error_bounds} holds with $x = X_k$}\ |\ \Fcal_k] &\geq 1 - \tfrac{S(\epsilonerror)}{b_k}, \\
      \text{where}\ \ S(\epsilonerror) &:= \tfrac{\gamma_f^2}{\epsilon_f^2} + \tfrac{\gamma_e^2}{\epsilon_e^2} + \tfrac{\gamma_c^2}{\epsilon_c^2} + \tfrac{\gamma_{\nabla f}^2}{\epsilon_{\nabla f}^2} + \tfrac{\gamma_{\nabla e}^2}{\epsilon_{\nabla e}^2} + \tfrac{\gamma_{\nabla c}^2}{\epsilon_{\nabla c}^2}.
    \ealigned
  \eequationNN
  Thus, \eqref{eq.error_bounds} holds with $x = X_k$ with conditional probability at least $p \in (0,1)$ whenever $b_k \geq S(\epsilonerror)/(1-p)$.
\elemma
\bproof{Proof.}
  Consider arbitrary $k \in \N{}$ and, without loss of generality, consider the approximation $\cgradapprox(X_k)$. (The other five approximations can be handled identically.) Let $\cgradapprox(X_k)$ be the average of $b_k$ estimators, say, $\{G_1,\dots,G_{b_k}\}$ satisfying the stated conditions. Writing $\langle \cdot,\cdot \rangle_F$ for the Frobenius inner product, conditional independence and unbiasedness give $\E[\langle G_i - \nabla c(X_k), G_j - \nabla c(X_k) \rangle_F \mid \Fcal_k] = 0$ for $i \neq j$ so that
  \bequationNN
    \E[\| \cgradapprox(X_k) - \nabla c(X_k)\|_F^2\ |\ \Fcal_k ] = \tfrac{1}{b_k^2} \sum_{j=1}^{b_k} \E[\|G_j - \nabla c(X_k)\|_F^2\ |\ \Fcal_k] \leq \tfrac{\gamma_{\nabla c}^2}{b_k}.
  \eequationNN
  Since for any matrix its spectral norm is bounded above by its Frobenius norm, it follows from the conditional Markov inequality applied to $\|\cgradapprox(X_k) - \nabla c(X_k)\|_F^2$ that
  \bequationNN
    \P[\|\cgradapprox(X_k) - \nabla c(X_k)\|_2 > \epsilon_{\nabla c} \ |\ \Fcal_k]
    \leq \P[\|\cgradapprox(X_k) - \nabla c(X_k)\|_F^2 \geq \epsilon_{\nabla c}^2 \ |\ \Fcal_k]
    \leq \tfrac{\gamma_{\nabla c}^2}{b_k\epsilon_{\nabla c}^2},
  \eequationNN
  where the last inequality uses that $\epsilon_{\nabla c}$ is $\Fcal_k$-measurable.  Summing the six analogous bounds (a union bound) gives the claim.  Finally, $1 - S(\epsilonerror)/b_k \geq p$ if and only if $b_k \geq S(\epsilonerror)/(1-p)$.
  \qed
\eproof

The form of the conditions in Lemma~\ref{lem.sample}, namely, that a deterministic accuracy condition holds with probability at least $p$ conditioned on the past, is the \emph{probabilistically accurate model} form that has appeared throughout the literature on stochastic optimization with random models; see, e.g., \cite{CartSche2018,BlanCartMeniSche2019,PaquSche2020,BeraCaoSche2021}.

Lemma~\ref{lem.sample} suggests how the two settings of \S\ref{sec.lagrange} and \S\ref{sec.primal_step} can be achieved.  Under the LICQ and strict complementarity, one should expect a fixed sample size for each iteration to be sufficient; see Theorem~\ref{cor.qp} below. On the other hand, under only the MFCQ and \eqref{eq:second_order_cond}, it has been shown in Corollary~\ref{cor.rate} that it suffices for the error bounds at iteration $k$ to satisfy $\|\epsilon_{{\rm err},k}\|_\infty \leq \Lambda'\|X_k - x_*\|_2$.  In that case $S(\epsilon_{{\rm err},k})$ should grow like $\|X_k - x_*\|_2^{-2}$, and by Lemma~\ref{lem.sample} a sample size $b_k \geq S(\epsilon_{{\rm err},k})/(1-p)$ growing at that rate suffices.  This is the adaptive-sampling regime studied in \cite{ByrdChinNoceWu2012,BollByrdNoce2018}. Variance-reduced estimators, such as those of SVRG~\cite{JohnZhan2013} and SARAH~\cite{NguyLiuScheTaka2017}, achieve a comparable effect at lower cost, and in the deterministic-noise setting the analogue is a simulation or rule whose precision can be increased successively~\cite{curtis2025interior}.  For instance, if $f$ is an average of functions whose gradients are Lipschitz continuous with constant $L_f$, then the SVRG estimator of $\nabla f(X_k)$, formed with an $\Fcal_k$-measurable reference point~$\tilde X_k$, is conditionally unbiased with conditional mean-squared error at most $L_f^2\|X_k - \tilde X_k\|_2^2$, so Lemma~\ref{lem.sample} applies with $\gamma_{\nabla f}^2 := L_f^2\|X_k - \tilde X_k\|_2^2$.  Consequently, if the reference point is updated often enough that $\|\tilde X_k - x_*\|_2 \leq r\|X_k - x_*\|_2$ for some $r \in \R{}_{>0}$, then the bound on $\epsilon_{\nabla f}$ required by Theorem~\ref{ass.lp} can be met with a \emph{fixed} sample size.  (The SARAH estimator, conditioned on $\Fcal_k$, is biased, so Lemma~\ref{lem.sample} does not apply.)  A method that uses plain sample averaging with a fixed sample size throughout does not satisfy the hypothesis of Theorem~\ref{ass.lp}, and, by Theorem~\ref{th.impossible}, no method using data with a fixed noise level can succeed uniformly when $x_*$ has a weakly active constraint.
}

\edit{Let us now state specific theorems that correspond to the settings of \S\ref{sec.lagrange} and \S\ref{sec.primal_step}, respectively. Our first setting, for identification through Lagrange multipliers presented in \S\ref{sec.lagrange}, requires that the MFCQ and the second-order condition \eqref{eq:second_order_cond} hold at a local minimizer $x_*$.  For all $k \in \N{}$, the LP \eqref{prob.lp-lpec} is solved with $x = X_k$ to yield a Lagrange multiplier pair $(Y_k,Z_k)$.  The conditions of Theorem~\ref{the.LP_main_theorem} involve $(Y_k,Z_k)$, and hence the noise realized at iteration $k$, so they cannot serve directly as a hypothesis on the error bounds. Fortunately, Corollary~\ref{cor.rate} removes this dependence by giving a sufficient bound in terms of $\|X_k - x_*\|_2$ alone.  The following theorem therefore assumes that, at each iterate, the error bounds are at most proportional to $\|X_k - x_*\|_2$ and that \eqref{eq.error_bounds} holds with these bounds with probability at least $p$, and it concludes that, with probability one, beyond some iteration the optimal active set is identified correctly with conditional probability at least $p$.  This requires the noise in the function and derivative estimates to be reduced as the iterates converge at the rate stated previously.  By Theorem~\ref{th.impossible}(e), a requirement of this kind, at this rate, cannot be avoided in general when $x_*$ has a weakly active constraint.}

{\editcolor
\btheorem\label{ass.lp}
  Suppose that, with respect to \eqref{eq.opt}, the conditions of Theorem~\ref{the.LP_main_theorem} hold, consider any $\beta \in \R{}_{>0}$ and $\sigma \in (0,1)$, and let $(\hat\epsilon,\Lambda')$ be defined as in Theorem~\ref{the.LP_main_theorem} and Corollary~\ref{cor.rate}.  Suppose also that, almost surely, $X_k \neq x_*$ for all $k \in \N{}$.  Suppose further that, for all $k \in \N{}$, the function and derivative approximations at $X_k$ are computed subject to positive error bounds $\epsilon_{{\rm err},k}$ that are $\Fcal_k$-measurable with
  \bequation\label{eq.thm43_bounds}
    \|\epsilon_{{\rm err},k}\|_\infty \leq \Lambda'\|X_k - x_*\|_2\ \ \text{whenever}\ \ \|X_k - x_*\|_2 \leq \hat\epsilon,
  \eequation
  and that, for some universal constant $p \in (0,1]$, one has, almost surely and for all $k \in \N{}$,
  \bequation\label{eq.thm43_hyp}
    \P\left[\text{\eqref{eq.error_bounds} holds with $x = X_k$ and $\epsilonerror = \epsilon_{{\rm err},k}$}\ \middle|\ \Fcal_k\right] \geq p \cdot \ones_{\{\|X_k - x_*\|_2 \leq \hat\epsilon\}},
  \eequation
  which is well defined since $X_k$ is $\Fcal_k$-measurable.  Let $(Y_k,Z_k)$ denote a solution of \eqref{prob.lp-lpec} with $x = X_k$, computed from these approximations.  Then, if $\{X_k\}$ converges almost surely to $x_*$, it holds with probability one that there exists $k_0 \in \N{}$ such that
  \bequationNN
    \P\left[\widetilde\Acal_{{\rm LP}}(Y_k,Z_k;X_k) = \Acal_*\ \middle|\ \Fcal_k\right] \geq p\ \ \text{for all}\ k \in \N{}\ \text{with}\ k \geq k_0.
  \eequationNN
  Moreover, with probability one, $\widetilde\Acal_{{\rm LP}}(Y_k,Z_k;X_k) = \Acal_*$ for infinitely many $k \in \N{}$.
\etheorem
\bproof{Proof.}
  Let $\Omega_0$ be the event, with probability one, on which $X_k \to x_*$ and $X_k \neq x_*$ for all $k \in \N{}$, and for $\tilde\omega \in \Omega_0$ let $k_0(\tilde\omega) := \min\{k \in \N{} : \|X_j(\tilde\omega) - x_*\|_2 \leq \hat\epsilon\ \text{for all}\ j \geq k\}$, which is finite.  For $k \in \N{}$, let $\Ecal_k$ denote the event in \eqref{eq.thm43_hyp}, let $\Bcal_k := \{0 < \|X_k - x_*\|_2 \leq \hat\epsilon\}$, which lies in $\Fcal_k$, and let $\Acal^{\rm id}_k$ denote the event that $\widetilde\Acal_{{\rm LP}}(Y_k,Z_k;X_k) = \Acal_*$.  On $\Ecal_k \cap \Bcal_k$, the approximations at $X_k$ satisfy \eqref{eq.error_bounds} with error bounds satisfying $\|\epsilon_{{\rm err},k}\|_\infty \leq \Lambda'\|X_k - x_*\|_2$ by \eqref{eq.thm43_bounds}, so Corollary~\ref{cor.rate}, via \eqref{eq.rate_lin}, applied with $x = X_k$ yields $\Ecal_k \cap \Bcal_k \subseteq \Acal^{\rm id}_k$.  Since $\Bcal_k \in \Fcal_k$, it follows from \eqref{eq.thm43_hyp} that, almost surely, $\P[\Acal^{\rm id}_k \mid \Fcal_k] \geq \P[\Ecal_k \cap \Bcal_k \mid \Fcal_k] = \ones_{\Bcal_k}\P[\Ecal_k \mid \Fcal_k] \geq p \cdot \ones_{\Bcal_k}$. On $\Omega_0$ one has $\ones_{\Bcal_k} = 1$ for all $k \geq k_0$, which gives the first conclusion, and also $\sum_{k \in \N{}} \P[\Acal^{\rm id}_k \mid \Fcal_k] \geq \sum_{k \geq k_0} p = \infty$.  The samples used to form the approximations at iteration $k$ are among the stochastic quantities determined prior to iteration $k+1$, so $\Acal^{\rm id}_k \in \Fcal_{k+1}$, and L\'evy's extension of the Borel--Cantelli lemma \cite[\S12.15]{Will1991} states that the events $\{\Acal^{\rm id}_k\ \text{infinitely often}\}$ and $\{\sum_{k} \P[\Acal^{\rm id}_k \mid \Fcal_k] = \infty\}$ are almost surely equivalent.  Hence $\Acal^{\rm id}_k$ occurs infinitely often almost surely on $\Omega_0$, and $\P[\Omega_0] = 1$.
  \qed
\eproof
}

\edit{
The quantity $\|X_k - x_*\|_2$ in \eqref{eq.thm43_bounds} is unknown in practice, so the error bounds there cannot be computed directly.  They are, however, equivalent up to constants to bounds in terms of the KKT error.  Indeed, for any $x$ with $\|x - x_*\|_2 \leq \hat\epsilon$, one has $\omega(x) \geq C^{-1}\|x - x_*\|_2$ as shown in the proof of Theorem~\ref{the.LP_main_theorem}, while for any $(y_*,z_*) \in \Dcal(x_*)$, Theorem~\ref{th.OberWrig2.1} gives $\omega(x) \leq \psi(x,y_*,z_*) \leq C\dist(\{(x,y_*,z_*)\},\{x_*\} \times \Dcal(x_*)) \leq C\|x - x_*\|_2$.  Hence \eqref{eq.thm43_bounds} holds if $\|\epsilon_{{\rm err},k}\|_\infty \leq (\Lambda'/C)\,\omega(X_k)$, i.e., if the noise is bounded by a constant multiple of the KKT error at $X_k$, and conversely a bound of the form \eqref{eq.thm43_bounds} implies one of this form. Observations of this type motivate various adaptive sampling strategies in practice.
}

\edit{Our second setting, for identification through a primal step as presented in \S\ref{sec.primal_step}, requires the LICQ and strict complementarity at a KKT point~$x_*$ rather than the MFCQ and the second-order condition \eqref{eq:second_order_cond}.  For this setting in the following theorem, we recall that for all $k \in \N{}$ the active-set estimate is obtained by solving the QP~\eqref{eq.sub_approx_nu}, which yields the subproblem solution we denote as $(D_k, R_k, S_k, T_k)$.}

\btheorem\label{cor.qp}
  Suppose that, with respect to \eqref{eq.opt}, the conditions of Theorem~\ref{th.PertGra} hold.  In addition, with $(\epsilon_x, \epsilon_e, \epsilon_c, \epsilon_{\nabla f}, \edit{\epsilon_{\nabla e}}, \epsilon_{\nabla c}, \underline\nu)$ as defined in Theorem~\ref{th.PertGra}, suppose \edit{that for some universal constant $p \in (0,1]$ one has, almost surely and for all $k \in \N{}$,}
  \bequation\label{eq.thm44_hyp}
    \P[\text{\eqref{eq.error_bounds} holds with $x = X_k$}\ |\ \Fcal_k] \geq p \cdot \ones_{\{\|X_k - x_*\|_2 \leq \epsilon_x\}},
  \eequation
  \edit{which is well defined because $X_k$ is $\Fcal_k$-measurable, and which holds with $p$ given explicitly by Lemma~\ref{lem.sample} under sampling with a \emph{fixed} sample size.}  Then, for any parameter $\theta \in \R{}_{>0}$ and $\nu \geq \underline\nu$, if $\{X_k\}$ converges almost surely to~$x_*$, it holds with probability one that there exists an index \edit{$k_0 \in \N{}$} such that for all $k \in \N{}$ with \edit{$k \geq k_0$} it holds with \edit{conditional} probability at least $p$ that the solution $(D_k, R_k, S_k, T_k)$ of~\eqref{eq.sub_approx_nu} with $x = X_k$ is unique and the active-set estimate defined by \eqref{eq.active_QP_noisy} yields $\widetilde\Acal_{{\rm QP}}(D_k; X_k) = \Acal_*$.  \edit{Moreover, with probability one, the active-set estimate has $\widetilde\Acal_{{\rm QP}}(D_k;X_k) = \Acal_*$ for infinitely many $k \in \N{}$.}
\etheorem
\bproof{Proof.}
  Same proof as Theorem~\ref{ass.lp}, with $\Omega_0$ the event on which $X_k \to x_*$, with $\Ecal_k$ the event in \eqref{eq.thm44_hyp}, with $\Bcal_k := \{\|X_k - x_*\|_2 \leq \epsilon_x\}$, and with Theorem~\ref{th.PertGra} in place of Corollary~\ref{cor.rate}.  Neither $X_k \neq x_*$ nor error bounds that vanish as $X_k \to x_*$ are needed, since the error bounds in Theorem~\ref{th.PertGra} are fixed.
  \qed
\eproof

By our previous arguments, a similar result would hold if the QP subproblem~\eqref{eq.sub_approx_nu} were replaced by an LP subproblem with a trust-region constraint on the primal step.

\section{Numerical Results}\label{sec.numerical}

The purpose of our experiments is to demonstrate the behavior of the active-set identification techniques analyzed in \S\ref{sec.noisy} in the context of a few test problems.  We present the results of two sets of experiments.  Our first set of experiments focuses on two related two-dimensional test problems that allow us to provide clear visualizations of active-set identification behavior in both deterministic and noisy settings.  Our second set of experiments focuses on behavior in the context of a larger-dimensional neural-network training problem.  We emphasize that any number of problems could be used to demonstrate the behavior of the active-set identification methods from \S\ref{sec.noisy}, each potentially leading to interesting conclusions about the techniques in general, the parameter choices for them, and so on.  We chose our first experiments since two-dimensional problems allow for nice visualizations, and we chose our second experiments since they involve the problem of training of a neural network, which is of immense interest in modern optimization.

\subsection{Visualization with 2D Problems}\label{sec.two-dimensional}

Consider a pair of two-dimensional problems with
\bequation\label{prob.two-dimensional}
  \min_{x \in \R{2}}\ f(x)\ \st\ c_1(x) := x_1^2 - x_2 \leq 0\ \text{and}\ c_2(x) := x_1^2 + x_2 - \tfrac12 \leq 0,
\eequation
where $c_i$ and $x_i$ for each $i \in \{1,2\}$ denote the $i$th component of the constraint function $c : \R{2} \to \R{2}$ and the $i$th component of the vector $x \in \R{2}$, respectively.  (Since there are no equality constraints in \eqref{prob.two-dimensional}, we denote the constraint function simply as $c$.)  The feasible region is thus the region in $\R{2}$ contained between two parabolas, one upward-facing and one downward-facing.  Let us also consider two objective functions $f_1 : \R{2} \to \R{}$ and $f_2 : \R{2} \to \R{}$ that are defined respectively as
\bequationNN
  f_1(x) = (x_1 + \thalf)^2 + 4(x_2 - \thalf)^2\ \ \text{and}\ \ f_2(x) = 4(x_1 + \tfrac35)^2 + (x_2 - \tfrac14)^2.
\eequationNN
The problems to minimize $f_1$ and minimize $f_2$, each subject to $c(x) \leq 0$, have convex feasible regions and strongly convex objectives.  Therefore, each problem has a unique solution, call it $x_*$ in each case.  It is straightforward to confirm that the problem to minimize $f_1$ subject to $c(x) \leq 0$ has one active constraint at $x_*$ such that $\Acal(x_*) = \{2\}$, and at $x_*$ the LICQ and strict complementarity hold.  On the other hand, for the problem to minimize $f_2$ subject to $c(x) \leq 0$, both constraints are active at $x_*$, i.e., $\Acal(x_*) = \{1,2\}$, and the LICQ and strict complementarity hold at $x_*$.

For each problem, we generated points in a grid around $x_*$.  For each point $x$ in the grid, we solved various instances of the LP-LPEC subproblem \eqref{prob.lp-lpec} with $M \gets 10^8$ and \eqref{eq.active_LP} applied with $\beta \gets 0.7071$ \edit{(determined as $\beta \gets 1/\max\{\|\nabla f(\xbar)\|_\infty, \|c(\xbar)\|_\infty, \|\nabla c(\xbar)\|_\infty\}$ at a point $\xbar$ chosen randomly in the grid, so that $\beta\bar\rho$ is scaled to the problem data)} and $\sigma \gets 0.7$; an LP with a trust-region constraint (that is \eqref{eq.sub_LP} with approximate function and derivative values) with $\nu \gets 100$ and trust-region radius $\Delta \gets 0.1$; and the QP subproblem \eqref{eq.sub_approx_nu} with $\nu \gets 100$ and $\theta \gets 5$.  The active-set estimates were determined accordingly, as described in \S\ref{sec.lagrange} and \S\ref{sec.primal_step}.  For each subproblem, we first considered the deterministic setting, denoted as $\epsilon = 0$.  For each point in the grid, we determined whether the active-set estimate exactly matched $\Acal(x_*)$ or not, yielding a value of 1 or 0, respectively.  Then, for two noise levels, denoted by $\epsilon \in \{10^{-2}, 10^{-1}\}$, respectively, we solved~8 instances of each subproblem with noise added to the function and derivative values, then determined the fraction of times that the active-set estimate exactly matched~$\Acal(x_*)$.  This yielded a fractional value between~0 and 1 for each point, where a value of 1 indicates that the active set was always identified correctly.  With respect to the noise, for each function value or component of a function or derivative value, we added a number generated from a uniform distribution over the interval $[-\epsilon,\epsilon]$, determined by the noise level.

\begin{figure}[ht]
    \centering
    \includegraphics[width=0.2\linewidth, trim=310 0 340 38, clip]{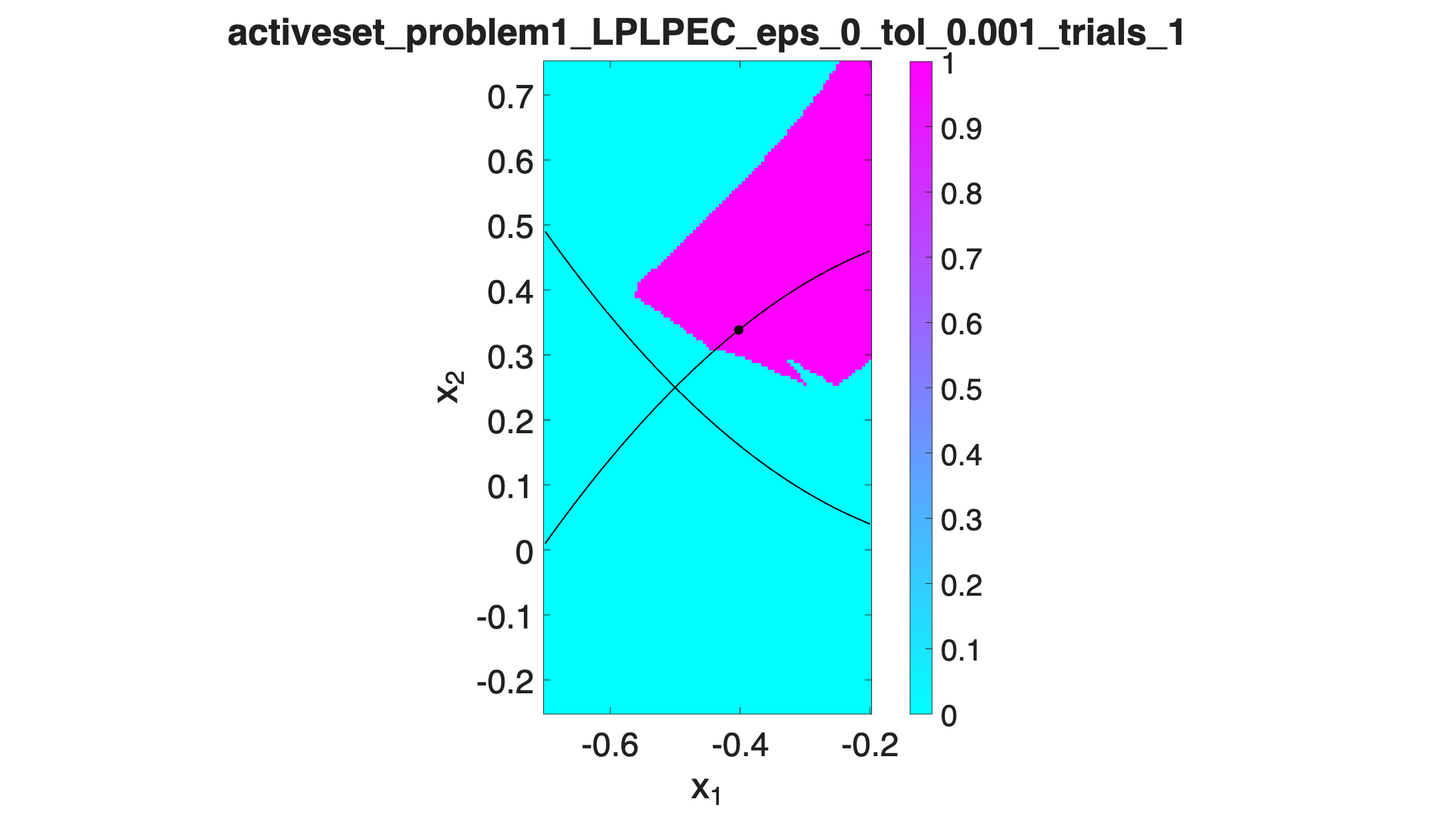}
    \hfill
    \includegraphics[width=0.2\linewidth, trim=310 0 340 38, clip]{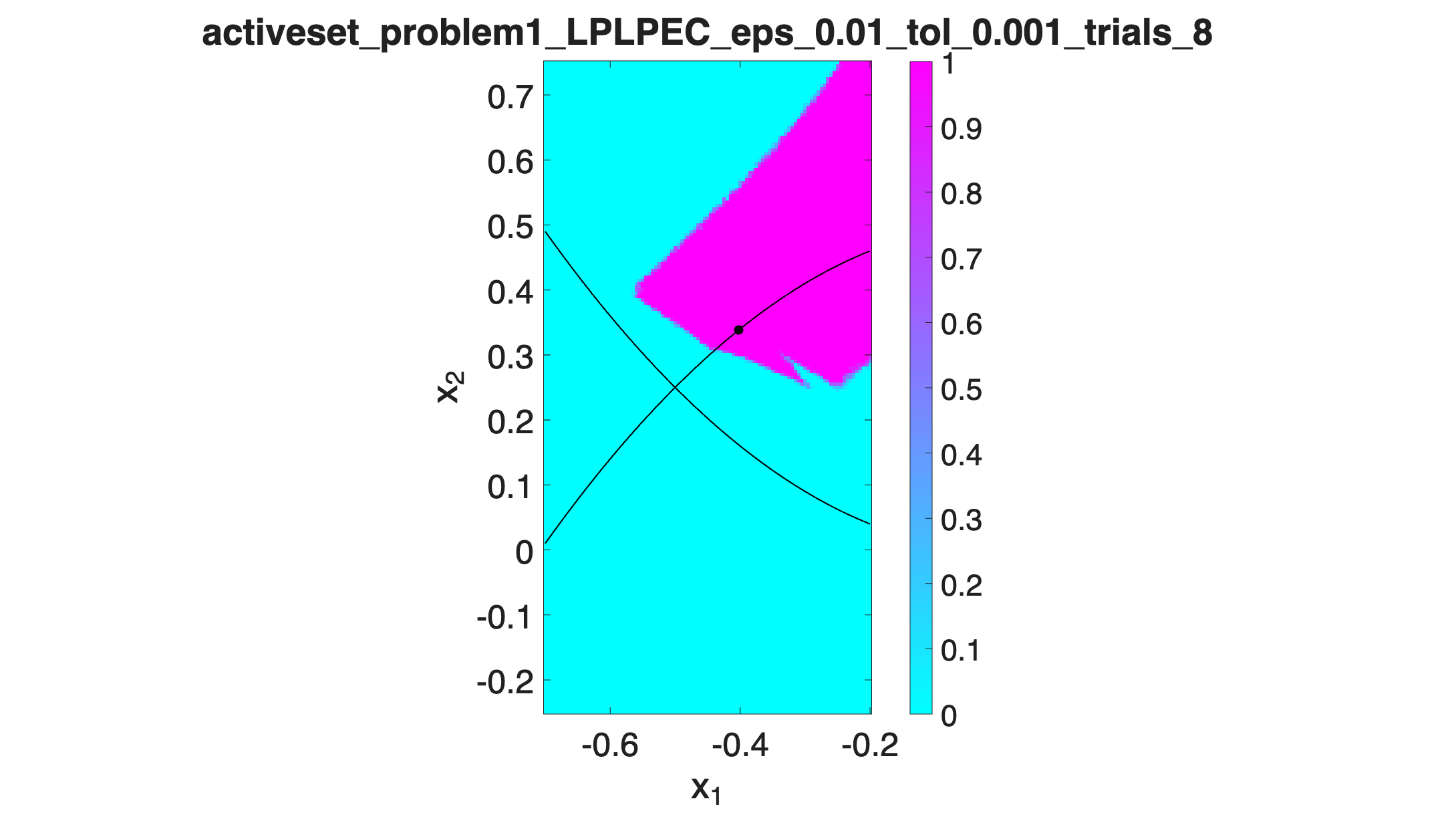}
    \hfill
    \includegraphics[width=0.2\linewidth, trim=310 0 340 38, clip]{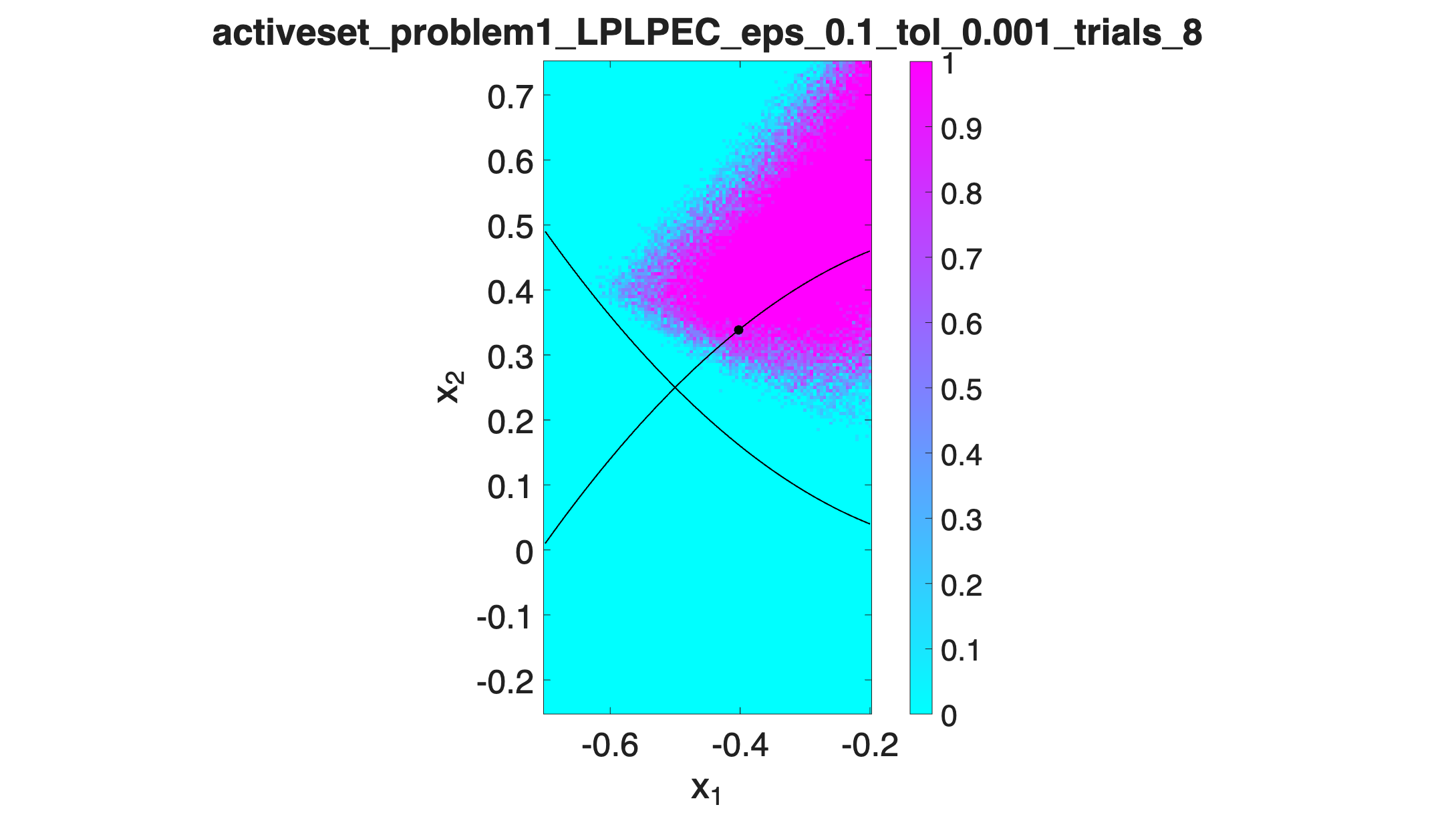}

    \includegraphics[width=0.2\linewidth, trim=310 0 340 38, clip]{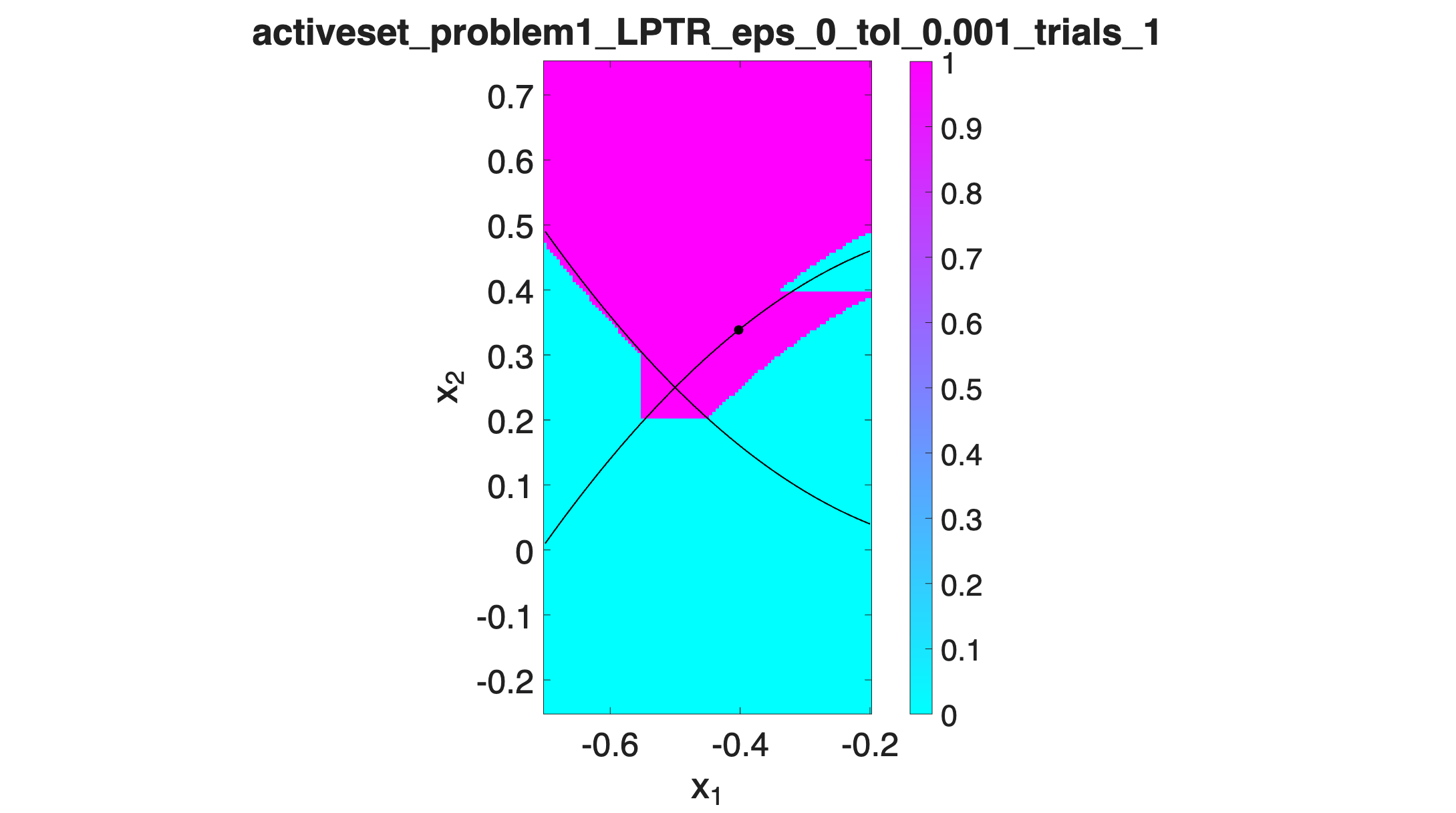}
    \hfill
    \includegraphics[width=0.2\linewidth, trim=310 0 340 38, clip]{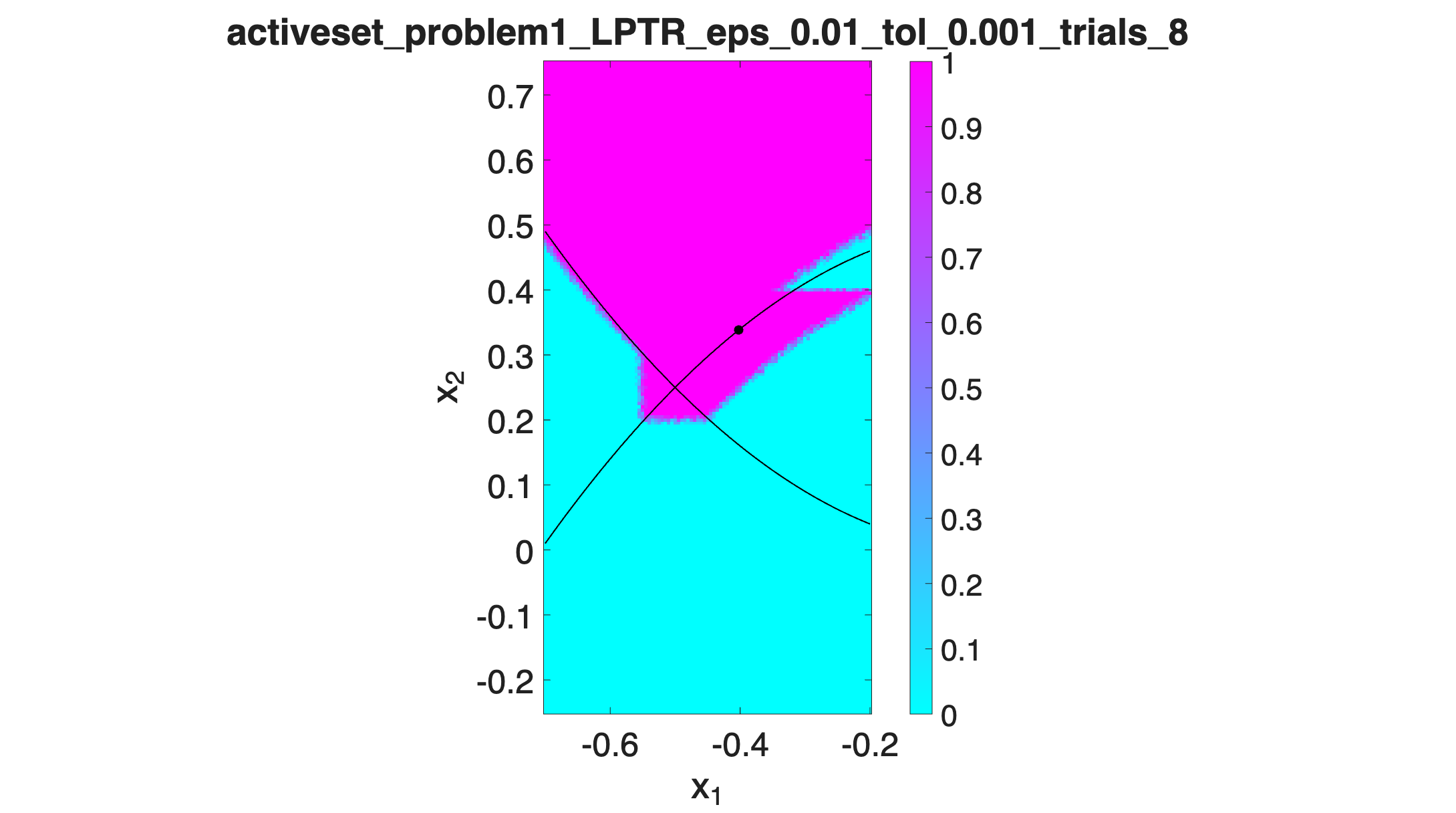}
    \hfill
    \includegraphics[width=0.2\linewidth, trim=310 0 340 38, clip]{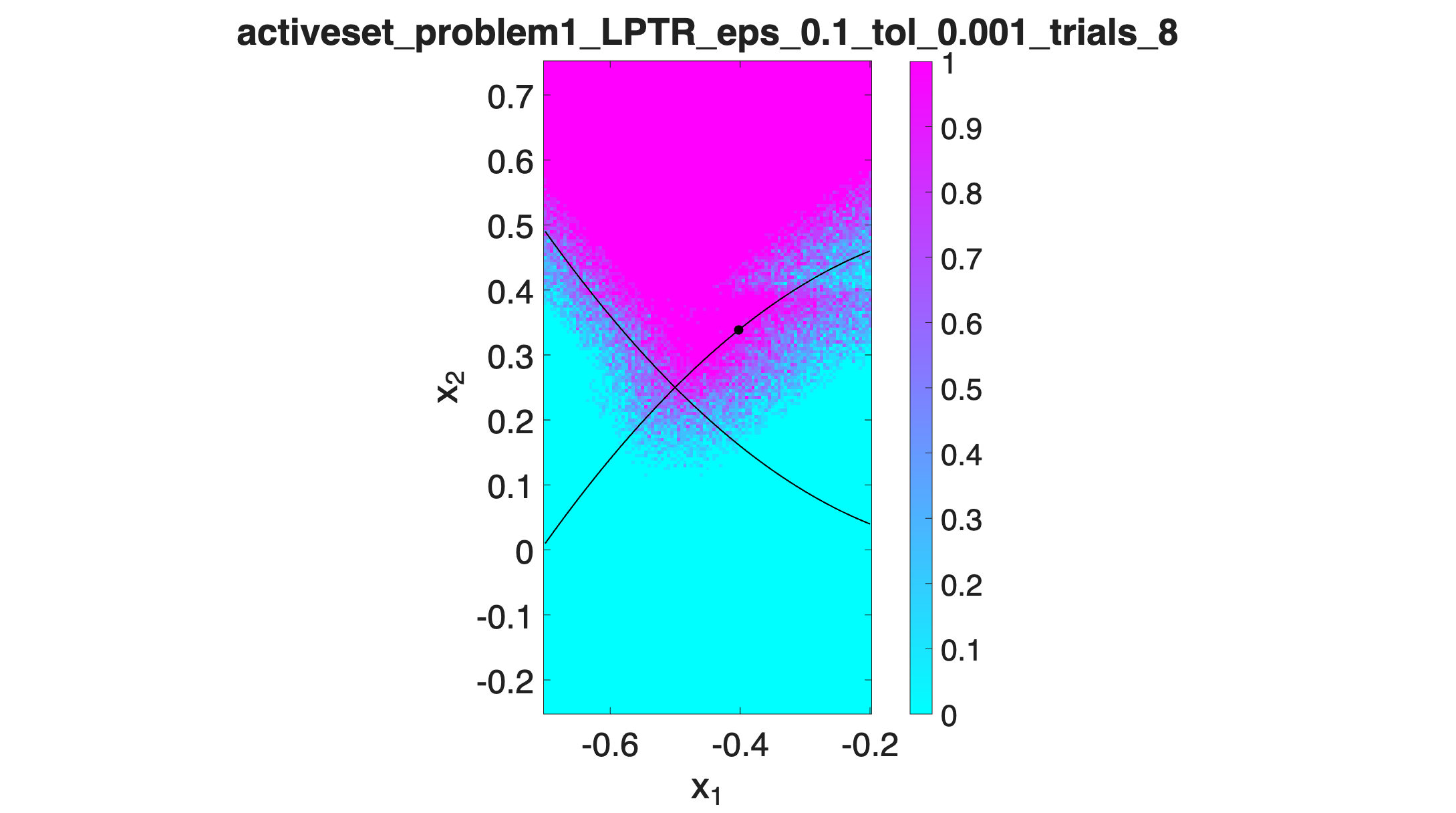}

    \begin{subfigure}[b]{0.2\linewidth}
        \centering
        \includegraphics[width=1\linewidth, trim=310 0 340 38, clip]{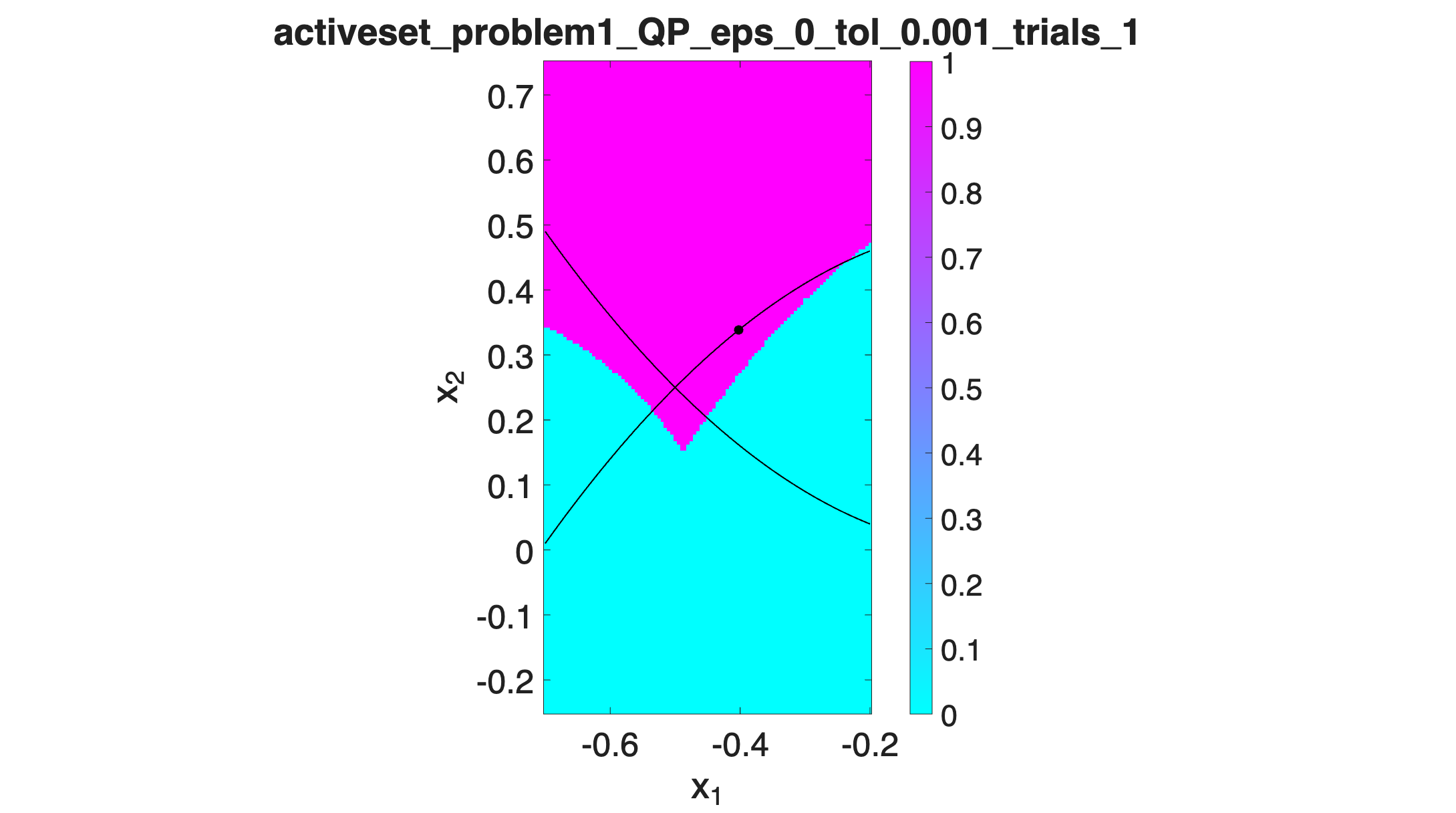}
        \caption*{(a) $\epsilon=0$}
    \end{subfigure}\hfill
    \begin{subfigure}[b]{0.2\linewidth}
        \centering
        \includegraphics[width=1\linewidth, trim=310 0 340 38, clip]{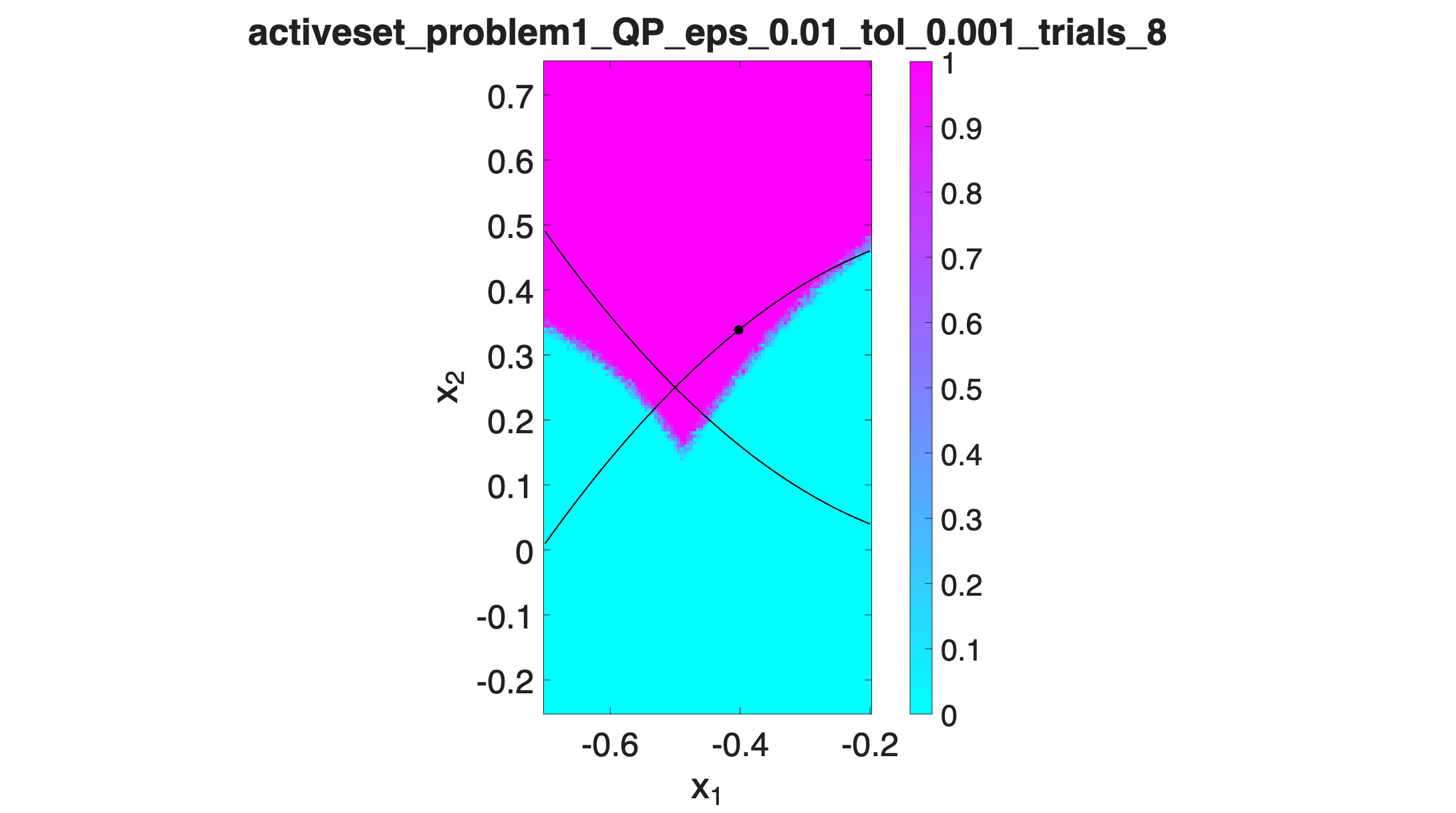}
        \caption*{(b) $\epsilon=10^{-2}$}
    \end{subfigure}\hfill
    \begin{subfigure}[b]{0.2\linewidth}
        \centering
        \includegraphics[width=1\linewidth, trim=310 0 340 38, clip]{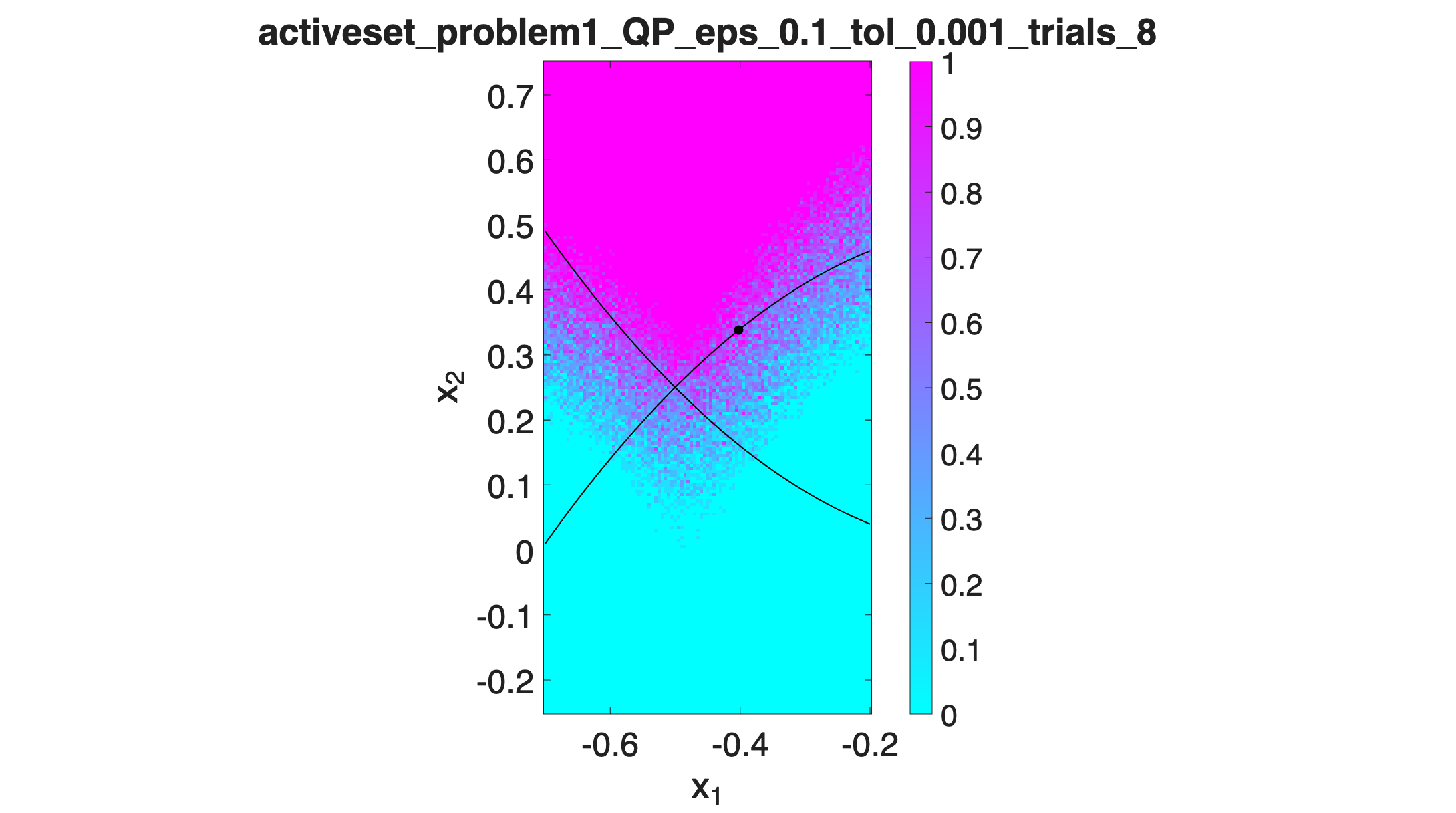}
        \caption*{(c) $\epsilon=10^{-1}$}
    \end{subfigure}
    \caption{Heatmaps for the problem to minimize $f_1$ subject to $c(x) \leq 0$ for noise levels $\epsilon = 0$ (deterministic), $10^{-2}$, and $10^{-1}$, respectively. The top row corresponds to the LP-LPEC approach based on computing Lagrange multiplier estimates, the middle row corresponds to the LP-TR approach based on computing a primal step, and the bottom row corresponds to the QP approach based on computing a primal step. The value of the heatmap at each point indicates the fraction of times that the active set $\Acal(x_*)$ was identified correctly by solving a subproblem at the point. (The unique optimal solution $x_*$ is indicated by a black dot.) In the case of positive noise levels, the value was determined as a fraction over 8 randomly generated subproblem instances.}
    \label{fig:heatmap_problem1_methods}
\end{figure}

\begin{figure}[ht]
    \centering
    \includegraphics[width=0.2\linewidth, trim=310 0 340 38, clip]{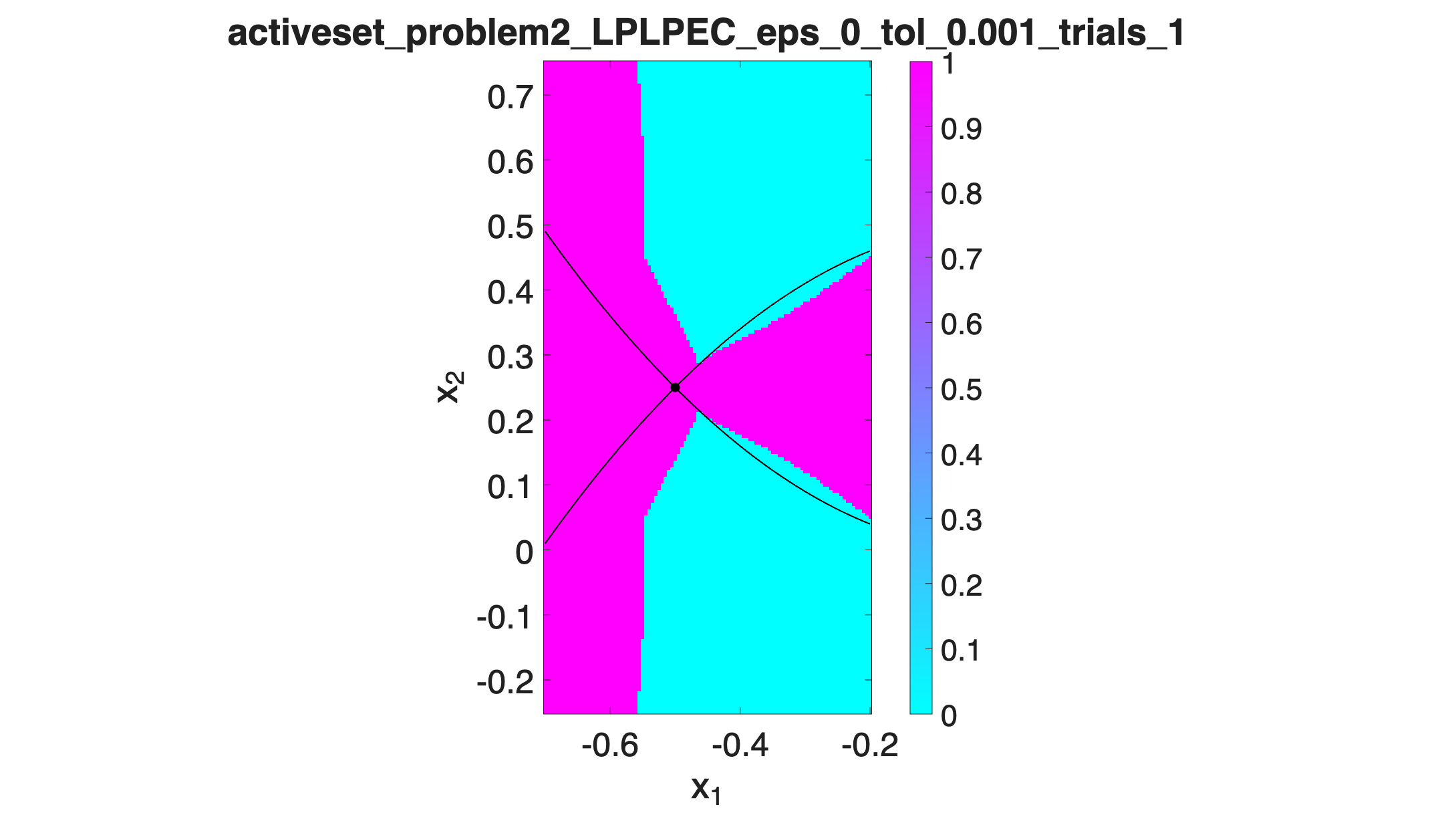}
    \hfill
    \includegraphics[width=0.2\linewidth, trim=310 0 340 38, clip]{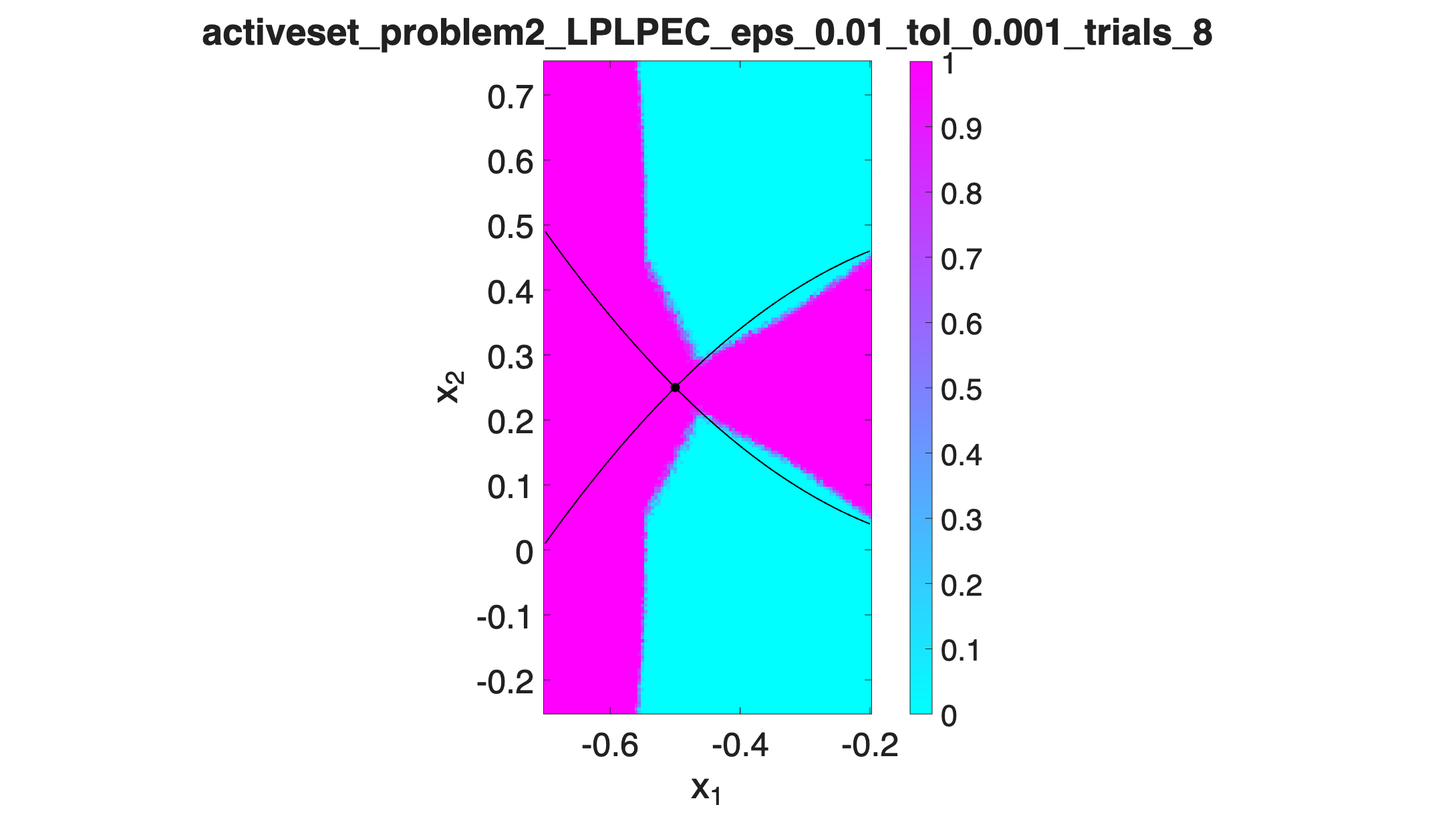}
    \hfill
    \includegraphics[width=0.2\linewidth, trim=310 0 340 38, clip]{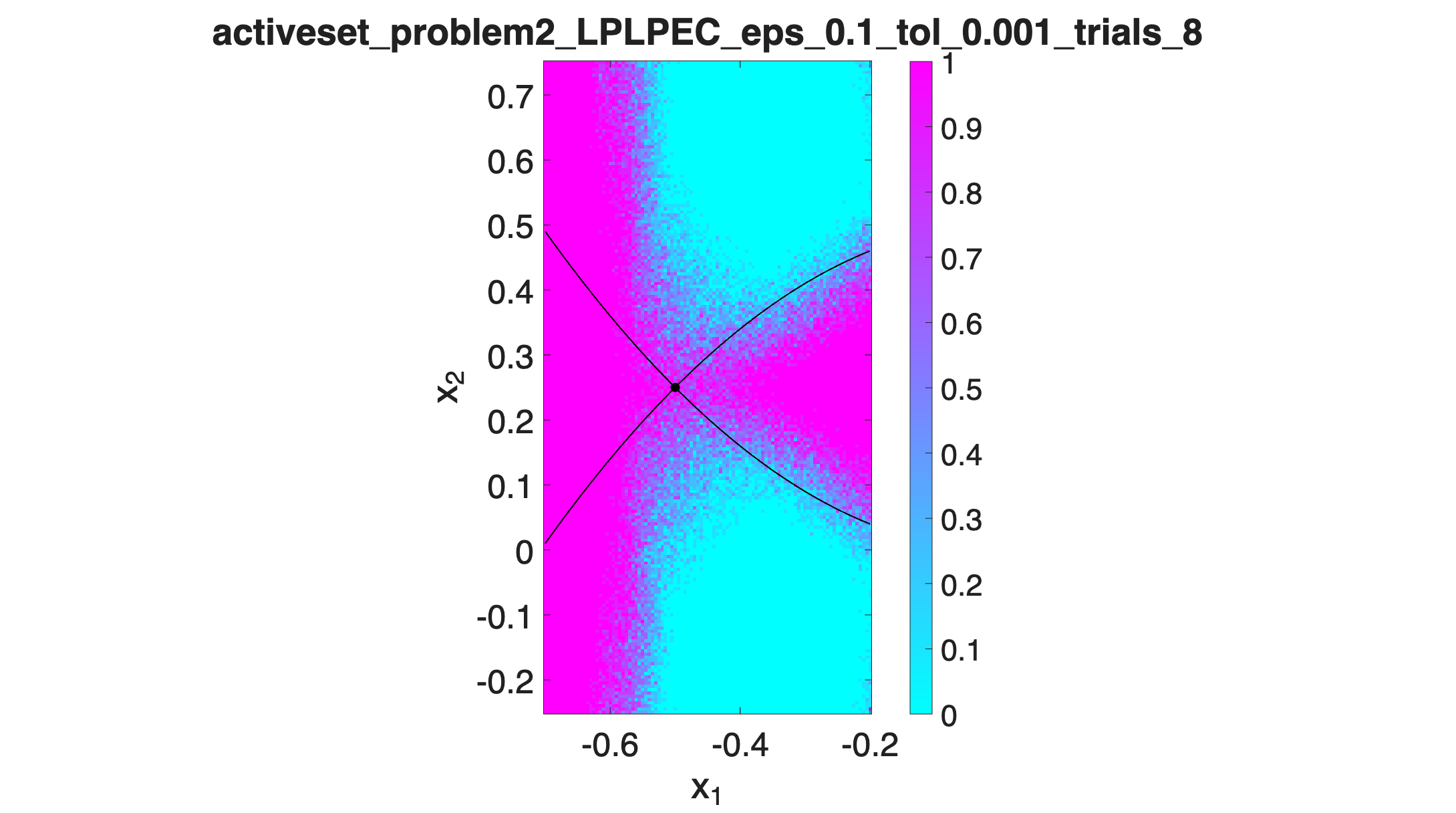}

    \includegraphics[width=0.2\linewidth, trim=310 0 340 38, clip]{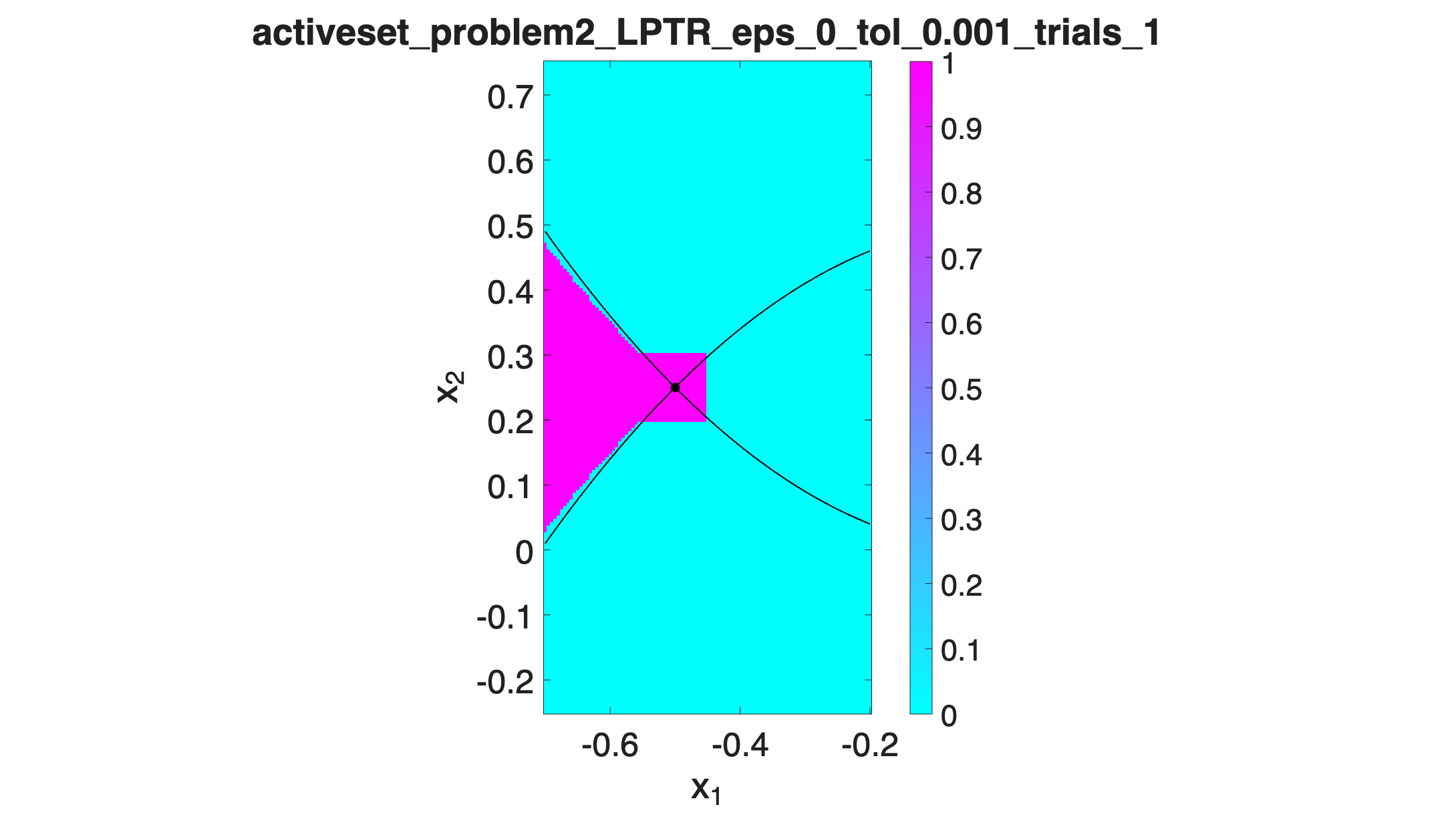}
    \hfill
    \includegraphics[width=0.2\linewidth, trim=310 0 340 38, clip]{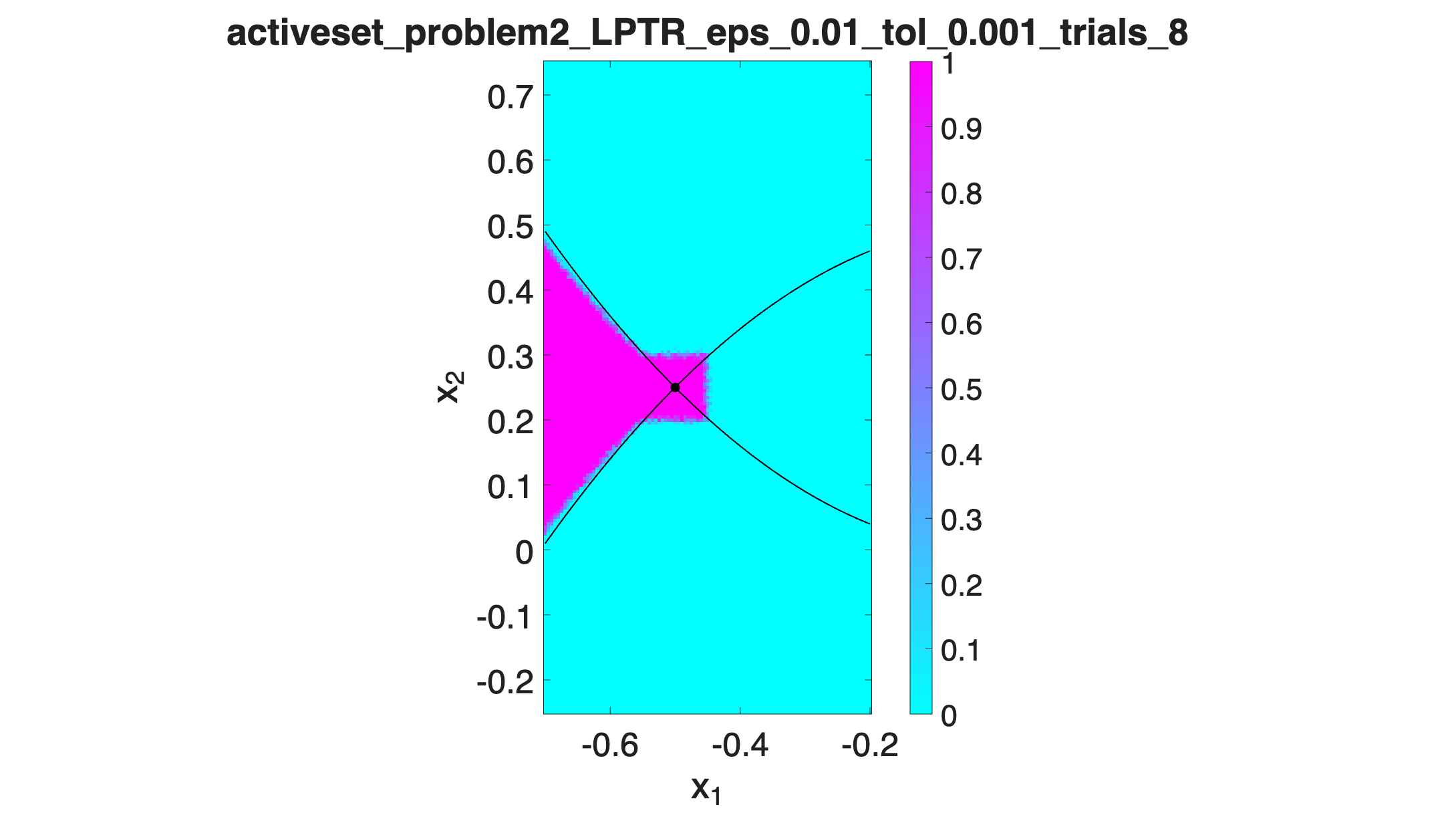}
    \hfill
    \includegraphics[width=0.2\linewidth, trim=310 0 340 38, clip]{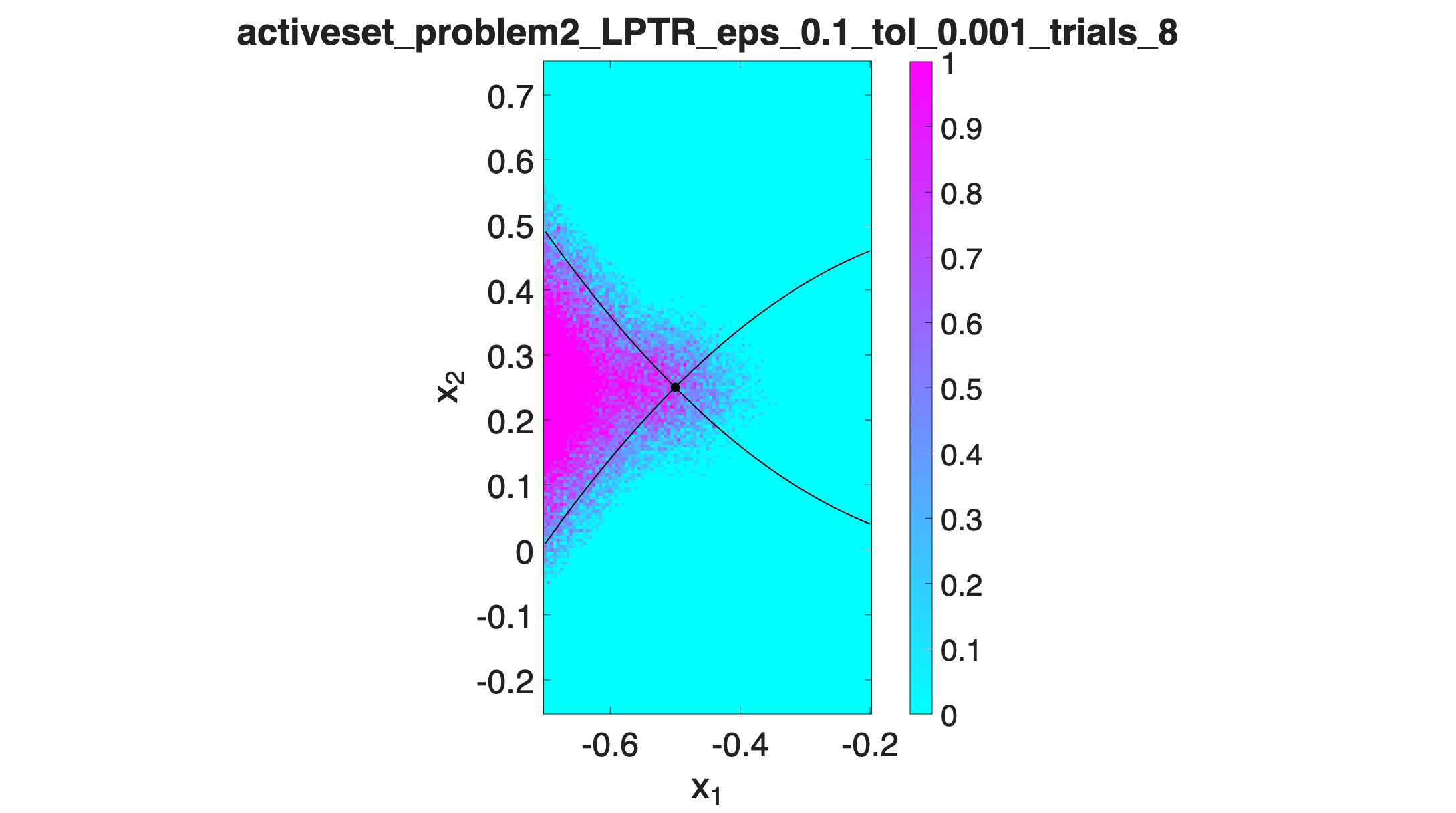}

    \begin{subfigure}[b]{0.2\linewidth}
        \centering
        \includegraphics[width=1\linewidth, trim=310 0 340 38, clip]{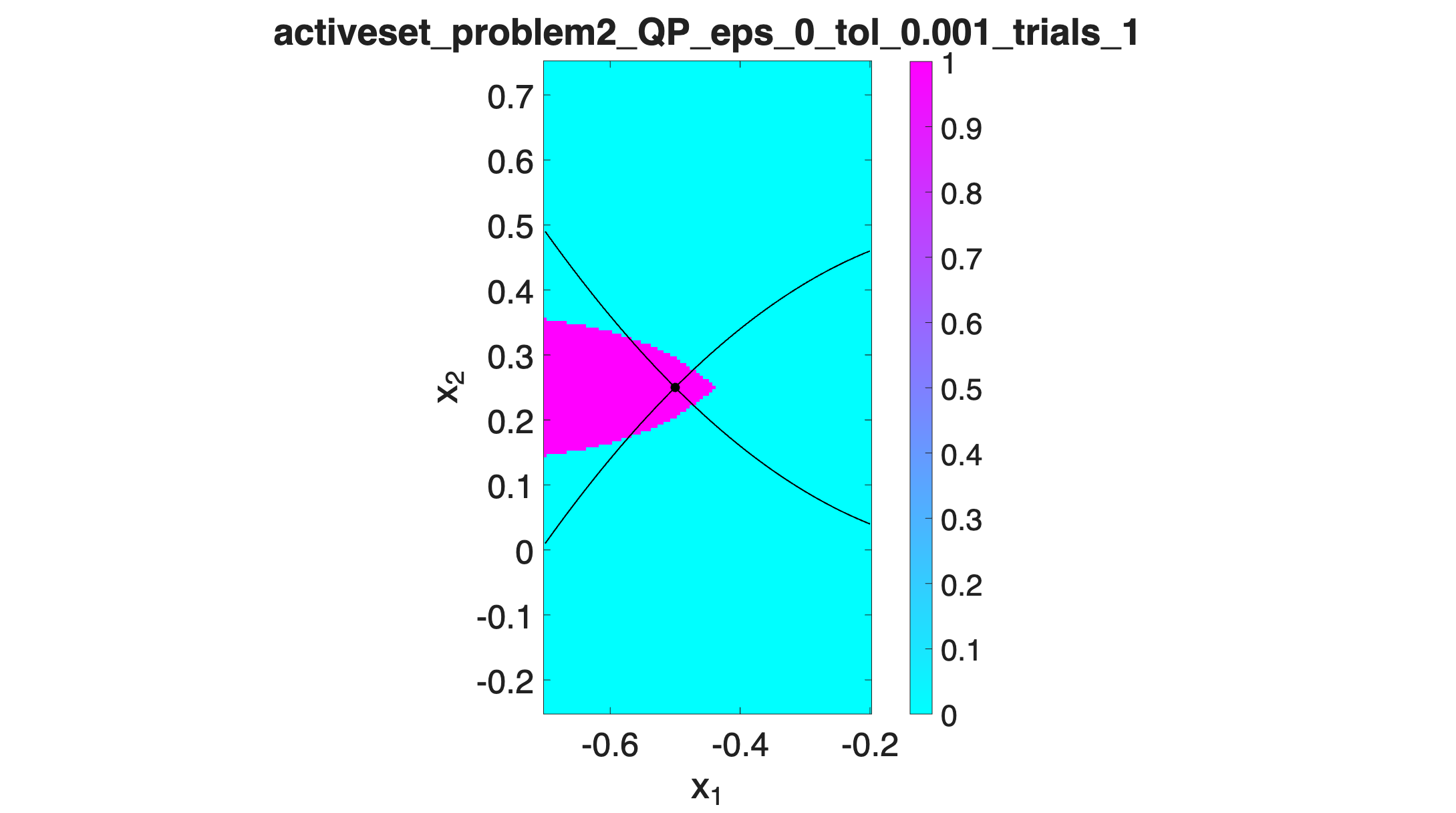}
        \caption*{(a) $\epsilon=0$}
    \end{subfigure}\hfill
    \begin{subfigure}[b]{0.2\linewidth}
        \centering
        \includegraphics[width=1\linewidth, trim=310 0 340 38, clip]{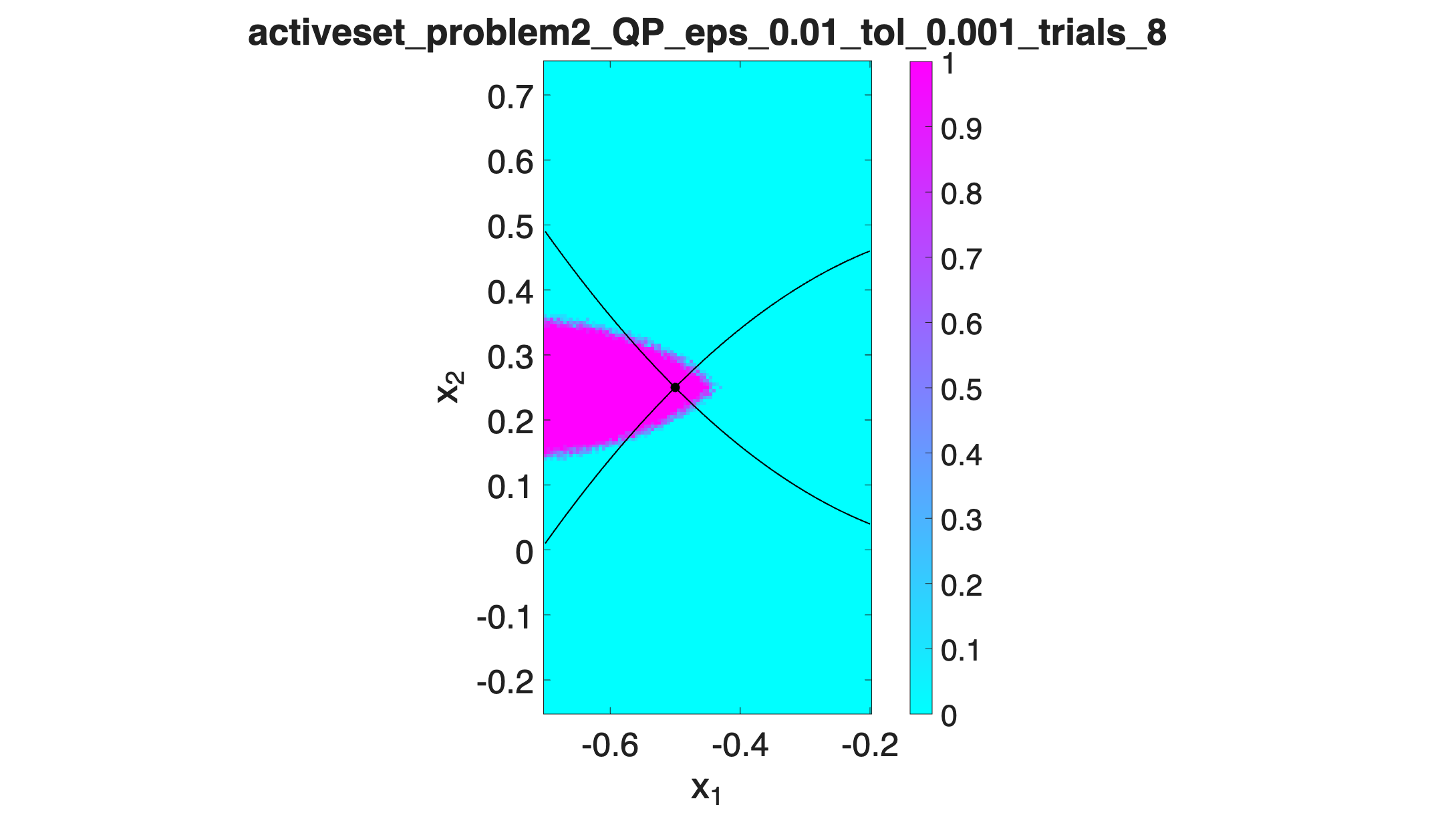}
        \caption*{(b) $\epsilon=10^{-2}$}
    \end{subfigure}\hfill
    \begin{subfigure}[b]{0.2\linewidth}
        \centering
        \includegraphics[width=1\linewidth, trim=310 0 340 38, clip]{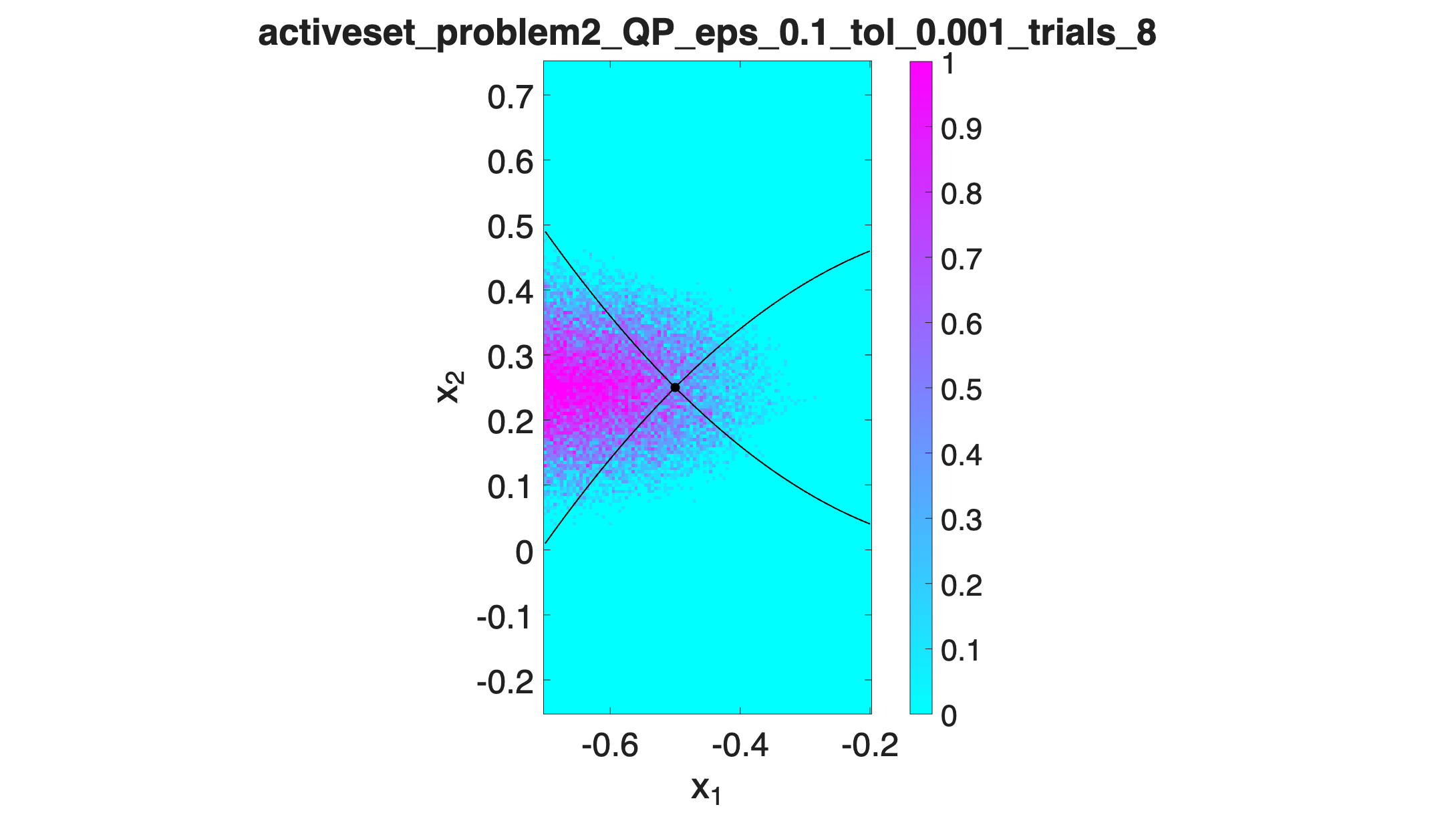}
        \caption*{(c) $\epsilon=10^{-1}$}
    \end{subfigure}
    \caption{Heatmaps for the problem to minimize $f_2$ subject to $c(x) \leq 0$.  The details are the same as in the caption for Figure~\ref{fig:heatmap_problem1_methods}.}
    \label{fig:heatmap_problem2_methods}
\end{figure}

The results associated with $f_1$ and $f_2$ are presented in Figures~\ref{fig:heatmap_problem1_methods} and \ref{fig:heatmap_problem2_methods}, respectively.  \edit{In each case, one finds that for noise level $10^{-2}$, the active-set identification at points near $x_*$ is relatively robust to noise.  For the LP-TR and QP approaches, this is consistent with Theorem~\ref{th.PertGra} and its LP counterpart, since the LICQ and strict complementarity hold at $x_*$ for both problems.  For the LP-LPEC approach, Corollary~\ref{cor.rate} guarantees correct identification only when the noise vanishes as $x \to x_*$, so the robustness observed near $x_*$ at a fixed noise level is an empirical finding rather than a consequence of our theory.}  On the other hand, for noise level $10^{-1}$, the active-set identification is less reliable for both methods, failing to identify $\Acal(x_*)$ in some cases.  This is also consistent with our theoretical guarantees, since our theory requires that the noise is sufficiently small in order to guarantee accurate active-set identification.  At least for these test problems, the LP-LPEC approach identifies the active set correctly less often than the LP-TR or QP approach when only one constraint is active at $x_*$, and identifies the active set correctly more often when both constraints are active, at least with our choices of the parameters $\beta$, $\sigma$, $\nu$, $\Delta$, and $\theta$.  Naturally, the behavior of the methods can vary with these parameter choices, so these comparisons do not necessarily hold in general.

A notable observation about the heatmaps in Figures~\ref{fig:heatmap_problem1_methods} and \ref{fig:heatmap_problem2_methods} is that, when the noise is nonzero, one should only assume that the errors in the function and derivative values are small with some probability such that active-set identification can be performed accurately, even within a small neighborhood of $x_*$.  This is an empirical demonstration of the comments that we made in \S\ref{sec.stochastic_algorithm}.

\subsection{Neural Network Training Problem}\label{sec.neural_network}

Now let us consider a constrained neural network training problem.  The problem that we consider is one of supervised learning of a prediction function that aims to predict, based on feature data for a given patient, whether the patient is at risk for heart disease.  The constraints are designed to ensure that the prediction function yields a higher probability of heart disease when the patient has a known high-risk factor.  We emphasize that the constraints that we impose are for demonstration purposes only, and may not be ones that lead to better predictions in reality.  That being said, the problem that we consider in this section might represent the type of constrained supervised learning problem that could potentially lead to better and more explainable prediction functions in the context of medical diagnoses, as opposed to prediction functions obtained purely through empirical loss minimization.

Let $\{(a_i,b_i)\}_{i=1}^N$ be a set of training data, where for each $i \in [N]$ the value $a_i \in \R{n_f}$ is an $n_f$-vector of feature data for a given patient (including age, sex, and cholesterol level)~and the value $b_i \in \{0,1\}$ indicates whether the patient eventually developed heart disease (indicated by $b_i = 1$) or not (indicated by $b_i = 0$).  In particular, let us suppose that for each $i \in [N]$ the $j$th component of $a_i$, namely, $a_{i,j}$, indicates the patient's cholesterol level.  For the purposes of devising a tractable problem, let us further suppose that the~$N$ data points are separated into~$K$ groups, where the indices in $\Ncal_1$ correspond to group~1, the indices in $\Ncal_2$ correspond to group~2, etc.  (In this manner, one has that $\{\Ncal_1, \Ncal_2, \dots, \Ncal_K\}$ are disjoint sets and $|\Ncal_1| + |\Ncal_2| + \cdots |\Ncal_K| = N$.  In practice, the groups might be determined randomly or according to certain features, such as age and/or sex.)  Given a prediction function $p$ parameterized by a vector $x \in \R{n}$ of trainable parameters---in our setting, defined by a neural network---and a loss function $\ell$, our problem is
\bequation\label{prob.neural_network}
  \min_{x \in \R{n}}\ \sum_{i=1}^N \ell(p(x,a_i),b_i)\ \st\ \max_{i \in \Ncal_k} \left\{-\tfrac{\partial p}{\partial a_{\cdot,j}} (x,a_i) \right\} \leq 0\ \text{for all}\ k \in [K],
\eequation
where $\partial p / \partial a_{\cdot,j}$ indicates the partial derivative of $p$ with respect to the $j$th element of the input of the neural network, which in our setting corresponds to cholesterol.  Here, the objective function is a typical one for empirical loss minimization that aims to maximize the accuracy of predictions on unseen data.  The constraints, on the other hand, aim to ensure that---all else being equal---an increase in a patient's cholesterol would increase (or at least not decrease) the prediction that the patient will develop heart disease.  Rather than impose a constraint for each data point, which may lead to an intractable problem in large-scale settings, a constraint is imposed for each group in the data set.  It should be noted that the maximum function is nonsmooth, and in practice one could impose a smooth constraint through a smoothed approximation of the maximum function.  However, for our experiments we did not smooth the functions since (a) we were able to solve the problem to sufficient accuracy despite nonsmoothness in the constraint functions and (b) it is often the case in modern practice that algorithms employed for neural network training ignore nonsmoothness due to maximum functions, as with rectified linear unit (ReLU) activation.  \edit{We also speculate that nonsmoothness of this type might actually help identification rather than hurt it.  After all, a max-function constraint is sharp at points where the maximum is attained uniquely, and sharpness is precisely the property that drives finite identification in nonsmooth analysis and optimization theory; see, e.g., \cite{BurkMore1988,Lewi2002,HareLewi2004}.}

We acquired a data set of 1025 data points with each point having $n_f = 13$ features~\cite{heart_disease_45}.  For the prediction function $p$, we used a neural network with one hidden layer with 6 nodes and tanh activation.  Including both neural network edge weights and bias terms as trainable parameters, this led to a problem with 91 variables.  For the loss function we used the log-loss function.  The data was split initially into training and testing data according to cholesterol level; in particular, we sorted the data points by cholesterol and used the top 20\% according to cholesterol as the testing data.  In this manner, training was conducted on those with lower cholesterol and testing was conducted on those with higher cholesterol.  As an interesting side-note to our study of the behavior of active-set identification strategies, the idea here was to test the effect of training without the highest-cholesterol individuals in order to see if the presence of the constraints in \eqref{prob.neural_network} was able to improve the predictions on the testing data.

The training data was split into 10 groups to define~10 constraints.  We began our experimentation by employing an implementation of an \mbox{L-BFGS} method with a backtracking Armijo line search to minimize the objective of \eqref{prob.neural_network}, without any constraints.  This resulted in an approximate stationary point---we used a threshold of $10^{-5}$ for the norm of the gradient---with a maximum constraint function value of \texttt{6.03e-03}; i.e., it resulted in a point that was infeasible for the constraints in \eqref{prob.neural_network}.  Corresponding to this point, the training accuracy was determined to be 88.05\% while the testing accuracy was determined to be 81.46\%.

We then employed the same L-BFGS method to minimize a quadratic penalty function corresponding to \eqref{prob.neural_network} with a penalty parameter of \texttt{1e+03}; see \cite{NoceWrig2006} for further details on penalty methods.  With a tolerance of $10^{-5}$ for the norm of the gradient (of the quadratic penalty function), this resulted in a point with a maximum constraint value of \texttt{8.57e-09}.  Corresponding to this point, the training accuracy was determined to be 92.20\% while the testing accuracy was determined to be 84.88\%.  This demonstrates a case in which the presence of the constraints during the training process was able to improve testing accuracy by over 3\%.

\edit{We emphasize that this experiment evaluates the three identification rules along a sequence of iterates produced by an algorithm that does not itself use active-set information. As such, it is not a demonstration of an algorithm that computes search directions from an estimated active set, which is a topic for future work (see \S\ref{sec.conclusion}).  Denote by $\hat x$ the final iterate of the L-BFGS method applied to the quadratic penalty function, which we treat as a solution estimate.  At $\hat x$, in order to determine the estimated optimal active set, we employed the active-set identification techniques from \S\ref{sec.noisy} to solve an LP-LPEC subproblem, an LP-TR subproblem with a trust-region constraint, and a QP subproblem.  Here and throughout the remainder of these experiments, we employed the parameters $M \gets 10^8$, $\beta \gets 10^{-10}$, $\sigma \gets 0.7$, $\nu \gets 10^2$, $\Delta \gets 10^{-2}$, and $\theta \gets 10^2$.  Solving all of these subproblems led to the same estimate, $\{1,\dots,10\}$, i.e., all constraints being active.  We take this common estimate as the reference active set in what follows. That is, we do not claim that $\hat x$ is a local minimizer of \eqref{prob.neural_network}, nor that the reference set is provably $\Acal(x_*)$ for a nearby minimizer.}

\begin{figure}[ht]
  \centering
  \includegraphics[width=0.325\linewidth]{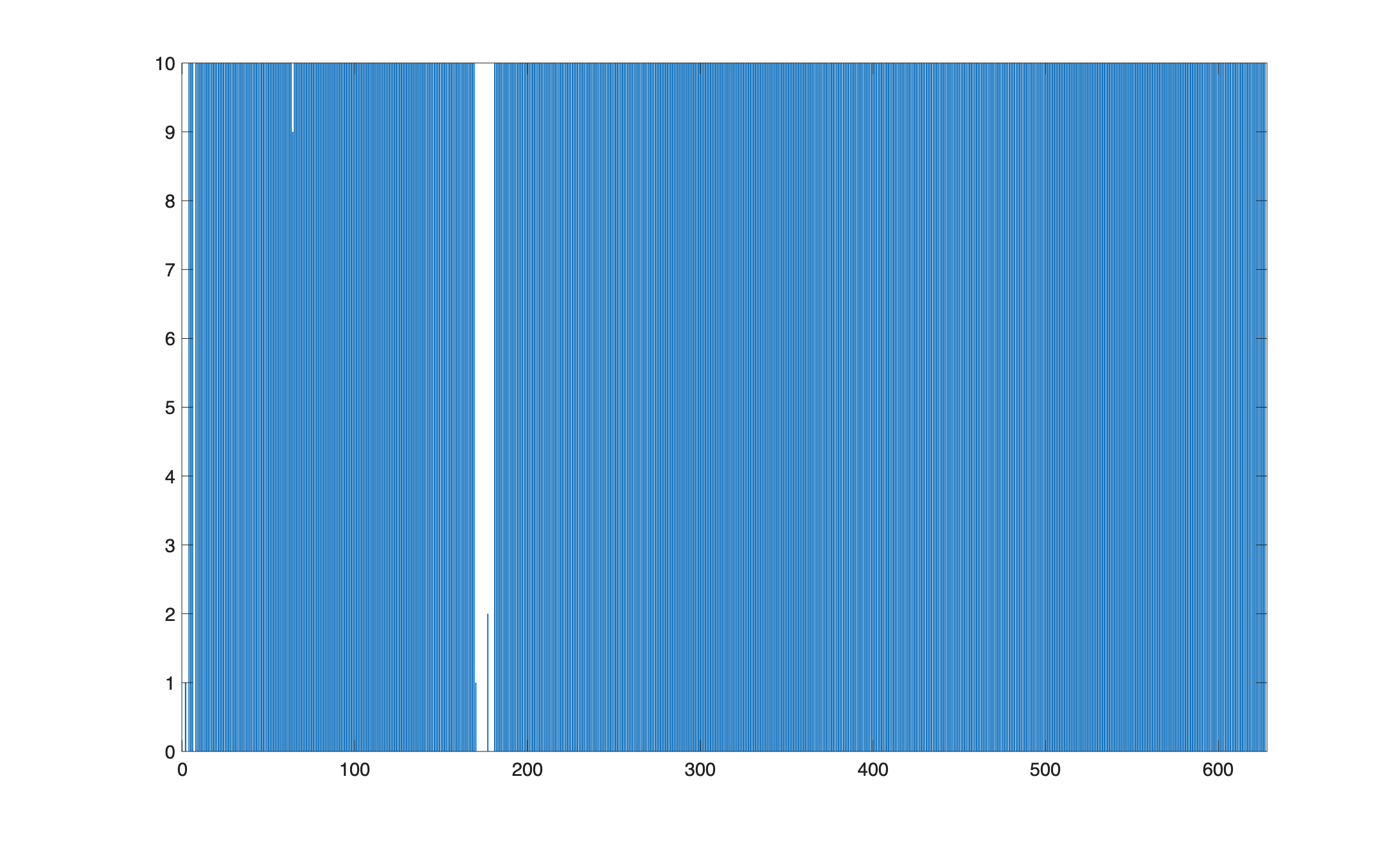}\ 
  \includegraphics[width=0.325\linewidth]{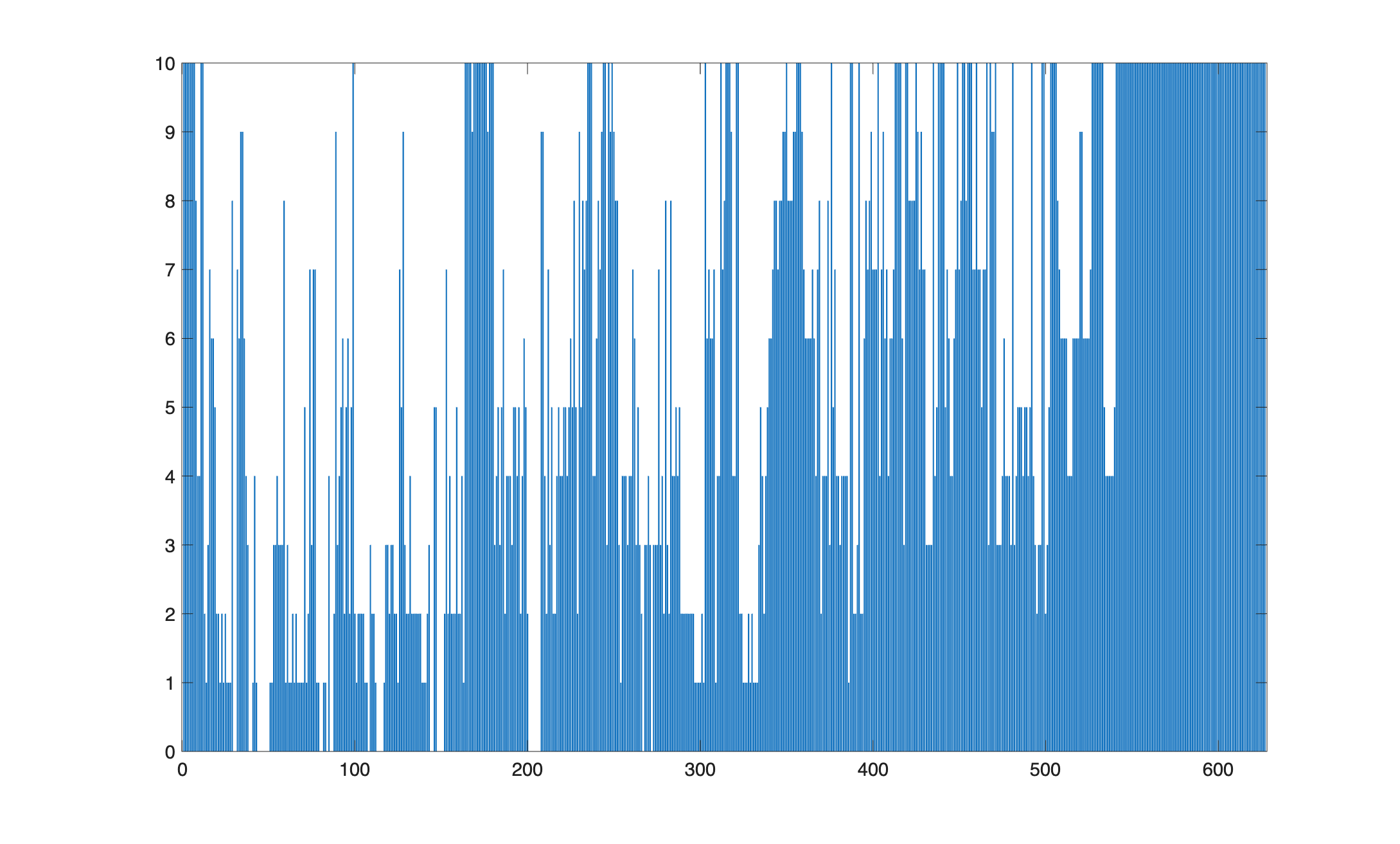}\ 
  \includegraphics[width=0.325\linewidth]{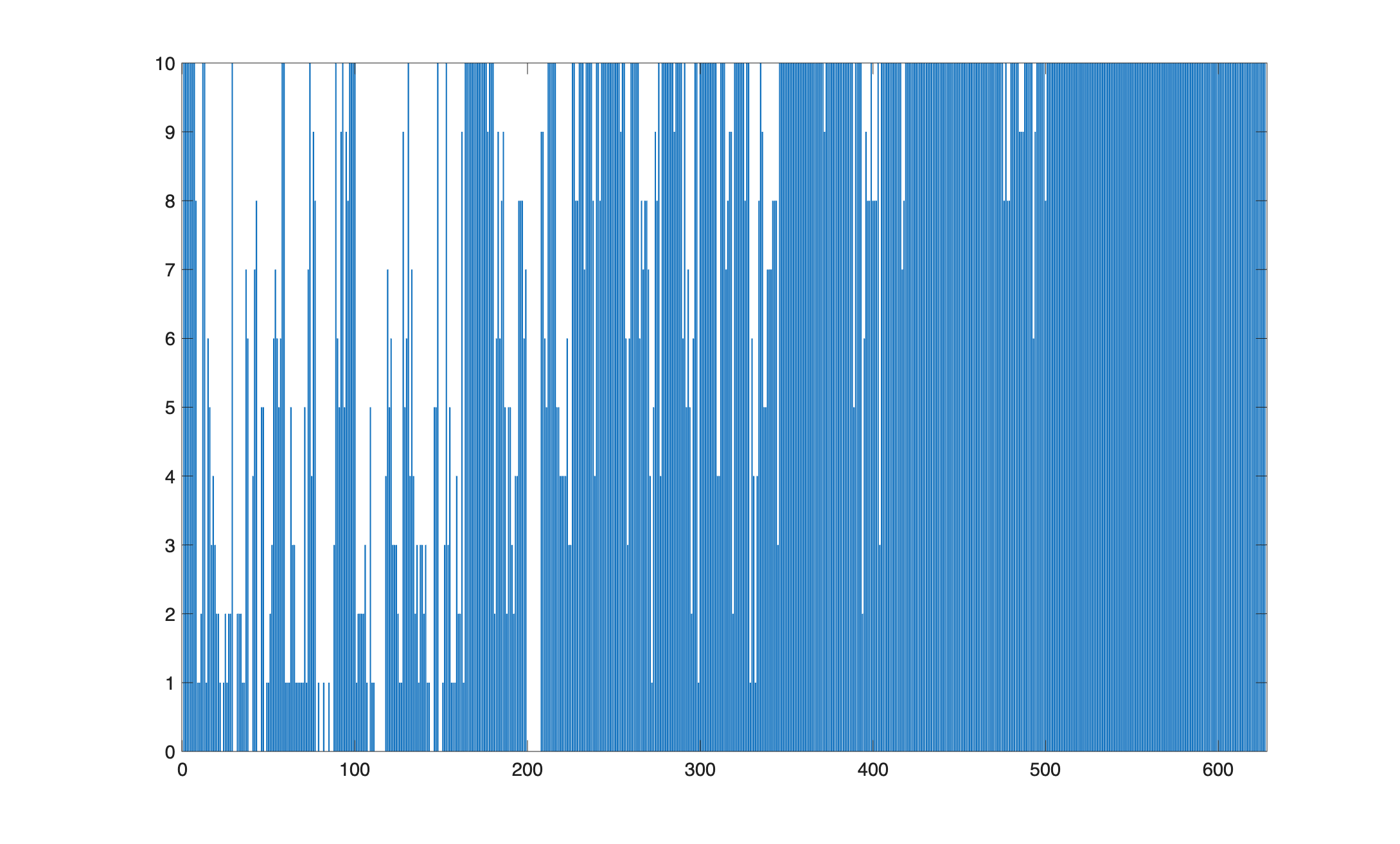}
  \caption{Over the iterates of an algorithm to solve \eqref{prob.neural_network}, the number of correctly identified elements of the optimal active set by the LP-LPEC (left), LP-TR (middle), and QP (right) approach when there \emph{is no noise} in the function and derivative values.}
  \label{fig.nn_determinstic}
\end{figure}

We studied the behavior of the active-set identification techniques from \S\ref{sec.noisy} by following the sequence of iterates in the L-BFGS method for minimizing a quadratic penalty function, as described previously.  At each iterate, Figure~\ref{fig.nn_determinstic} shows the number of correctly identified active constraints by the LP-LPEC, LP-TR, and QP approaches when there \emph{is no noise} in the function and derivative values used to define the subproblems for each approach.  One finds that the LP-LPEC approach identifies the correct active set even in early iterations.  We found this to be due to the fact that the iterates are consistently points where all constraints are violated.  Combined with the fact that the right-hand side of the inequality in \eqref{eq.active_LP} is always negative, this demonstrates a situation in which the LP-LPEC approach is lucky always to identify violated constraints as being active at the optimal solution.  The LP-TR and QP approaches, on the other hand, do not identify the optimal active set until the iterates are closer to~$\xhat$.  That said, it should be emphasized that all approaches correctly identify the optimal active set in all later iterations.

\begin{figure}[ht]
  \centering
  \includegraphics[width=0.325\linewidth]{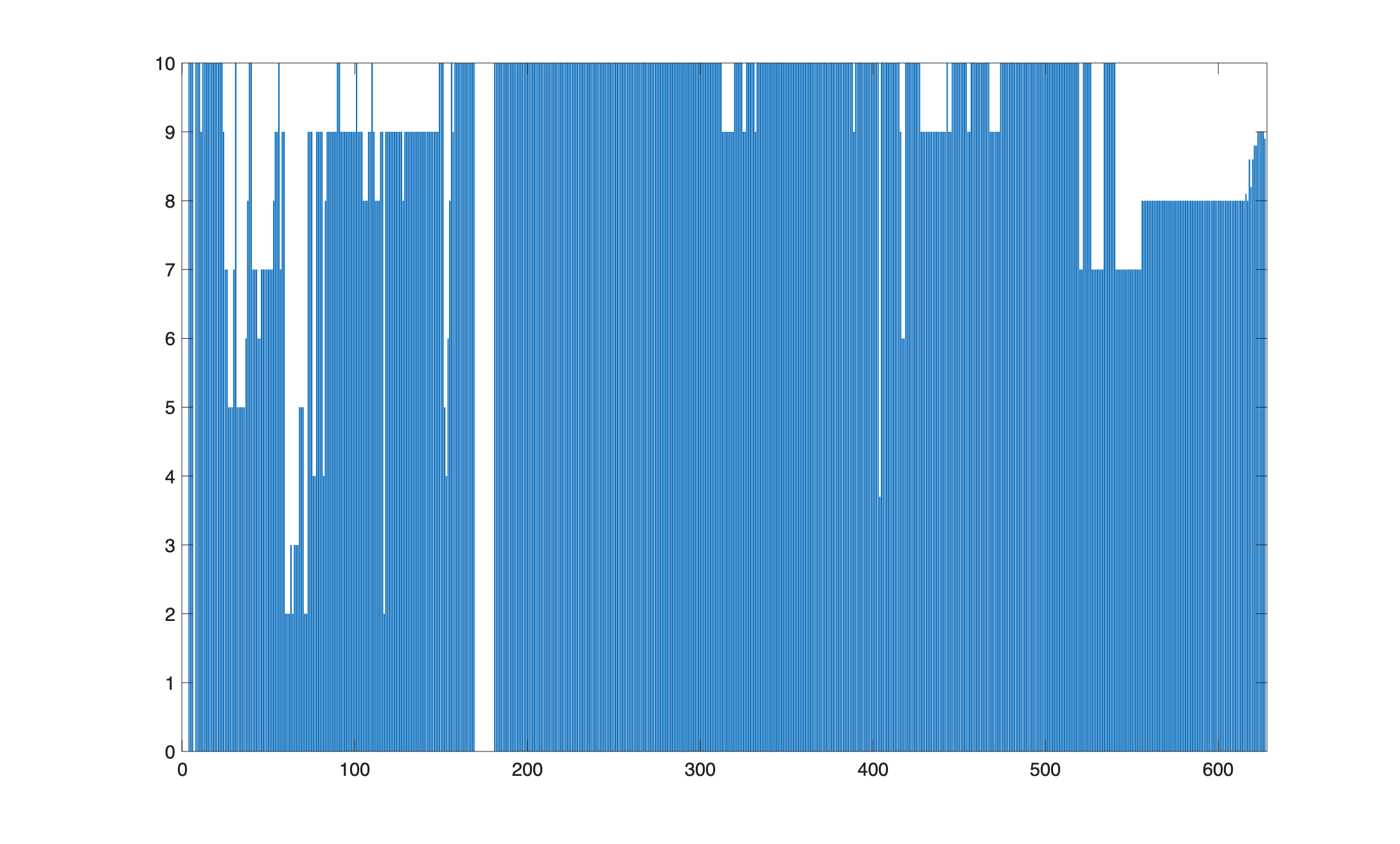}\ 
  \includegraphics[width=0.325\linewidth]{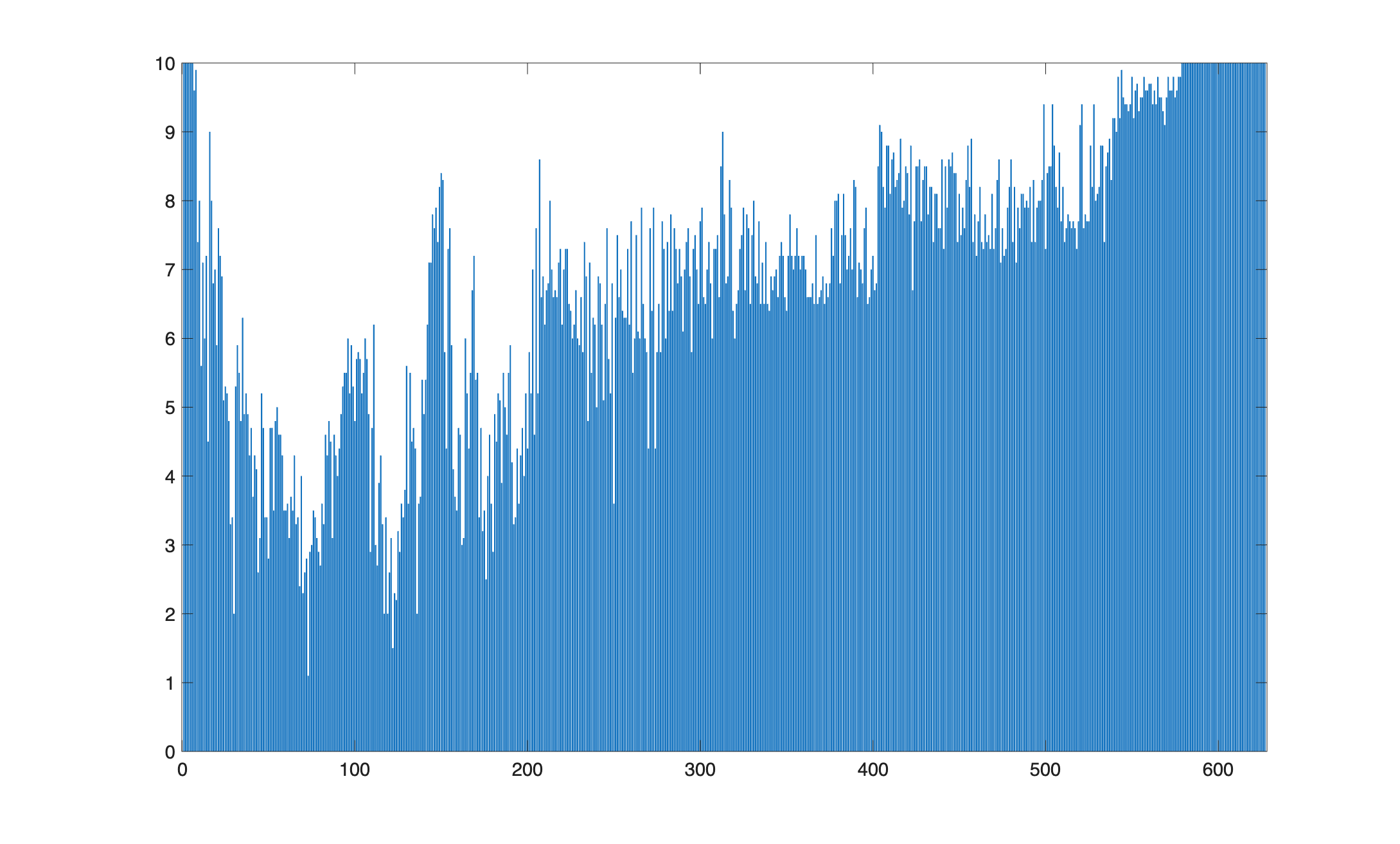}\ 
  \includegraphics[width=0.325\linewidth]{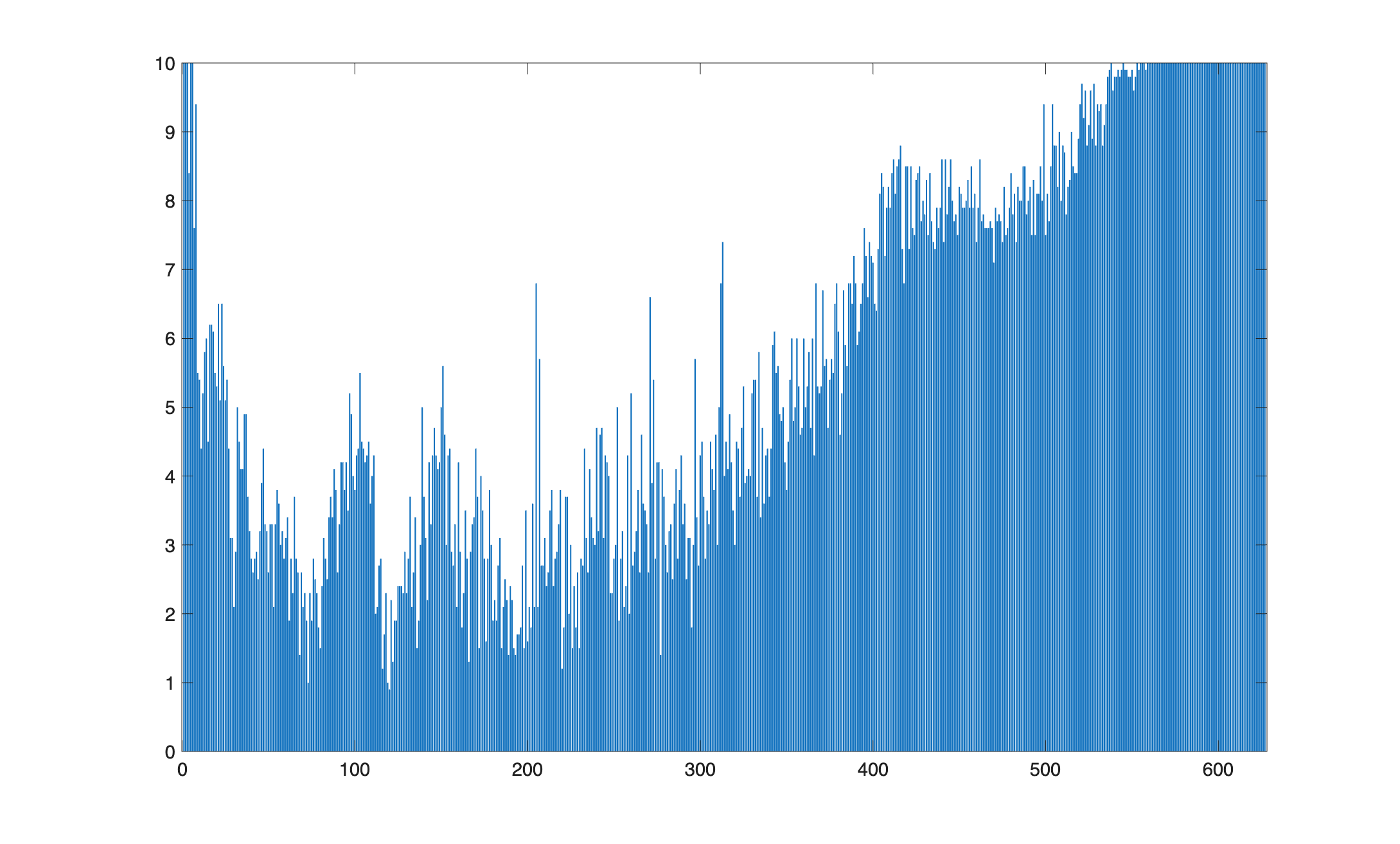}
  \caption{Over the iterates of an algorithm to solve \eqref{prob.neural_network}, the number of correctly identified elements of the optimal active set by the LP-LPEC (left), LP-TR (middle), and QP (right) approach when there \emph{is noise} in the function and derivative values.}
  \label{fig.nn_stochastic}
\end{figure}

To observe the effect of noise in the function and derivative values on active-set identification, we ran the algorithm with exact function and derivative values, but at each iterate we generated 10 different sets of noisy objective and constraint function and derivative values for active-set identification by sampling the data.  In particular, at each iterate, we used a random sample of data points with sample size equal to 32 in order to evaluate the function and derivative values employed in the subproblems for active-set identification.  (We ran the code with other similar sample sizes and the results were qualitatively similar to those presented here.)  The results are shown in Figure~\ref{fig.nn_stochastic}, where for each iteration index we plot the average number of correctly identified active constraints over the 10 different samples that were generated.  One finds that the LP-LPEC approach does a poorer job (compared to the results in Figure~\ref{fig.nn_determinstic}) of identifying the optimal active set in earlier iterations, but it still identifies most of the optimal active set correctly in later iterations.  The LP-TR and QP approaches, on the other hand, happen to do a better job of identifying the active set.

As for the experiments in \S\ref{sec.two-dimensional}, the behavior of the methods can vary with different problems and parameter choices, so these  specific comparisons do not necessarily hold in general.  That said, the important thing is that we observe accurate active-set identification in later iterations for all of the subproblem techniques, even with noise in the problem function and derivative values.

\section{Conclusion}\label{sec.conclusion}

We have shown that state-of-the-art techniques for active-set identification when solving constrained continuous optimization problems can also be employed in contexts when function and derivative values are corrupted by noise\edit{. We have also quantified the noise that can be tolerated.  Two of our conclusions deserve explicit emphasis.  First, the admissible noise level can be made explicit, as seen in our analysis. Under the LICQ and strict complementarity it is proportional to the product of the conditioning of the active constraint Jacobian with the degeneracy margins of the minimizer (Lemma~\ref{lem.explicit}), and under the MFCQ and the second-order condition \eqref{eq:second_order_cond} it suffices that the noise vanish linearly in the distance to the minimizer (Corollary~\ref{cor.rate}).  Second, the contrast between these regimes is intrinsic rather than an artifact of the strategies analyzed. In particular, at a minimizer with a weakly active constraint, no identification rule can succeed under a fixed noise level, and the noise must in general vanish at least linearly, so the rate of Corollary~\ref{cor.rate} cannot be improved (Theorem~\ref{th.impossible}). Accordingly, the results of this paper may be read both as delimiting what active-set identification can achieve in the presence of noise and as a contribution to the perturbation theory of nonlinear programming, in which the perturbed object is a computable subproblem at a point that is merely near a solution}.  We have also shown how such techniques can be employed to offer active-set identification guarantees in the context of a stochastic algorithm for solving constrained continuous optimization problems.  We have also provided the results of numerical experiments that demonstrate the behaviors of these techniques on a pair of illustrative problems and a  neural-network training problem.

Interesting follow-up investigations to our work in this paper include studying how active-set identification techniques can be employed within noisy and stochastic algorithms for solving constrained continuous optimization problems both for computing certain types of primal search directions and for more accurate Lagrange multiplier estimation.  For example, one may consider whether, despite noise in problem function and derivative values, one can design algorithms with convergence guarantees that operate by (a) estimating the active set, then (b) computing each search direction with respect to the estimated active set, as has been done in the deterministic setting \cite{ByrdGoulNoceWalt2003,ByrdGoulNoceWalt2005}.  It would also be interesting to investigate how, e.g., the theoretical guarantees for Lagrange multiplier estimation presented in~\cite{CurtJianWang2024} (for the case of equality constraints only) can be extended to the setting of an algorithm for solving inequality-constrained problems.  This may involve employing the active-set identification techniques and theoretical guarantees from this paper in order to provide guarantees for optimal Lagrange multiplier estimation.

\appendix
\section{Proof of Theorem~\ref{th.OberWrig}}\label{app.OberWrig}

\bproof{Proof.}
  Consider part (a).  First, by continuity of $c$, it follows that one can choose $\epsilon \in \R{}_{>0}$ small enough such that for any $x \in \R{n}$ with $\|x - x_*\|_2 \leq \epsilon$ and any $i\notin \Acal(x_*)$ one has that $c_i(x)\leq \tfrac{1}{2}c_i(x_*) < 0$.  Second, given such an $\epsilon$ and for any such $x$ and positive $\nu$, it follows from \eqref{eq.FEEL} that $\|d_{x,\theta,\nu}\|_2 \to 0$ as $\theta \nearrow \infty$.  Thus, one can choose $\underline\theta$ large enough such that $i\notin\Acal(x_*)$ implies $c_i(x) + \nabla c_i(x)^Td_{x,\theta,\nu}<0$.  Thus, the claim in part (a) follows from the definition \eqref{eq.active_QP}.

  Next, consider part (b).  Let us first suppose that $(\epsilon,\underline\theta)$ satisfies the conditions of part (a).  Furthermore, to prove the claim, let us introduce
  \bequation\label{eq.def_nu}
    \underline{\nu} := \max\{ \max\{\| y_* \|_\infty, \| z_* \|_\infty\} + 1 : (y_*, z_*) \in \Dcal(x_*) \}.
  \eequation
  The value $\underline\nu$ is well-defined and positive since $x_*$ is a primal KKT point and the MFCQ assumption ensures that $\Dcal(x_*)$ is bounded. For any pair $(y_*,z_*)\in\Dcal(x_*)$, the dual subproblem~\eqref{eq.sub_dual} at $x_*$ has an objective value of 0, as previously stated. Otherwise, for the dual subproblem at $x \neq x_*$ with $\|x - x_*\|_2 \leq \epsilon$ and $\nu \geq \bar\nu$, a feasible point is given by $(\alpha,\beta) = (y_*,z_*)$, for which
  \begin{align*}
      &e(x)^Ty_* + c(x)^Tz_* - \tfrac{1}{2} \theta^{-1} \|\nabla f(x) + \nabla e(x) y_* + \nabla c(x) z_*\|_2^2\\
      &=(e(x) - e(x_*))^Ty_* + (c(x) - c(x_*))^Tz_* \\
      &\quad - \tfrac{1}{2} \theta^{-1} \|\nabla f(x)-\nabla f(x_*) + (\nabla e(x)-\nabla e(x^*)) y_* + (\nabla c(x)-\nabla c(x_*)) z_*\|_2^2\\
      &\geq -L(\|x-x_*\|_2^2+\|x-x_*\|_2)
  \end{align*}
  for a sufficiently large constant $L \in \R{}_{>0}$.  This latter inequality follows due to the continuous differentiability of the problem functions and $\|x - x_*\|_2 \leq \epsilon$.
  
  Now let us assume, for the sake of deriving a contradiction, that no matter how small one chooses $\epsilon \in \R{}_{>0}$, there exists a point $x \in \R{n}$ satisfying $\|x - x_*\|_2 \leq \epsilon$ such that the resulting active-set estimate $\Acal_{{\rm QP}}(d_{x,\theta,\nu}; x)$ yields $-\nabla f(x_*)\notin \text{range}[\nabla e(x_*)]+\text{pos}[\{\nabla c_i(x_*) : i \in \Acal_{{\rm QP}}(d_{x,\theta,\nu}; x)\}]$. Since there are only finitely many active set estimates (i.e., subsets of $[q]$), it follows that there exists a sequence of points converging to $x_*$ such that along this sequence the active-set estimate is the same set of indices, call it $\widehat\Acal$.  At the same time, observe that by Rockafellar \cite[Theorem 19.1]{ROCKAFELLAR1970} the set $\text{range}[\nabla e(x_*)]+\text{pos}[\{\nabla c_i(x_*) : i\in\widehat\Acal\}]$ is finitely generated, meaning that it is a closed set.  Now, by the definition of the distance function in equation \eqref{eq.def_dist}, one can define the quantity $\tau$ as
  \bequation\label{eq.def_tau}
    \tau := \thalf \dist(-\nabla f(x_*), \text{range}[\nabla e(x_*)] + \text{pos}[\{\nabla c_i(x_*) : i\in\widehat\Acal\}]) > 0.
  \eequation
  Reducing $\epsilon$, if necessary, one has $\dist(-\nabla f(x), \text{range}[\nabla e(x)] + \text{pos}[\{\nabla c_i(x) : i\in\widehat\Acal\}]) \geq \tau$ for the given $\widehat\Acal$ and all $x$ satisfying $\|x - x_*\|_2 \leq \epsilon$. The justification for this claim follows from standard continuity arguments; see \cite[Appendix 1]{OberWrig2006}.
  
  Now consider any such $x$ such that $\Acal_{{\rm QP}}(d_{x,\theta,\nu}; x) = \widehat\Acal$.  Let $(d_x, r_x, s_x, t_x)$ and $(\alpha_x, \beta_x)$ denote the solutions to~\eqref{eq.sub} and \eqref{eq.sub_dual}, respectively.  For any index $i \notin \widehat\Acal$, the KKT conditions for \eqref{eq.sub}--\eqref{eq.sub_dual} and equation \eqref{eq.active_QP} together imply that
  \bequation\label{eq.lemma4.2_proof_comp}
    c_i(x) + \nabla c_i(x)^Td_x < 0,\ \ [r_x]_i=0,\ \ \text{and}\ \  [\beta_x]_i=0.
  \eequation
  One can now analyze the objective of the dual subproblem~\eqref{eq.sub_dual} by separating it into two parts. First, using property \eqref{eq.lemma4.2_proof_comp}, one can observe that
  \begin{align*}
      &\ \tfrac{1}{2\theta}\|\nabla f(x) + \nabla e(x) \alpha_x + \nabla c(x) \beta_x\|_2^2 \\
      \geq&\ \tfrac{1}{2\theta} \min_{(\alpha,\beta) \in \R{p} \times \R{q}_{\geq0}} \|\nabla f(x) + \nabla e(x) \alpha + \nabla c(x) \beta\|_2^2\\
      =&\ \tfrac{1}{2\theta} \min_{(\alpha,\beta) \in \R{p} \times \R{q}_{\geq0}} \left\|\nabla f(x) + \nabla e(x) \alpha + \sum_{i\in\widehat\Acal}\nabla c_i(x) \beta_i \right\|_2^2 \\
      =&\ \tfrac{1}{2\theta} \dist(-\nabla f(x), \text{range}[\nabla e(x)] + \text{pos}[\{\nabla c_i(x) : i\in\widehat\Acal\}])^2 \geq \tfrac{1}{2\theta} \tau^2.
  \end{align*}
  From $\beta_x \geq 0$ and \eqref{eq.sub_dual}, one also has that $e(x)^T\alpha_x + c(x)^T\beta_x \leq e(x)^T\alpha_x + (c(x))_+^T\beta_x \leq \nu \|e(x)\|_1 + \nu \|(c(x))_+\|_1$. Combining these inequalities into the dual objective \eqref{eq.sub_dual}, one has
  \begin{multline}\label{eq.lara_loves_olives_1}
        e(x)^T\alpha_x + c(x)^T\beta_x - \tfrac{1}{2\theta}\|\nabla f(x) + \nabla e(x) \alpha_x + \nabla c(x) \beta_x\|_2^2 \\ \leq \nu||e(x)||_1 + \nu||(c(x))_+||_1 -\tfrac{1}{2\theta}\tau^2.
  \end{multline}
  Finally, one can decrease $\epsilon$ further, if necessary, to have that 
  \bequation\label{eq.lara_loves_olives_2}
    \nu\|e(x)\|_1 + \nu\|(c(x))_+\|_1 - \tfrac{1}{2\theta}\tau^2 < - \tfrac{\|x - x_*\|_2^\zeta}{4}\tau^2.
  \eequation
  This choice is possible since $e(x)$ and $(c(x))_+$ are both $\Ocal(\|x - x_*\|_2)$ and under the conditions of the theorem one has that $\theta \leq \|x - x_*\|_2^{-\zeta}$. Since $\zeta \in (0,1)$, the bounds \eqref{eq.lara_loves_olives_1}--\eqref{eq.lara_loves_olives_2} stand in direct contradiction to our earlier result, which established that the optimal value of the dual has
  \bequationNN
    -L(\|x-x_*\|_2^2+\|x-x_*\|_2) \leq  e(x)^T\alpha_x + c(x)^T\beta_x - \tfrac{1}{2\theta}\|\nabla f(x) + \nabla e(x) \alpha_x + \nabla c(x) \beta_x\|_2^2.
  \eequationNN
  This inconsistency forces the conclusion that $\tau$ defined in equation \eqref{eq.def_tau} cannot be positive, and in fact must be 0.  Thus, the inclusion in \eqref{eq.ran_pos_stat} must hold, as desired.
  
  Finally, consider part (c).  Since the LICQ holds, one has that $[z_*]_i>0$ for all $i\in \Acal(x_*) \setminus \Acal_0(x_*)$. Thus, for the inclusion \eqref{eq.ran_pos_stat} to hold, any index $i$ for which $[z_*]_i>0$ must satisfy $i\in\Acal_{{\rm QP}}(d_{x,\theta,\nu}; x)$; otherwise, \eqref{eq.ran_pos_stat} would be violated. Therefore, one can conclude the left-hand side inclusion $\Acal(x_*) \setminus \Acal_0(x_*) \subseteq \Acal_{{\rm QP}}(d_{x,\theta,\nu}; x)$. The final conclusion follows directly from the fact that strict complementarity implies that $\Acal_0(x_*)=\emptyset$, meaning there are no active inequality constraints at the solution for which the corresponding multiplier is zero. The proof is complete.
  \qed
\eproof

\bibliographystyle{plain}
\bibliography{references}

\end{document}